\documentclass[11pt,letterpaper]{article}
\usepackage[left=0.88in,right=0.88in,top=0.80in,bottom=0.74in,footskip=24pt]{geometry}
\usepackage[T1]{fontenc}
\usepackage[utf8]{inputenc}
\IfFileExists{newtxtext.sty}{\usepackage{newtxtext}}{\usepackage{times}}
\usepackage{amsmath,amssymb,amsthm,mathtools}
\IfFileExists{newtxmath.sty}{\usepackage{newtxmath}}{\usepackage{txfonts}}
\usepackage[scaled=0.92]{helvet}
\DeclareMathAlphabet{\mathsf}{OT1}{phv}{m}{n}
\SetMathAlphabet{\mathsf}{bold}{OT1}{phv}{b}{n}
\usepackage{bm,mathrsfs}
\DeclareFontFamily{U}{rsfs}{}
\DeclareFontShape{U}{rsfs}{m}{n}{<-6>rsfs5<6-8>rsfs7<8->rsfs10}{}
\usepackage{booktabs,tabularx,array,longtable}
\usepackage{microtype,graphicx}
\usepackage[dvipsnames]{xcolor}
\usepackage{placeins,float,needspace}
\usepackage[shortlabels]{enumitem}
\usepackage{titlesec,etoolbox}
\usepackage[font=small,labelfont=bf,labelsep=period]{caption}
\usepackage[numbers,sort&compress]{natbib}
\usepackage{xurl}
\usepackage[hidelinks,breaklinks=true]{hyperref}
\hypersetup{pdftitle={Dynamics of Protected State Evolution},pdfauthor={Qinyou Wang},pdfsubject={Controlled evolution under versioned protection contracts},pdfkeywords={controlled dynamics, protected state evolution, intrinsic geometry, irreducible cooperation, state quotients, resource accounting, generation promotion},pdfdisplaydoctitle=true}
\numberwithin{equation}{section}
\titleformat{\section}{\normalfont\fontsize{12}{14}\selectfont\bfseries}{\thesection}{0.75em}{}
\titleformat{\subsection}{\normalfont\normalsize\bfseries}{\thesubsection}{0.75em}{}
\titleformat{\subsubsection}{\normalfont\normalsize\itshape}{\thesubsubsection}{0.75em}{}
\titlespacing*{\section}{0pt}{2.0ex plus .4ex minus .2ex}{0.8ex plus .2ex}
\titlespacing*{\subsection}{0pt}{1.6ex plus .3ex minus .2ex}{0.6ex plus .2ex}
\titlespacing*{\subsubsection}{0pt}{1.3ex plus .3ex minus .2ex}{0.5ex plus .2ex}
\titlespacing*{\paragraph}{0pt}{1.4ex plus .2ex}{0.75em}
\makeatletter
\renewcommand{\@maketitle}{%
  \begin{center}
    {\fontsize{19}{22}\selectfont\bfseries\@title\par}
    \vskip 12pt
    {\normalsize\@author\par}
  \end{center}\vskip 9pt}
\makeatother
\setlist{nosep,leftmargin=*}
\newcolumntype{Y}{>{\raggedright\arraybackslash}X}
\newcolumntype{L}[1]{>{\raggedright\arraybackslash}p{#1}}
\allowdisplaybreaks[1]
\newtheorem{theorem}{Theorem}[section]
\newtheorem{proposition}[theorem]{Proposition}
\newtheorem{lemma}[theorem]{Lemma}
\newtheorem{corollary}[theorem]{Corollary}

\theoremstyle{definition}
\newtheorem{definition}[theorem]{Definition}
\newtheorem{example}[theorem]{Example}
\theoremstyle{remark}
\newtheorem{remark}[theorem]{Remark}
\newcommand{\protectenv}[1]{\AtBeginEnvironment{#1}{\Needspace{3\baselineskip}}}
\forcsvlist{\protectenv}{theorem,proposition,lemma,corollary,counterexample,definition,example,remark}
\newcommand{\R}{\mathbb R}

\newcommand{\calX}{\mathcal X}
\newcommand{\calY}{\mathcal Y}
\newcommand{\calK}{\mathcal K}
\newcommand{\calE}{\mathcal E}
\newcommand{\calO}{\mathcal O}
\newcommand{\calP}{\mathcal P}

\newcommand{\calA}{\mathcal A}

\newcommand{\cI}{\mathcal I}

\newcommand{\op}{\mathrm{op}}
\newcommand{\eff}{\mathrm{eff}}
\newcommand{\F}{\mathrm F}

\DeclareMathOperator{\Ran}{Ran}
\DeclareMathOperator{\Ker}{ker}

\DeclareMathOperator{\diag}{diag}

\DeclareMathOperator{\Tr}{Tr}
\DeclareMathOperator{\dist}{dist}
\DeclareMathOperator{\pr}{pr}
\DeclareMathOperator{\vecop}{vec}
\DeclareMathOperator{\KL}{KL}

\DeclareMathOperator{\softmax}{softmax}

\definecolor{FFBlue}{HTML}{4C78A8}
\definecolor{FFBlueLight}{HTML}{EAF1F8}
\definecolor{FFGold}{HTML}{D3A52F}
\definecolor{FFGoldLight}{HTML}{FBF4DF}
\definecolor{FFGray}{HTML}{6A7179}
\definecolor{FFRule}{HTML}{D3D7DB}
\title{Dynamics of Protected State Evolution}
\author{Qinyou Wang}
\date{}
\begin{document}
\maketitle
\begin{abstract}
States with identical protected observations can admit different future interventions. Each transition specifies its successor, response increment, required initial reserve and net consumption. Composing these records in execution order distinguishes protocols with equal net consumption and different initial reserve requirements. At a fixed state, a weighted Gram operator determines the linear response ellipsoid and a quotient identifies directions attainable only through cooperation between blocks. Nonlinear systems with identical complete affine data can have incomparable finite-budget response sets. In regular charts, the first correction to quadratic minimum cost is a covariant cubic form involving the third derivative of cost and the Hessians of protection and response. At singular gates, sequential feasibility depends on the peak reserve along the path. A history-weighted operator detects responses outside the old channel's range closure. Controlled congruence characterizes exact reductions of deterministic labelled dynamics, and a two-cut value identity separates the effect of paid verification from that of subsequent action selection.

On Qwen2.5-7B Base, two projected policies preserved 2,048 reference top tokens under a per-position KL ceiling of $2\times10^{-5}$. Over 200 arrivals they recovered 391 of 396 initial rejections. Each yielded over $1{,}200$ times the held-out cross-entropy reduction of historical-channel-only AdamW. Under an expanded 2,560-position contract, fixed continuation and neutral renewal began on the same protected logit fiber in distinct structural strata. Both admitted all 200 subsequent writes. Fixed continuation yielded $1.74$ times the held-out loss reduction per unit squared write cost. Renewal had a smaller KL excess on the first candidate tested at every arrival and used fewer full verification scans.
\end{abstract}
\section{Introduction}\label{sec:introduction}
A system that preserves an output can still change what it can do next. Its next intervention depends on the physical configuration and on the state left by earlier actions, including the controller's memory, the current protection contract, and the available reserve. Write $q_t$ for the physical and controller state, including its contract, and $b_t$ for the available reserve. Together they form the execution state. A protocol $\pi_t$ with required reserve $m(\pi_t)$ and net consumption $c(\pi_t)$ acts by
\[
  (q_t,b_t)\ \xrightarrow{\ \pi_t\ }\ (q_{t+1},\,b_t-c(\pi_t)),
  \qquad m(\pi_t)\preceq b_t .
\]
The state carries the reference values that the successor must preserve. At a promotion it also acquires new obligations. Execution order already matters for a scalar reserve. A check that requires the entire initial reserve and releases it can precede a positive debit. After the same debit the check is unaffordable. Section~\ref{m:sec:protocols} gives the composition law for arbitrary finite protocols.

The geometry at one state describes one part of this evolution. The linearized protection and response maps determine the affordable response ellipsoid of a synchronous chart. Curvature changes the cost of a finite response. The order in which blocks can move determines which joint corrections are executable. A committed action then supplies the parent for a new chart. Changes of coordinates transport the metric and the admissible actions together with the response maps.

In a synchronous linear chart the compensation equations are
\begin{equation}\label{eq:intro-compensation}
 Cx+By=0,\qquad Hx+\Psi y=D.
\end{equation}
The first equation preserves the protected response. The second specifies the target. In cost-orthogonal coordinates for the two source blocks, Theorem~\ref{m:thm:weighted} gives
\begin{equation}\label{eq:intro-gram}
 G=G_0+\Lambda W_{\mathcal E}^{-1}\Lambda^*,
 \qquad
 h_{\mathcal A_b}(z)=\sqrt{b\langle z,Gz\rangle}.
\end{equation}
On the compensable domain $\mathcal E$, $W_{\mathcal E}$ is the induced cost and $\Lambda$ the response including compensation. The image of the full cost ball has support function $h_{\mathcal A_b}$.

\begin{figure}[tbp]
\centering
\includegraphics[trim=4pt 4pt 4pt 4pt,clip]{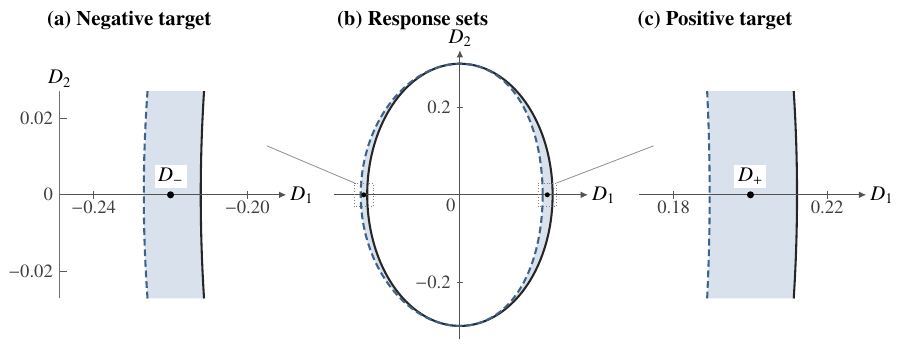}
\caption{Identical affine compensation data and incomparable finite response sets. In this extension of Theorem~\ref{m:prop:tuplematched}, protection is $x+y+ay^2=0$, response is $(D_1,D_2)=(x,z)$, and cost is $x^2+y^2+z^2$. Both systems admit the radius-$0.45$ chart and use squared budget $b=0.09$. The common tuple is $C=(1,0)$, $B=Q=1$, $H=I_2$, $\Psi=0$, $N=\operatorname{diag}(0,1)$, $L=(-1,0)^\top$, $W=2$, and $G=\operatorname{diag}(1/2,1)$. The center shows both sets. The flanking views magnify the boxed neighborhoods at identical scales. The solid curve is the common linear prediction ($a=0$), and the dashed curve gives $a=0.8$. Shading marks responses admitted by exactly one system. At $D_-=(-0.22,0)$ the two minimum costs are $0.0968$ and $0.0848$, respectively. At $D_+=(0.20,0)$ they are $0.0800$ and $0.1025$. These values follow from $E_a(D)=D_1^2+y_a(D_1)^2+D_2^2$, with $y_a(u)=-2u/(1+\sqrt{1-4au})$.}
\label{fig:nonlinear-capacity}
\end{figure}

Figure~\ref{fig:nonlinear-capacity} exhibits two systems with the same complete affine tuple. Each can reach a finite target that the other misses at the same budget.

A causal Transformer provides a concrete setting in which nonlinear readouts and controller memory interact. Its parameters fix the logits on a given input, and its AdamW moments enter the next update. Earlier work described observation-relative state fibers \citep{paper1} and regulation and acquisition operators at a fixed parent \citep{paper4}. Sections~\ref{t:sec:source}--\ref{exp:protected-adaptation} realize the state and reserve dynamics in this architecture and follow two trajectories across a change of contract.

\paragraph{From local response to closed-loop evolution.}
At a single parent, the Gram determines which linear responses fit a given budget. The quotient $V_{\rm joint}/V_{\rm sep}$ identifies those that lie outside the response space of separately feasible blocks (Section~\ref{m:sec:geometry}). Curvature and the order of execution enter at the next step. Regular charts that agree to first order can separate in their cubic cost, and the reserve needed to pass a singular gate depends on the route taken (Section~\ref{m:sec:finite}). Across parents, controlled congruence characterizes exact reductions of labelled dynamics. In a Hilbert-space chart, the history-weighted operator measures how far supported responses extend beyond the old range closure (Section~\ref{m:sec:closure}). For stochastic policies, the change in terminal value splits at the moment a paid observation becomes available to the decision rule (Theorem~\ref{e:thm:two-cut}). Table~\ref{tab:claims} summarizes the quantities that enter at each level.

\begin{table}[H]
\centering\small
\caption{Quantities carried from a single transition to an entire controlled trajectory.}\label{tab:claims}
\begin{tabularx}{\linewidth}{@{}p{0.18\linewidth}Y Y@{}}
\toprule
Object & Quantity determined & Structure governing the next comparison\\
\midrule
State and reserve & Successor, response increment and remaining reserve & Composition in execution order\\
Protected linear chart & Responses within budget and joint response classes & Physical decoder, curvature and execution order\\
Finite intervention & Native response and objective change at one parent & Continuation kernel and controller memory\\
Matched continuation & Terminal contrast and drift from a fixed anchor & Selection rule and accumulated certification cost\\
Complete policy & Expected terminal utility at a common initial reserve & Contract version and inherited obligations across generations\\
\bottomrule
\end{tabularx}
\end{table}

\paragraph{Mathematical connections and realization.}
For a scalar reserve, a record acts by partial translation between upper intervals, as in the bicyclic extensions of a totally ordered group's positive cone \citep{gutik}. Coupling the record to an intervention also tracks its response and the state from which the next action starts. Operator factorization and shorted operators describe the geometry at a fixed state \citep{douglas,anderson}. Perturbation theory relates the derivatives of the constraints to the cost of finite motion \citep{dontchev,bonnans,jackson}. Predictive representations and invariant observable spaces give related finite realization classes \citep{psr,koopman}. For stochastic transitions, probabilistic bisimulation provides a corresponding notion of state equivalence with rewards \citep{givan}.

Low-rank and latent adapters supply finite-dimensional decoders of changes to physical weights \citep{lora,erilora,metaucf,latentlora}. In our realization, a historical channel and an acquisition channel act on the same two internal projections. Replaying a write from a fixed parent isolates the historical channel's effect at a fixed acquisition write. Running the controllers on a common stream tests how those interventions accumulate. At the change of contract the controllers share every protected logit, providing a common starting observation for the comparison of inherited and renewed acquisition frames. Appendix~\ref{app:specializations} gives the fixed-parent affine specialization and its final-head interface.

\paragraph{Notation and physical metrics.}
All inner-product spaces are real. Adjoints use the reference inner products unless a physical metric is explicitly included. For a closed subspace $S\subseteq U$ with inclusion $\iota_S:S\hookrightarrow U$, the cost induced by a bounded coercive self-adjoint $M:U\to U$ is represented by the compression $M_S=\iota_S^*M\iota_S:S\to S$. In general, $M_S^{-1}\ne\iota_S^*M^{-1}\iota_S$. Whitening transports both the cost and the admissible action sets. The structural signature map $\varsigma:\mathcal S\to\Sigma$ partitions the state space into strata $\mathcal S_\sigma=\varsigma^{-1}(\{\sigma\})$. At time $t$ the signature is $\sigma_t=\varsigma(q_t)$. The contract version and memory relevant to future transitions belong to $q$. A reserve $b$ can be vector-valued. In a quadratic endpoint chart it is a scalar squared-norm budget. We write $G$ for the compensated response Gram and its matrix in a specified basis. In policy comparisons, $\mathcal G_t$ denotes the information $\sigma$-algebra after certification and selection. We write $\overline{\Ran T}$ for range closure.
% Complete-state relations, compensation geometry and historical response.
\section{Complete states and affordable protected protocols}\label{m:sec:protocols}
The state $q\in\mathcal S$ specifies the physical variables, controller memory, retained structure, protection references and active contract that determine the next actions and their successors. In a learning system it includes parameters and optimizer memory. A fixed-parent affine chart \cite{paper4} describes one intervention at $q$. A sequence of interventions also depends on the state each action leaves behind.

Write $x=(q,b)$ for the \emph{execution state}, where $b$ is the current reserve. The controller memory in $q$ includes any past charges used to choose an action or compute a readout.

\paragraph{Admissibility and information.}
At state $q$ a local action chart $(U_q,\Phi_q,c_q)$ offers each endpoint in $U_q$ at cost $c_q$. The active contract selects the protected endpoints through $\Phi_q$. For a pathwise contract the executor also specifies a family of paths and a resource rule for traversal. A continuous protected lift in that family gives an execution at its path cost. Decisions use the information available at their respective cuts. Fixed-parent trials use the same values of the physical and controller registers that define the local chart. Each branch carries the cost of its own trials.

\begin{definition}[Complete protected affordable relation]\label{m:def:global}
Let $V$ be a real vector lattice with cone $V_+$ and $\calO$ a real vector space. A protocol $\pi:q\leadsto q'$ admissible under the action and protection rules at $q$ has response increment $D(\pi)\in\calO$ and reserve record $(m(\pi),c(\pi))\in V_+\times V$ with $c(\pi)\preceq m(\pi)$. Each execution branch has the principal reserve domain $m(\pi)+V_+$ and leaves $b-c(\pi)$ from reserve $b$. Alternative executions contribute the union of their transition relations. For composable protocols the response obeys $D(\pi_2\circ\pi_1)=D(\pi_1)+D(\pi_2)$. Define
\[
\mathscr C(q,b)=\{(D(\pi),q',b-c(\pi)):\pi:q\leadsto q',\;m(\pi)\preceq b\},\qquad
\calA_b(q)=\pr_{\calO}\mathscr C(q,b).
\]
For trace-valued responses, trace concatenation supplies the composition law.
\end{definition}

\begin{proposition}[Exact record and relation composition]\label{m:prop:record}
The reserve records form an associative monoid under
\[
(m_2,c_2)\circ(m_1,c_1)=\big(m_1\vee(c_1+m_2),\ c_1+c_2\big),
\]
with identity $(0,0)$. For complete states, the relation of a two-stage protocol is the ordinary relational composite on the intermediate \emph{state and reserve}, and its response increments add. The record acts faithfully as the partial translation $T_{m,c}(b)=b-c$ on the domain $m+V_+$.
\end{proposition}
\begin{proof}
A two-stage execution requires $b\succeq m_1$ and $b-c_1\succeq m_2$, hence exactly $b\succeq m_1\vee(c_1+m_2)$. Its final reserve is $b-c_1-c_2$. Three successive stages require $b$ to dominate $m_1$, $c_1+m_2$, and $c_1+c_2+m_3$, which is independent of parenthesization. The composite threshold is nonnegative because it dominates $m_1\succeq0$, and $c_1+c_2\preceq c_1+m_2\preceq m_1\vee(c_1+m_2)$, so the composite is again a record. The unit laws follow from $c\preceq m$ and $m\succeq0$. A partial translation determines its domain's least element $m$, and evaluation at $m$ then determines $c$. The suffix starts from the same intermediate state and reserve.
\end{proof}
\begin{remark}[Budget semantics]\label{m:remark:2_3}
A single endpoint-write cost $\|u\|^2$, cumulative sums of squared microscopic increments, and the peak anchored deviation $\sup_t\|u_t-u_{\rm anchor}\|^2$ are different contracts. For example, the stages $0\to1\to2$ each have local squared peak one, but their common-anchor squared peak is four. A common-anchor peak record has zero net consumption. An irreversible write charge reduces the remaining reserve. Continuous vector-valued peaks require existence of the relevant order supremum (for example, finite-dimensional continuous compact paths, or an order-bounded family in a Dedekind-complete lattice).
\end{remark}

\subsection{Versioned protection and generation promotion}
\begin{definition}[Protection contract and nested promotion]\label{m:def:contract-promotion}
A generation carries a contract
\[
 \mathfrak C_g=\bigl(g,\mathcal J_g,
       \{(P_j,a_j,\mathcal T_j)\}_{j\in\mathcal J_g},
       \mathcal F_g,z_g\bigr),
\]
stored in the complete state. Each $P_j$ is a readout into a real vector space, $a_j$ is its fixed reference, and $\mathcal T_j$ is an allowed error set containing zero. The set $\mathcal F_g$ identifies physical registers whose committed values are $z_g$; these may include controller-memory registers. On a common semantic state domain define
\[
 \mathcal Q_g=\{q:P_j(q)-a_j\in\mathcal T_j\ (j\in\mathcal J_g),\quad
                       \operatorname{val}_{\mathcal F_g}(q)=z_g\}.
\]
This is the protection predicate evaluated on a state, without requiring its active version field to equal $g$.

A \emph{nested promotion} is an atomic transition $q^-\leadsto q^+$ from version $g$ to $g+1$, satisfying the following conditions. All old triples $(P_j,a_j,\mathcal T_j)$ and physical register identities are retained, with $\mathcal J_g\subseteq\mathcal J_{g+1}$, $\mathcal F_g\subseteq\mathcal F_{g+1}$ and $z_{g+1}|_{\mathcal F_g}=z_g$. New obligations use readouts and reference data available at promotion, with anchors fixed by a specified state or rule, and $q^+\in\mathcal Q_{g+1}$. The transition carries its reserve record $(m_{\rm prom},c_{\rm prom})$ under Definition~\ref{m:def:global}, including required formation, checking, storage and commit costs. Any staged implementation must preserve the old contract throughout its execution before the commit.

Throughout generation $g+1$, every admissible transition preserves $\operatorname{val}_{\mathcal F_{g+1}}(q)=z_{g+1}$, including indirect updates induced by the controller. A rejected promotion leaves $\mathfrak C_g$ active and retains the reserve record of the attempt.
\end{definition}

\begin{proposition}[Preservation under promotion and same-parent kernel nesting]\label{m:prop:contract-promotion}
For a nested promotion, $\mathcal Q_{g+1}\subseteq\mathcal Q_g$. Thus the commit and any subsequent protocol that enforces $\mathcal Q_{g+1}$ throughout also preserve the old obligations. Promotion and its continuation compose in that order under Proposition~\ref{m:prop:record}.
For a separate linearized comparison, evaluate the retained and enlarged \emph{exact-equality} obligations at one common parent $q_*$ satisfying both lists of exact equalities, in the same differentiable action chart with the same control and output Hilbert spaces, bounded response derivative $H$ and bounded coercive self-adjoint physical metric $M$. Let $C_g$ contain the old protection and frozen-register rows, $C_{\rm add}$ the added readout rows, and $E_{\rm fr}$ the newly frozen-register rows, all bounded. Then
\[
 C_{g+1}=\begin{bmatrix}C_g\\ C_{\rm add}\\ E_{\rm fr}\end{bmatrix},
 \qquad
 \ker C_{g+1}=\ker C_g\cap\ker C_{\rm add}\cap\ker E_{\rm fr}
                  \subseteq\ker C_g.
\]
Writing $\bar C_i=C_iM^{-1/2}$ and $\bar H=HM^{-1/2}$, the corresponding linearized protected response Grams obey
\[
 G_i=\bar H P_{\ker\bar C_i}\bar H^*,\qquad
 G_{g+1}\preceq G_g.
\]

\end{proposition}
\begin{proof}
Every condition defining $\mathcal Q_g$ occurs unchanged among those defining $\mathcal Q_{g+1}$, including the old frozen values. This proves predicate inclusion. The commit condition supplies the initial witness; applying inclusion at every protected intermediate state proves preservation along a subsequent enforcing protocol. 

In the common chart, an exact protected first-order action must satisfy every retained, added and newly frozen equality. The kernel of a stacked bounded operator is the intersection of its row-block kernels. The same inclusion holds after the single cost-whitening map $M^{-1/2}$. These kernels are closed; for nested closed subspaces $K'\subseteq K$, $P_K-P_{K'}=P_{K\cap(K')^\perp}\succeq0$. Multiplication by $\bar H$ and its adjoint proves the Gram order. A tolerance set $\mathcal T_j$ with nonzero width gives a feasibility inequality for the corresponding readout.
\end{proof}

\subsection{The finite-budget primitive and its attainment boundary}
At a fixed state $q$, let $u\mapsto \Upsilon_q(u)$ be an action chart with endpoint domain $U_q\ni0$ and physical cost $c_q:U_q\to[0,\infty)$ satisfying $c_q(0)=0$. Put
$\Phi_q(u)=(p_q(u),r_q(u))$ with $\Phi_q(0)=0$ and define
\[
\begin{aligned}
\mathcal A_b^{\Upsilon}(q)&=\{D:\exists\,u\in U_q,\ p_q(u)=0,\ r_q(u)=D,\ c_q(u)\le b\},\\
E_q(D)&=\inf\{c_q(u):u\in U_q,\ p_q(u)=0,\ r_q(u)=D\},\qquad\inf\varnothing=+\infty.
\end{aligned}
\]
Here $p_q$ measures the change in the protected readout. We suppress the chart superscript when the chart is fixed.
\begin{proposition}[Attainability at the cost boundary]\label{m:prop:attain}
For any $b>0$,
\[
\{D:E_q(D)<b\}\subseteq\mathcal A_b^{\Upsilon}(q)\subseteq\{D:E_q(D)\le b\}.
\]
If every target whose infimum is at most $b$ has a minimizer in $U_q$, then the second inclusion is equality.
\end{proposition}
\begin{proof}
If $E_q(D)<b$, the definition of infimum supplies a witness of cost below $b$. Any witness of cost at most $b$ proves $E_q(D)\le b$. If the infimum is attained, its minimizer is itself an affordable witness. If the possible costs of $D$ are exactly $1+1/n$ for $n\ge1$, then $E_q(D)=1$ and $D\notin\mathcal A_1^{\Upsilon}(q)$.
\end{proof}

\subsection{Regulation, acquisition, and admissible action words}
\begin{definition}[Structural signature and the two action categories]\label{m:def:frtfat}
Let $\mathcal S$ be the complete protected-state space of Definition~\ref{m:def:global}.
A \emph{structural-signature map} $\varsigma:\mathcal S\to\Sigma$ identifies existing writable components, their physical identities, admissible write interfaces and retained structural objects. Write $\sigma_t=\varsigma(q_t)$ and $\mathcal S_\sigma=\varsigma^{-1}(\{\sigma\})$ for the realized signature and its state stratum.
A \emph{regulation action} preserves the signature, protection and action-cost contract. It may change existing coefficients, controller memory or future response operators.
An \emph{acquisition action} changes the signature by creating, replacing or renewing writable structure. Replacement can change the available directions at fixed dimension.
A \emph{joint protocol} is a finite word over the two categories, with each label available at its intermediate state. Proposition~\ref{m:prop:record} composes its resource records in order. Regulation followed by acquisition and acquisition followed by regulation can have different successors.
\end{definition}
\begin{proposition}[Fixed-signature regulation can change the response Gram]\label{m:prop:frtdistinct}
A regulation at fixed structural signature can change the Gram of a subsequent action.
\end{proposition}
\begin{proof}
Use complete states $(h,z)$ with a stored gain $h\in[0,1]$, a response accumulator $z\in\R$, a fixed structural signature, and protected present readout $P(h,z)=0$. At any state with $h=0$, define the regulation $(0,z)\mapsto(1,z)$ with reserve record $(1,1)$. An existing probe $y$ has successor $(h,z)\mapsto(h,z+hy)$, response increment $hy$, and record $(y^2,y^2)$. Response increments telescope in $z$ and all protection conditions are satisfied. The current probe Gram is $h^2$. Starting at $(0,0)$, regulation followed by one probe has $|D|\le\sqrt{b-1}$ for $b>1$, although the original probe alone has zero response.
\end{proof}
The image of a protocol family grows under inclusion. In particular, a mixed family that admits each pure protocol at the same cost contains both pure response images at every reserve. Theorem~\ref{m:thm:congruence} gives the corresponding criterion for reducing the states on which these protocols act.

\Needspace{6\baselineskip}
\section{Weighted compensation and algebraic cooperation}\label{m:sec:geometry}
\begin{theorem}[Weighted orthogonal compensation]\label{m:thm:weighted}
Let $\calX,\calY,\calP,\calO$ be real Hilbert spaces; let $C,B,H,\Psi$ be bounded with the domains and codomains of the block equations $Cx+By=0$, $Hx+\Psi y=D$. Assume $\Ran C$ is closed. Equip the action space with the block metric $M=M_X\oplus M_Y$, with bounded coercive self-adjoint blocks. Put $K_0=\ker C$, $\calE=B^{-1}(\Ran C)$ and let $\iota_0,\iota_E$ be their inclusions. Define
\[
 C_{M_X}^{\dagger}=M_X^{-1/2}(CM_X^{-1/2})^{\dagger},\quad
 T=C_{M_X}^{\dagger}B\iota_E,\quad
 \Lambda=\Psi\iota_E-HT,
\]
\[
 M_0=\iota_0^*M_X\iota_0,\qquad
 W_{\calE}=\iota_E^*M_Y\iota_E+T^*M_XT.
\]
Then $\calE$ is closed, $M_0$ and $W_{\calE}$ have bounded inverses, and
\[
\ker[C\ B]=(K_0\oplus\{0\})\oplus_M\{(-Ty,\iota_Ey):y\in\calE\}.
\]
In the cost Hilbert space, the complete protected response Gram is
\begin{equation}\label{m:eq:gram}
\boxed{G=G_0+\Lambda W_{\calE}^{-1}\Lambda^*,\qquad
G_0=H\iota_0 M_0^{-1}\iota_0^*H^*.}
\end{equation}
The additional cost-whitened channel is $\widetilde{\cal L}=\Lambda W_{\calE}^{-1/2}$. In the Euclidean finite-dimensional case, an isometric coordinate map $Q:\R^r\to\calY$ onto $\calE$ gives the fixed-chart formulas of \cite[Theorem~5.3]{paper4}
\[
 L=(\Psi-HC^{\dagger}B)Q,\qquad
 W=I+Q^*B^*(C^{\dagger})^*C^{\dagger}BQ,\qquad
 G=NN^*+LW^{-1}L^*,\quad N=HP_{K_0}.
\]
\end{theorem}
\begin{proof}
The map $M_X^{-1/2}$ is onto, so $CM_X^{-1/2}$ has closed range and bounded Moore--Penrose inverse. Solving $Cx=d$ with minimum $M_X$-norm is equivalent to solving $(CM_X^{-1/2})w=d$ with minimum reference norm and taking $x=M_X^{-1/2}w$. Hence $CC_{M_X}^{\dagger}=P_{\Ran C}$ and $C_{M_X}^{\dagger}d$ is $M_X$-orthogonal to $K_0$ for $d\in\Ran C$. The inverse image $\calE$ of the closed range is closed. Every protected pair has the unique representation $(\iota_0z-Ty,\iota_Ey)$, and
\[
 \|(\iota_0z-Ty,\iota_Ey)\|_M^2
 =\langle z,M_0z\rangle+\langle y,W_{\calE}y\rangle,
 \qquad \bar H(\iota_0z-Ty,\iota_Ey)=H\iota_0z+\Lambda y.
\]
Coercivity of the compressions follows by applying coercivity to the included vector; $W_{\calE}$ adds a nonnegative form. Whiten the two orthogonal cost coordinates. The resulting response row operator is
$[H\iota_0M_0^{-1/2}\ \Lambda W_{\calE}^{-1/2}]$, whose product with its adjoint is \eqref{m:eq:gram}. Euclidean coordinates give the displayed $L,W,N$. The block elimination takes place in the cost inner product, as in range-factorization geometry \cite{douglas}.
\end{proof}

\begin{theorem}[Exact budget image and support function]\label{m:thm:support}
In the synchronous chart $U_q=\calK=\ker[C\ B]$ with cost $\|u\|_M^2$, let $\bar H=[H\ \Psi]$ and $\calA_b=\{\bar Hu:u\in\calK,\|u\|_M^2\le b\}$. For $b\ge0$, $\calA_b$ is convex, weakly compact, and
\begin{equation}\label{m:eq:support}
\boxed{h_{\calA_b}(z)=\sqrt{b\langle z,G z\rangle}.}
\end{equation}
For each fixed $b>0$, within one common linear output/cost metric, Gram equality is equivalent to equality of budget images, and Loewner Gram order is equivalent to inclusion. At $b=0$ every budget image is $\{0\}$. For a finite-dimensional surjective output restriction, minimum cost for target $D$ is $\langle D,G^{-1}D\rangle$.
\end{theorem}
\begin{proof}
The bounded weighted Hilbert ball is weakly compact and its bounded linear image remains so. Riesz representation gives the support function by maximizing the linear functional $\langle z,\bar H u\rangle$. Support functions characterize closed convex sets; polarization converts equal quadratic forms into operator equality. The minimum-norm preimage formula gives the cost in the positive-definite case.
\end{proof}

For $C=1$, $B=\beta$, $H=0$ and $\Psi=1$ in Euclidean coordinates, the compensated column is always $L=1$, while the response Gram is $G=(1+\beta^2)^{-1}$.

\Needspace{6\baselineskip}
\subsection{Joint conditioning and the minimum-cost lift}
\begin{lemma}[Joint conditioning identity]\label{m:lem:conditioning}
Let $J=[\bar C;\bar H]:\R^d\to\R^{p+o}$ and $G=\bar H(I-\bar C^\dagger\bar C)\bar H^*$ in Euclidean coordinates. Then $J$ is surjective if and only if $\bar C$ is surjective and $G\succ0$. In that case the minimum \emph{squared} norm of $Jz=(a,b)$ equals
\[
\|\bar C^\dagger a\|^2+
\langle b-\bar H\bar C^\dagger a, G^{-1}(b-\bar H\bar C^\dagger a)\rangle,
\]
and in particular $\|J^\dagger(0,D)\|^2=\langle D,G^{-1}D\rangle$.
\end{lemma}
\begin{proof}
Surjectivity of $J$ implies that of $\bar C$. Given $\bar Cz=a$, write $z=\bar C^\dagger a+w$ where $w\in\ker\bar C$; the components are orthogonal. The second equation requires $\bar Hw=b-\bar H\bar C^\dagger a$, which is solvable for all targets precisely when $\bar H|_{\ker\bar C}$ is onto. Its Gram is $G$ and its minimum squared preimage norm is the displayed quadratic form. Conversely, onto $\bar C$ and positive-definite $G$ make both equations solvable for every $(a,b)$.
\end{proof}

\begin{theorem}[Metric-free residual quotient]\label{m:thm:alg}
For real vector spaces and finitely many blocks, form the algebraic quotient $\mathfrak D=\calP/\Ran C$ and put $\delta_i(y_i)=[B_i y_i]$, $E_{\rm joint}=\ker(\sum_i\delta_i)$, $E_{\rm sep}=\bigoplus_i\ker\delta_i$, and
\[
 Z=\ker\left(\bigoplus_i\Ran\delta_i\xrightarrow{\Sigma}\mathfrak D\right).
\]
Let $K_{\rm joint}=\{(x,y):Cx+\sum_iB_i y_i=0\}$, $K_{\rm sep}=\{(x,y)\in K_{\rm joint}:y\in E_{\rm sep}\}$, $\bar H(x,y)=Hx+\sum_i\Psi_i y_i$, and $V_\square=\bar H(K_\square)$. There is an exact sequence
\[
0\longrightarrow E_{\rm sep}\longrightarrow E_{\rm joint}
\xrightarrow{\omega}Z\longrightarrow0,\qquad \omega(y)=(\delta_i y_i)_i.
\]
For $y\in E_{\rm joint}$, any lift $Cx_y=-\sum_iB_i y_i$ defines the same quotient response $[Hx_y+\sum_i\Psi_i y_i]\in\calO/H(\ker C)$. Furthermore
\[
 \varpi:Z\to V_{\rm joint}/V_{\rm sep},\quad
 \varpi(\omega(y))=[\bar H(x_y,y)]
\]
is a well-defined surjection and $V_{\rm joint}/V_{\rm sep}\simeq Z/\ker\varpi$.
\end{theorem}
\begin{proof}
The residual map $\omega$ has kernel $E_{\rm sep}$ and is onto $Z$ by choosing a preimage in each block of any residual tuple summing to zero. Two old-channel lifts differ in $\ker C$, giving the first quotient's independence. If $\omega(y)=\omega(y')$, then $y-y'\in E_{\rm sep}$ and the difference of any two protected lifts lies in $K_{\rm sep}$. Thus $\varpi$ is well defined. Every joint output has such a lift, so $\varpi$ is onto, and the first-isomorphism theorem gives the last identity.
\end{proof}
The quotient measures directions outside $V_{\rm sep}$; a cost advantage within $V_{\rm sep}$ is measured by the Gram. For example, impose $y_1-y_2=0$ with no protection constraint on $x$, response $x+y_1+y_2$ and Euclidean squared cost. Separate compatibility forces $y_1=y_2=0$, while joint compatibility allows $y_1=y_2$. Both response spaces are $\mathbb R$. The target $1$ costs $1$ separately and $1/3$ jointly.

\subsection{Irreducible cooperation and response-versus-control filtering}
The algebraic residual quotient records which protection deficits cancel across modules. The response map then selects the directions that deliver task output and may annihilate the entire cooperative subspace. For the following Hilbert projection identity, require the two protected subspaces to be closed (as when $\Ran C$ is closed and the block maps are bounded). First whiten the fixed physical cost and transport both subspaces and $\bar H$ by that same isometry; the following symbols denote these transported objects. Let $K_{\rm sep}\subseteq K_{\rm joint}$ be the separately compatible and jointly compatible protected control spaces. The orthogonal complement
\[
K_{\rm coop}=K_{\rm joint}\cap K_{\rm sep}^{\perp}
\]
yields the exact identity
\[
G_{\rm joint}-G_{\rm sep}
=\bar H P_{K_{\rm coop}}\bar H^*\succeq0.
\]
The increment measures the part of $K_{\rm coop}$ visible through $\bar H$.
\begin{theorem}[Persistence of an irreducible coalition]\label{m:thm:coalition}
Fix $m\ge2$. Let $B:\R^m\to\R^{m-1}$ have columns $e_1,\ldots,e_{m-1},-\sum_{j<m}e_j$. On a neighborhood of zero let
\[
p(y)=By+f(y),\quad r(y)=y_1+g(y),
\]
where $f,g$ are $C^2$, vanish at zero, and have vanishing first derivatives there. Then on a sufficiently small common neighborhood, every proper subset of structural coordinates, with all coordinates outside that subset fixed at zero, has only the trivial exactly protected action, while the full $m$-block system reaches every sufficiently small scalar target $D$ with a locally unique action and squared cost
\[
y(D)=D\mathbf 1+O(D^2),\qquad E(D)=mD^2+O(|D|^3).
\]
\end{theorem}
\begin{proof}
Every proper subset of the columns of $B$ is linearly independent. On each corresponding coordinate subspace $S$, its minimum singular value $c_S$ is positive. Since $f(y)=O(\|y\|^2)$, an exactly protected $y\in S$ satisfies $c_S\|y\|\le C\|y\|^2$, hence is zero in a sufficiently small ball. Finitely many proper coalitions permit one common radius. In the full system, the derivative of the joint protection--response map is $[B;e_1^*]$, which is square and invertible because $\ker B=\operatorname{span}\{\mathbf1\}$ and $e_1^*\mathbf1=1$. The inverse-function theorem gives the unique local action delivering $(0,D)$, with derivative $y'(0)=\mathbf1$. Squared norm gives the cost expansion.
\end{proof}
\begin{proposition}[Non-submodular structural opportunity]\label{m:prop:nonsubmodular}
For the linear family above with $f=g=0$, let $\nu(S)$ denote the dimension of the protected obligation-response space when only modules in $S\subseteq\{1,\ldots,m\}$ are writable. Then $\nu(S)=0$ for every proper subset $S$, while $\nu(\{1,\ldots,m\})=1$. This family violates submodularity and has its first nonzero coalition at order $m$.
\end{proposition}
\begin{proof}
For proper $S$, the protected kernel is zero since the corresponding columns are linearly independent. For all $m$ blocks the protected kernel is $\operatorname{span}\{\mathbf1\}$, on which $r(y)=y_1$ is onto. Put $A=\{1,\ldots,m\}\setminus\{1\}$ and $B=\{1,\ldots,m\}\setminus\{2\}$. Then $\nu(A)+\nu(B)=0<1=\nu(A\cup B)+\nu(A\cap B)$ for every $m\ge2$.
\end{proof}
\Needspace{6\baselineskip}
\subsection{Sequential and synchronous execution}
\begin{figure}[tbp]
\centering
\includegraphics[trim=4pt 4pt 4pt 4pt,clip]{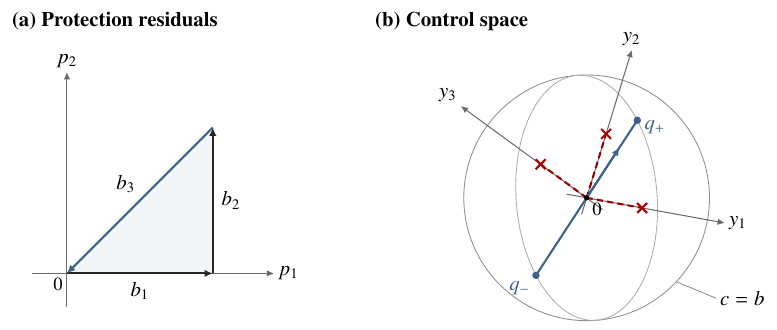}
\caption{Three-block cooperation and sequential obstruction. In (a) the protection columns $b_1=(1,0)^\top$, $b_2=(0,1)^\top$ and $b_3=(-1,-1)^\top$ close a triangle. Every proper subset is linearly independent. The protected response dimension is zero on each proper subset and one on the triple (Theorem~\ref{m:thm:coalition}). Adding block 3 to $\{1\}$ therefore has zero marginal gain, while adding it to $\{1,2\}$ has gain one (Proposition~\ref{m:prop:nonsubmodular}). In (b) protection is $p(y)=(y_1-y_3,y_2-y_3)=0$, response is $D=y_1$, and cost is $c=\|y\|^2$. Crosses mark single-coordinate first steps with nonzero residual. The blue diameter is the intersection of $\ker B$ with the budget ball in an orthographic view. Its endpoints are $q_\pm=\pm\sqrt{b/3}\,(1,1,1)$. Synchronous execution follows $y(s)=st(1,1,1)$, $0\le s\le1$. The path preserves protection and reaches $D=t$ at squared cost $3t^2$ (Proposition~\ref{m:prop:sequential}).}
\label{fig:irreducible-cooperation}
\end{figure}

\begin{proposition}[Sequential obstruction]\label{m:prop:sequential}
Assume fixed linear protection $Cx+\sum_iB_i y_i=0$, zero initial new blocks, and an executor changing at most one $y_i$ per stage while possibly adjusting $x$ at the same stage. Every finite exactly protected sequence ends in $E_{\rm sep}=\bigoplus_i\ker\delta_i$. Suppose that at every intermediate protected state each paired straight segment $(x,y_i)\mapsto(x+t\Delta x,y_i+t\Delta y_i)$, $0\le t\le1$, with $C\Delta x+B_i\Delta y_i=0$ is path-admissible, with unrestricted budget, sign and domain. Then the sequentially reachable new-block set equals $E_{\rm sep}$. Under synchronous admissibility, the corresponding algebraic new-block space is $E_{\rm joint}$.
\end{proposition}
\begin{proof}
For a stage touching block $i$, protection before and after gives $C\Delta x+B_i\Delta y_i=0$, hence $\delta_i(\Delta y_i)=0$. Summing from zero proves the necessary inclusion. For the converse, take arbitrary $y\in E_{\rm sep}$ and choose $x_i$ with $Cx_i=-B_i y_i$. Starting from the protected old state, successively follow the paired segments with increments $(x_i,y_i)$. Each segment stays exactly protected by linearity and is path-admissible by hypothesis. Their final new coordinate is $y$.
\end{proof}
Figure~\ref{fig:irreducible-cooperation} illustrates the irreducible three-block construction and its fixed-constraint execution boundary.

\section{Finite nonlinear compensation and singular execution}\label{m:sec:finite}
\subsection{Exact finite repair from a joint regularity certificate}

\begin{theorem}[Finite protected repair with explicit budget bounds]\label{m:thm:finite}
Let $\Phi=(p,r):\overline B_{r_0}(0)\subseteq\R^d\to\R^{p+o}$ extend to a $C^{1,1}$ map on a neighborhood, with $\Phi(0)=0$ and a surjective $J=D\Phi(0)$. Use the closed action ball with cost $c(u)=\|u\|^2$. Suppose $D\Phi$ is Lipschitz with constant $M$. Set
\[
\sigma_*=\|J^\dagger\|^{-1}>0,\quad
r_* =\min\{r_0,\sigma_*/(2M)\},\quad
 d_* =\sigma_* r_*/4,\quad k_*=2M/\sigma_*^3,
\]
where $r_*=r_0$ for $M=0$. For every $\|D\|\le d_*$ there exists a continuous local selection $u(D)\in\operatorname{Ran}J^*$ with
\[
\Phi(u(D))=(0,D),\qquad
\|u(D)-J^\dagger(0,D)\|\le k_*\|D\|^2,\qquad
\|u(D)\|\le r_*/2.
\]
The curve $t\mapsto u(tD)$ preserves protection throughout. The minimum squared endpoint cost $E(D)$ on the ball exists and, with $e(D)=\sqrt{\langle D,G^{-1}D\rangle}$, satisfies
\[
\big[e(D)-k_*\|D\|^2\big]_+^2\le E(D)\le
\big[e(D)+k_*\|D\|^2\big]^2.
\]
In particular $E(D)=\langle D,G^{-1}D\rangle+O(\|D\|^3)$.
\end{theorem}
\begin{proof}
For $D=0$, choose $u=0$, which is the unique zero-cost minimizer. For $D\ne0$, write $N(u)=\Phi(u)-Ju$. Then $\|N(u)\|\le M\|u\|^2/2$ and $\|DN(u)\|\le M\|u\|$. For $\|D\|\le d_*$ set $a=2\|D\|/\sigma_*\le r_*/2$ and consider on $\overline B_a\cap\operatorname{Ran}J^*$ the self-map
\[
F_D(u)=J^\dagger\big((0,D)-N(u)\big).
\]
Its Lipschitz constant is at most $Ma/\sigma_*\le1/4$ and its norm is at most $\|D\|/\sigma_*+Ma^2/(2\sigma_*)\le5a/8<a$. The contraction argument \cite{dontchev} gives a fixed point. Applying $J$ gives the exact target because $JJ^\dagger=I$. The remainder estimate gives $\|u-J^\dagger(0,D)\|\le Ma^2/(2\sigma_*)=k_*\|D\|^2$. The contraction estimates are uniform for targets in the same small target ball, yielding continuity and an exact-protection selected path. A closed feasible subset of the compact physical action ball attains its minimum. The candidate above gives an upper bound for its norm. Any minimizer $u_*$ has $\|u_*\|\le a$, and
\[
\|u_*\|\ge\|J^\dagger Ju_*\|
\ge\|J^\dagger(0,D)\|-\|J^\dagger N(u_*)\|
\ge e(D)-k_*\|D\|^2.
\]
The constructed candidate gives the corresponding upper bound. Squaring yields the displayed bounds and the cubic-order cost remainder.
\end{proof}
\begin{corollary}[Budget inner/outer tests]\label{m:cor:budgettests}
Within the target ball above,
\[
e(D)+k_*\|D\|^2\le\sqrt b\ \Longrightarrow\ D\in\mathcal A_b,
\qquad
e(D)-k_*\|D\|^2>\sqrt b\ \Longrightarrow\ D\notin\mathcal A_b.
\]
Figure~\ref{fig:nonlinear-capacity} gives two finite nonlinear budget images with the same tangent ellipsoid.
\end{corollary}
\begin{proof}
The first inequality and the attained upper cost bound give a feasible point of cost at most $b$; the second makes every feasible action too expensive by the lower cost bound.
\end{proof}
For $p(x,y)=x$ and $r_{\kappa}(x,y)=y-\kappa y^2$ the tangent Gram equals one at every $\kappa>0$, whereas the largest positive response equals $1/(4\kappa)$. Uniform aperture estimates therefore require control beyond the tangent Gram.

\begin{theorem}[Uniform finite-response aperture]\label{m:thm:uniformaperture}
Let $\Phi_n=(p_n,r_n):U_n\to P_n\oplus\calO$ be a family of $C^1$ maps between real Hilbert spaces on open neighborhoods of the closed radius-$r_0$ balls. Suppose $\Phi_n(0)=0$, $J_n=D\Phi_n(0)$ is bounded and onto for every $n$, each closed radius-$r_0$ ball is the action domain with cost $\|u\|^2$, and a common $K\ge0$ satisfies
\[
 \|J_n^\dagger(D\Phi_n(u)-J_n)\|_{\op}\le K\|u\|\qquad(\|u\|\le r_0).
\]
For fixed $d\ne0$, put $a_n(d)=\|J_n^\dagger(0,d)\|>0$ and
\[
 \alpha_n(b;d)=\sup\{t\ge0:\exists u\in\overline B_{\sqrt b},\ \Phi_n(u)=(0,td)\}.
\]
For $\rho=\sqrt b\in(0,r_0]$ with $K\rho<1$,
\begin{equation}\label{m:eq:aperture}
\boxed{\frac{\rho-K\rho^2/2}{a_n(d)}\le\alpha_n(b;d)\le\frac{\rho+K\rho^2/2}{a_n(d)}.}
\end{equation}
The protected linear Gram has bounded inverse, and
\[
 a_n(d)^2=\langle d,G_n^{-1}d\rangle,\qquad
 \sup_n\left|\frac{\alpha_n(b;d)a_n(d)}{\sqrt b}-1\right|\le\frac K2\sqrt b.
\]
The dependence on $n$ enters the leading aperture through $a_n(d)^{-1}$.
\end{theorem}
\begin{proof}
Surjectivity gives a bounded Moore--Penrose inverse and a closed row space $\Ran J_n^*$. Fix $n$ and write $R(u)=\Phi_n(u)-J_nu$, $F(u)=J_n^\dagger R(u)$ and $v=J_n^\dagger(0,d)$. The fundamental theorem of calculus on line segments gives
\[
 F(0)=0,\qquad \|F(u)\|\le K\|u\|^2/2,
 \qquad \|F(u)-F(w)\|\le K\rho\|u-w\|
\]
for $u,w$ in the closed $\rho$-ball. If $\Phi_n(u)=(0,td)$ and $\|u\|\le \rho$, then
$tv=J_n^\dagger J_nu+F(u)$; since $J_n^\dagger J_n$ is an orthogonal projection, $ta_n(d)\le \rho+K\rho^2/2$. This proves the upper bound even when the supremum is not attained.
For $t_-=(\rho-K\rho^2/2)/a_n(d)>0$, the map
\[
 T(u)=t_-v-F(u)
\]
takes the complete space $\overline B_\rho\cap\Ran J_n^*$ into itself, has Lipschitz constant $K\rho<1$, and therefore has a fixed point. Applying $J_n$ and using $J_nJ_n^\dagger=I$ gives $\Phi_n(u)=(0,t_-d)$. The same construction for $t\in[0,t_-]$ uses the common ball and satisfies
$\|u(t)-u(t')\|\le a_n(d)|t-t'|/(1-K\rho)$, so the selected exact-protection curve is continuous. The two bounds are uniform in $n$. The minimum-norm identity follows from the orthogonal solution of the joint linear equations, as in Lemma~\ref{m:lem:conditioning}.
\end{proof}

\subsection{Cubic cost and coordinate covariance}
The protected Gram determines the quadratic cost of a prescribed response. At the next order, the curvature of the constraints enters together with the third derivative of the physical cost.

\begin{theorem}[Covariant cubic minimum-cost invariant]\label{m:thm:cubic}
Let $\Phi=(p,r)$ and $c$ be $C^3$ near $0$ in a finite-dimensional action space, with $\Phi(0)=0$, $J=D\Phi(0)=[C;H]$ onto, $c(0)=0$, $Dc(0)=0$, and $M=\tfrac12D^2c(0)\succ0$. On a sufficiently small closed action ball let $E(D)$ be the minimum of $c(u)$ subject to $\Phi(u)=d_D=(0,D)$. Set
\[
 \eta_D=(JM^{-1}J^*)^{-1}d_D,\quad v_D=M^{-1}J^*\eta_D,\quad
 G=HM^{-1/2}P_{\ker(CM^{-1/2})}M^{-1/2}H^*.
\]
For small $D$ the minimizer is unique, and
\[
 \begin{aligned}
 E(D)&=\langle D,G^{-1}D\rangle+\chi(D)+o(\|D\|^3),\\
 \chi(D)&=\tfrac16D^3c(0)[v_D,v_D,v_D]
          -\langle\eta_D,D^2\Phi(0)[v_D,v_D]\rangle .
 \end{aligned}
\]
The remainder is $O(\|D\|^4)$ if the degree-three Taylor polynomial of $c$ at zero has error $O(\|u\|^4)$, as for a quadratic cost or locally Lipschitz $D^3c$. At fixed target coordinates, $\chi$ is invariant under origin-preserving $C^3$ changes of action coordinates transporting both $c$ and $\Phi$. For the squared geodesic distance $c(u)=d_g(0,u)^2$ of a smooth Riemannian metric, $\chi(D)=-\langle\eta_D,(\nabla d\Phi)_0[v_D,v_D]\rangle$, with the source Levi--Civita and flat target connections.
\end{theorem}
\begin{proof}
On a sufficiently small ball $c(u)$ lies between positive multiples of $\|u\|^2$. Surjectivity gives a local right inverse of $\Phi$ through zero and hence a feasible $O(\|D\|)$ action. Compactness gives a minimizer, and the two cost bounds put every minimizer within $O(\|D\|)$ of zero. Each is interior for small $D$. Stationarity reads $\tfrac12\nabla c(u)=D\Phi(u)^*\lambda$ and $\Phi(u)=d_D$. Full row rank near zero gives $\lambda=O(\|u\|)$, since $(D\Phi D\Phi^*)^{-1}D\Phi$ is locally bounded. The linearization in $(u,\lambda)$ at zero is
$\bigl[\begin{smallmatrix}M&-J^*\\J&0\end{smallmatrix}\bigr]$, invertible by $M\succ0$ and surjectivity. Thus all minimizers lie on the unique local $C^2$ stationary branch \cite{bonnans,jackson}, with
$u=v_D+w_D$, $w_D=O(\|D\|^2)$, and $\lambda=\eta_D+O(\|D\|^2)$.

Expanding the constraint gives $Jw_D=-\tfrac12D^2\Phi(0)[v_D,v_D]+O(\|D\|^3)$. In the quadratic part of $c$, the cross term is $2\langle Mv_D,w_D\rangle=2\langle\eta_D,Jw_D\rangle$. Terms containing $w_D$ in the cubic polynomial, as well as $\langle w_D,Mw_D\rangle$, are $O(\|D\|^4)$. Lemma~\ref{m:lem:conditioning} in $M^{1/2}$-whitened coordinates gives $\langle v_D,Mv_D\rangle=\langle D,G^{-1}D\rangle$. Taylor's formula for $c$ now gives both remainder estimates.

For $u=\zeta(w)$ set $A=D\zeta(0)$, $B_2=D^2\zeta(0)$ and $\tilde v=A^{-1}v_D$. The pulled-back metric is $A^*MA$ and the multiplier stays $\eta_D$. The cubic cost term gains $\langle Mv_D,B_2[\tilde v,\tilde v]\rangle$, while the constraint term gains $-\langle\eta_D,JB_2[\tilde v,\tilde v]\rangle$. They cancel. The minimum value is identical on transported neighborhoods, so its quantitative remainder also persists. In geodesic normal coordinates the squared distance is quadratic and the ordinary Hessian of $\Phi$ at the origin equals its covariant Hessian.
\end{proof}
\begin{corollary}[Hidden-state separation beyond the Gram]\label{m:cor:cubicseparation}
For two behavior-matched states satisfying Theorem~\ref{m:thm:cubic} in comparable charts with the same $G$, suppose $\chi_1(d)\ne\chi_2(d)$ in some direction $d$. For all sufficiently small $t>0$, the sign of $E_1(td)-E_2(td)$ equals that of $\chi_1(d)-\chi_2(d)$. A budget strictly between the two costs distinguishes the exact affordable sets.
\end{corollary}
\begin{proof}
The common quadratic term cancels, leaving $(\chi_1(d)-\chi_2(d))t^3+o(t^3)$. Both minima are attained.
\end{proof}
\begin{corollary}[Ray aperture expansion]\label{m:cor:cubicaperture}
In the chart of Theorem~\ref{m:thm:cubic}, fix $d\ne0$, $q=\langle d,G^{-1}d\rangle$ and $\chi=\chi(d)$. Let $\alpha(b;d)$ be the supremum of $t\ge0$ reachable in the chart at cost $c(u)\le b$ with $\Phi(u)=(0,td)$. Then
\[
 \alpha(b;d)=\sqrt{b/q}-\frac{\chi}{2q^2}b+o(b)\qquad (b\downarrow0).
\]
The error is $O(b^{3/2})$ under the fourth-order remainder bound of the theorem.
\end{corollary}
\begin{proof}
Differentiating the attained cost along the stationary branch gives $\frac{d}{dt}E(td)=2\langle\lambda(td),(0,d)\rangle=2qt+O(t^2)$. Small budgets confine all actions to an $O(\sqrt b)$ neighborhood of zero, so the aperture is the local positive inverse of $E(td)=b$. Both the inverse and $t_0=\sqrt{b/q}-\chi b/(2q^2)$ are asymptotic to $\sqrt{b/q}$. Substitution gives $E(t_0d)-b=o(b^{3/2})$, or $O(b^2)$ under the stronger bound. Dividing by the derivative at an intermediate point gives the stated errors.
\end{proof}
\begin{theorem}[Finite nonlinear separation at a common affine tuple]\label{m:prop:tuplematched}
Fix a common closed action ball $\overline B_{r_0}(0,0)\subset\mathbb R^2$ of radius $r_0>0$, with physical cost $x^2+y^2$. For $a\in[0,A]$ define $p_a(x,y)=x+y+ay^2$ and $r_a(x,y)=x$. At the origin all charts have the same affine compensation data,
\[
C=B=H=1,\quad\Psi=0,\quad Q=1,\quad N=0,\quad L=-1,\quad W=2,\quad G_+=1/2,
\]
yet for small positive $D$,
\[
\boxed{E_a(D)=2D^2+2aD^3+O(D^4),\quad
\alpha_a(b;1)=\sqrt{b/2}-\tfrac a4 b+O(b^{3/2}).}
\]
For any $0\le a_1<a_2\le A$, the two scalar response sets at budget $b$ are incomparable for all sufficiently small $b>0$.
\end{theorem}
\begin{proof}
For $b<r_0^2$, a write of cost at most $b$ lies in the common ball. As $b\downarrow0$ these writes approach the origin. The obligation fixes $x=D$. The unique protected root near zero of $D+y+ay^2=0$ is $y=-D-aD^2+O(D^3)$, so $D^2+y^2=2D^2+2aD^3+O(D^4)$. The common Jacobian at zero gives the displayed affine tuple. The aperture correction follows from Corollary~\ref{m:cor:cubicaperture}. On the ray $-1$ the cubic cost coefficient changes sign, giving $\alpha_a(b;-1)=\sqrt{b/2}+ab/4+O(b^{3/2})$. Thus increasing $a$ shrinks the positive extent and increases the negative extent for all sufficiently small $b$ at fixed $a_1<a_2$. Additive physical writes $(u,v)$ on the state $(x,y,a)$ with readout $x+y+ay^2$ realize the family in common cost coordinates.
\end{proof}
\begin{remark}[Response curvature and target coordinates]\label{m:remark:5_10}
For $p=x$ and cost $x^2+y^2$, take $r=y+ay^2$ with $a\ge0$. At $G=1$, local inversion gives $E_a(D)=D^2-2aD^3+O(D^4)$. For small $b>0$ the positive-ray aperture is exactly $\sqrt b+ab$, attained at $y=\sqrt b$. Changing the response to $r=y+ay^3$ makes $\chi=0$ and gives $E_a(D)=D^2-2aD^4+O(D^6)$. Both families admit behavior-matched parents $(0,a,0)$ under the write $(\theta,a,z)\mapsto(\theta+x,a,z+r(y))$ with protected readout $\theta$ and response readout $z$.

For two costs with equal quadratic part, a common target-coordinate change $D=A_O z+O(\|z\|^2)$ transports their cubic difference to $\Delta\chi(A_O z)$. Indeed, substitution in $E_1(D)-E_2(D)=\Delta\chi(D)+o(\|D\|^3)$ gives $\Delta\chi(A_O z)+o(\|z\|^3)$.
\end{remark}

\subsection{Singular structure and sequential execution}
\paragraph{Degenerate tangent geometry.}
On the protected line $p(x,u)=x+u=0$, the response $r(x,u)=u^3$ has zero linear Gram and exact squared cost $E(D)=2|D|^{2/3}$. Conversely $p(y)=y_1^2+y_2^2$ has $\ker Dp(0)=\mathbb R^2$ and an isolated protected point. Even agreement through an arbitrary finite jet permits different exact responses. For an order $j$ choose $2n>j$. The maps $p_0(x,a)=x^2$ and $p_1(x,a)=x^2+a^{2n}$ have the same $j$-jet at zero, but the common response $r=a$ ranges over a neighborhood of zero for $p_0$ and equals zero on $p_1^{-1}(0)$.

\begin{proposition}[Exact singular cooperative cost]\label{m:prop:singularcost}
For $m\ge2$ let
\[
p(y)=(y_1^2-y_m^3,\ldots,y_{m-1}^2-y_m^3),\qquad r(y)=y_1.
\]
From zero, no \emph{finite sequence} of writes, each changing at most $m-1$ coordinates and preserving exact protection at every stage, can reach a nonzero endpoint. A synchronous protected write delivering $D$ exists and has exact squared endpoint cost
\[
\boxed{E_m(D)=(m-1)D^2+|D|^{4/3}.}
\]
With independent coordinate squared budgets $b_i$, the endpoint is affordable precisely for
\[
|D|\le\min\{\sqrt{b_1},\ldots,\sqrt{b_{m-1}},b_m^{3/4}\}.
\]
\end{proposition}
\begin{proof}
If the final coordinate is fixed at zero, protection forces all first $m-1$ coordinates to vanish. If any of those coordinates is fixed at zero, its constraint forces $y_m=0$, and hence all others vanish. Thus no nontrivial protected endpoint can omit any coordinate. For a synchronous write, $y_1=D$ implies $y_m^3=D^2$, so $y_m=|D|^{2/3}$ and $y_i=\pm D$ for all $i<m$. Squaring gives the exact total cost and the independent-coordinate constraints. Induction proves the sequential impossibility from the zero state.
\end{proof}
\begin{remark}[A singular zero-gap discontinuity]\label{m:remark:6_5}
On a ball containing $(1,t)$, $p(x,y)=x(x-1)$ and $r(x,y)=xy$ give $E(t)=1+t^2$ for small $t>0$, while $E(0)=0$. A semialgebraic cost germ \cite{basu} with a positive leading power on a feasible target ray satisfies $E(td)\to0$ there.
\end{remark}

\begin{definition}[Related-rank and block-activation gate sets]\label{m:definition:6_6}
In a finite-dimensional Euclidean block model, let $C,B_i,H,\Psi_i$ be fixed linear maps. The protection is $p(t,x,y)=\Phi(t)(Cx+\sum_iB_iy_i)$ from $(t_*,0,0)$. At the returned gate $t_*$, the obligation increment is $\bar H(x,y)=Hx+\sum_i\Psi_i y_i$. Use the joint/separate spaces of Theorem~\ref{m:thm:alg}. Put
\[
R_{\rm rel}=\Ran C+\sum_i\Ran B_i,\qquad
\Sigma_{\rm rel}=\{t:\ker\Phi(t)\cap R_{\rm rel}\ne\{0\}\},
\]
\[
\Sigma_{\rm act}=\bigcup_i\{t:\exists y_i,\ \Phi(t)B_iy_i=0,\ [B_iy_i]_{\calP/\Ran C}\ne0\},\quad
\Sigma_{\rm full}=\{t:R_{\rm rel}\subseteq\ker\Phi(t)\}.
\]
Write $\kappa_\square=\dist(t_*,\Sigma_\square)^2$, with $+\infty$ for an empty set. If $t\mapsto\Phi(t)|_{R_{\rm rel}}$ is operator-norm continuous near $t_*$, $R_{\rm rel}$ is finite dimensional, and $\Phi(t_*)|_{R_{\rm rel}}$ is injective, then $\kappa_{\rm rel}>0$: injectivity persists in a neighborhood by the positive smallest singular value.
\end{definition}
\begin{theorem}[Necessary block activation and a sufficient gate path]\label{m:thm:gate}
Assume the initial gate restriction on $R_{\rm rel}$ is injective, the protected protocol may change at most one of $t,x,y_1,\ldots,y_m$ at each stage, and every intermediate state is exactly protected. Use an anchored peak budget $\sup_s(\|t_s-t_*\|^2+\|x_s\|^2+\sum_i\|y_{i,s}\|^2)\le b$. Then
\[
b\le\kappa_{\rm act}\quad\Longrightarrow\quad
\text{every reachable terminal new-block residual is separately compatible.}
\]
If the joint obligation-response space strictly exceeds the separately compatible response space, then $\kappa_{\rm rel}\le\kappa_{\rm act}\le\kappa_{\rm full}$. Conversely, choose a jointly protected witness $v=(x_v,y_v)$ with $\bar Hv\notin V_{\rm sep}$. Suppose a full-collapse gate $t_\circ\in\Sigma_{\rm full}$ has a route $\gamma:[0,1]\to\calY_0$ from $t_*$ to $t_\circ$ with peak $K_\gamma=\sup_s\|\gamma(s)-t_*\|^2<b$. Require a number $\varepsilon_0>0$ such that, for every $0<\varepsilon\le\varepsilon_0$, the outward gate path at zero write, the individual straight writing segments to $\varepsilon v$ at $t_\circ$, and the reversed gate path at the \emph{new frozen write} $\varepsilon v$ are all path-admissible. Then a nonseparable response can be constructed. The construction uses at most $m+3$ block-level moves, counting the outward and return gate paths as one move each.
\end{theorem}
\begin{proof}
For each $Y_i$ stage at fixed $t$, protection of before and after states forces $\Phi(t)B_i\Delta y_i=0$. Unless $t\in\Sigma_{\rm act}$, its residual class modulo $\Ran C$ is zero. At $b\le\kappa_{\rm act}$, any stage exactly at an activating gate exhausts the peak budget in gate displacement and so must have $x=y=0$, preventing nonzero activation. Induction preserves separate compatibility. The set inclusions $\Sigma_{\rm full}\subseteq\Sigma_{\rm act}\subseteq\Sigma_{\rm rel}$ hold whenever a nontrivial residual channel exists, giving the distance inequalities. For the sufficient direction travel along $\gamma$ at zero write, write a small multiple of a nonseparable jointly compatible witness one block at a time at full collapse, then return along the reversed route. Choose $0<\varepsilon\le\varepsilon_0$ with $\varepsilon^2\|v\|^2<b-K_\gamma$. On each writing segment, the accumulated block squared norm is at most $\varepsilon^2\|v\|^2$, and the protection map vanishes on the relevant range at $t_\circ$. On return, $Cx_v+\sum_iB_i y_{v,i}=0$, so protection holds for every gate value. The peak is at most $K_\gamma+\varepsilon^2\|v\|^2<b$ and the final output $\varepsilon\bar Hv$ lies outside $V_{\rm sep}$.
\end{proof}
\begin{example}[Rank loss without execution]\label{m:counterexample:6_8}
Take $t_*=1$, $Cx=(x,0)$, $B_1y_1=(0,y_1)$, $B_2y_2=(0,-y_2)$ and
\[
\Phi(t)=t\begin{pmatrix}1&-1\\t-0.9&0\end{pmatrix}.
\]
At $t=0.9$ the relevant operator loses rank ($\kappa_{\rm rel}=0.01$), but neither new block is individually activatable until $t=0$, so $\kappa_{\rm act}=\kappa_{\rm full}=1$.
\end{example}

\begin{remark}[Activation and return at a pure-block gate]\label{m:remark:6_10}
For the zero-parent, pure-block executor, $\Sigma_{\rm rel}$ detects loss of injectivity on the residual range and $\Sigma_{\rm act}$ detects a writable block with a nonzero residual class. At a gate in $\Sigma_{\rm full}$, the admissible return path carries the jointly protected write back to $t_*$, where its obligation response is evaluated.
\end{remark}

\section{Structural growth and target-conditioned improvement}\label{m:sec:growth}
\paragraph{Objectives determined by the protected readout.}
A loss of the form $\ell=\widetilde\ell\circ P$ is constant on every fiber of $P$. If $P$ is differentiable at $q$ and $\widetilde\ell$ at $P(q)$, then $d\ell_q=d\widetilde\ell_{P(q)}\circ dP_q$. Hence $d\ell_q$ vanishes on $\ker dP_q$.

\begin{proposition}[Structural capacity can grow or saturate]\label{m:prop:scaling}
For nested \emph{closed} protected subspaces $\calK_n$ of a common cost Hilbert space and a fixed bounded total response $H$, set $G_n=HP_{\calK_n}H^*$. Then $G_{n+1}-G_n\succeq0$. For a rank-one increment $G_{n+1}=G_n+w_nw_n^*$ and a fixed nonzero target $d$, assume a finite-dimensional common response space with $G_n\succ0$ for $n\ge n_0$ (or, in Hilbert output, boundedly invertible $G_n$). For $b>0$, set
\[
q_n=\langle d,G_n^{-1}d\rangle,\quad
\lambda_n(d)=\frac{|\langle d,G_n^{-1}w_n\rangle|^2}{q_n(1+\langle w_n,G_n^{-1}w_n\rangle)}\in[0,1).
\]
Then the linear aperture obeys the exact product law
\[
\frac{\alpha_N(b;d)}{\alpha_{n_0}(b;d)}=\prod_{j=n_0}^{N-1}(1-\lambda_j(d))^{-1/2}.
\]
The single-operator bound is $\sup_{D\in\calA_{n,b}}\|D\|\le\sqrt b\|H\|$. Across finite systems with varying response operators, any nonnegative, nondecreasing sequence of scalar squared apertures is possible at fixed positive budget.
\end{proposition}
\begin{proof}
Orthogonal projection onto $\calK_{n+1}=\calK_n\oplus(\calK_{n+1}\cap\calK_n^\perp)$ gives the positive Gram increment. Sherman--Morrison applied to $G_n+w_nw_n^*$ yields $q_{n+1}=q_n(1-\lambda_n)$ and hence the aperture product. Cauchy--Schwarz in the $G_n^{-1}$ metric bounds the numerator of $\lambda_n$ by $q_n\langle w_n,G_n^{-1}w_n\rangle$, hence $0\le\lambda_n<1$. The uniform norm bound follows from $\|Hu\|\le\|H\|\|u\|$. For any nondecreasing $s_n$ with $s_0=0$ use finite-dimensional modules $p_i=u_i+v_i$, $r_n=\sum_i\sqrt{2(s_i-s_{i-1})}\,v_i$ and cost $\sum_i(u_i^2+v_i^2)$. The Gram is exactly $s_n$, and $\alpha_n(b;1)=\sqrt{bs_n}$. Unbounded $s_n$ yields unbounded operator norms $\sqrt{s_n}$ in the protected coordinates.
\end{proof}

A uniform bound on the correction normal to the protected tangent space controls the loss of descent at finite radius. We measure descent against the objective pulled back to the action chart \cite{boumal}.
\begin{theorem}[Uniform exact-protection opportunity bound with explicit constants]\label{m:thm:uniform}
Fix $R>0$ and let $\{(U_i,P_i,\calO,p_i,r_i)\}_{i\in I}$ be any finite or countable family of finite-dimensional Euclidean action charts. Suppose $p_i(0)=r_i(0)=0$ and for every $i$:
\begin{enumerate}[label=(\roman*),leftmargin=*]
\item $p_i,r_i$ are $C^{1,1}$ on an open neighborhood of the closed $R$-ball; $C_i=Dp_i(0)$ is surjective and $\|C_i^\dagger\|\le a$;
\item $H_i=Dr_i(0)$ obeys $\|H_i\|\le h$, and for $\|u\|\le R$,
\[
\|Dp_i(u)-C_i\|\le M_p\|u\|,\qquad
\|Dr_i(u)-H_i\|\le M_r\|u\|;
\]
\item the action domain is the closed $R$-ball, with cost $\|u\|^2$ and protection equation $p_i(u)=0$. The common current objective $\ell:\calO\to\R$ is defined and differentiable on an open set containing the response segments, has $g=\nabla\ell(0)$ with $\|g\|\le g_0$, and its gradient is $\beta_\ell$-Lipschitz there.
\end{enumerate}
Set $k=aM_p$, $h_R=h+M_rR/2$ and
\begin{equation}\label{m:eq:Kexplicit}
K_{\rm unif}:=\frac{g_0M_r}{2}+\frac{\beta_\ell h_R^2}{2}+2g_0hk.
\end{equation}
Choose $\rho_*=\min\{R/2,1/(8k)\}$, with $1/0=+\infty$ when $k=0$. Write
\[
G_i=H_iP_{\ker C_i}H_i^*,\qquad
\kappa_i=\sqrt{\langle g,G_i g\rangle},
\]
and let $\Delta_i(\rho)=\ell(0)-\min\{\ell(r_i(u)):p_i(u)=0,\|u\|\le\rho\}$. Set $t(\rho)=\rho/(1+k\rho)$. For $\kappa_i>0$ the lower bound uses the protected normal-correction curve constructed below. Then, \emph{uniformly in $i$},
\begin{equation}\label{m:eq:uniform}
\boxed{\left|\Delta_i(\rho)-\rho\kappa_i\right|\le K_{\rm unif}\rho^2,\qquad 0<\rho\le\rho_*.}
\end{equation}
Consequently for two charts $i=0,1$ with $\kappa_1>\kappa_0$,
\[
\Delta_1(\rho)-\Delta_0(\rho)\ge(\kappa_1-\kappa_0)\rho-2K_{\rm unif}\rho^2,
\]
so strict improvement holds for every $0<\rho\le\rho_*$ with $2K_{\rm unif}\rho<\kappa_1-\kappa_0$.
\end{theorem}
\begin{proof}
Fix a chart and drop its index. The derivative bound implies $\|r(u)-Hu\|\le M_r\|u\|^2/2$ and $\|r(u)\|\le h_R\|u\|$ for $\|u\|\le R$. The objective descent lemma gives
\begin{equation}\label{m:eq:objRem}
|\ell(r(u))-\ell(0)-\langle g,Hu\rangle|
\le E_0\|u\|^2,\quad E_0=\tfrac12g_0M_r+\tfrac12\beta_\ell h_R^2.
\end{equation}
For an exactly protected $u$, $Cu=-[p(u)-Cu]$, so with $K=\ker C$,
\[
\|P_{K^\perp}u\|=\|C^\dagger Cu\|\le \frac{k}{2}\|u\|^2.
\]
Writing $u=P_Ku+P_{K^\perp}u$, equations \eqref{m:eq:objRem} and Cauchy--Schwarz show every feasible $\|u\|\le\rho$ has decrease at most
\[
\ell(0)-\ell(r(u))\le\rho\kappa+(E_0+\tfrac12g_0hk)\rho^2.
\]
For the converse suppose $\kappa>0$ and choose the unit descending tangent $v=-P_KH^*g/\kappa$. Define $N(w)=p(w)-Cw$. A normal correction $w_t\in K^\perp$ must solve
\[
w_t=-C^\dagger N(tv+w_t).
\]
Take $w_0=0$. If $k=0$, the derivative bound gives $C^\dagger N=0$ on the ball, so $N=0$ by $CC^\dagger=I$, and $w_t=0$ throughout. For $k>0$ and $0<t\le t(\rho_*)$, we have $kt\le1/8$. The map on $\{w\in K^\perp:\|w\|\le kt^2\}$ sends the ball into itself because
\[
\|C^\dagger N(tv+w)\|\le\frac k2(t+kt^2)^2
\le\frac{81}{128}kt^2<kt^2.
\]
The bound $t+kt^2\le9R/16<R$ keeps the entire correction ball inside the original chart. Its Lipschitz constant is at most $k(t+kt^2)\le9/64<1$, so the correction exists uniquely in this normal ball, with $\|w_t\|\le kt^2$. Since $CC^\dagger=I$, the fixed-point equation gives $p(tv+w_t)=0$. For continuity, put $T=t(\rho_*)$ and $q=k(T+kT^2)\le9/64$. On the radius-$(T+kT^2)$ ball, $C^\dagger N$ is $q$-Lipschitz. The fixed-point equations for $s,t\in[0,T]$ therefore give
\[
 (1-q)\|w_s-w_t\|\le q|s-t|.
\]
At budget radius $\rho$ use $t=t(\rho)$. Every intermediate point $u_s=sv+w_s$ with $0\le s\le t$ satisfies $p(u_s)=0$ and $\|u_s\|\le s+ks^2\le t+kt^2\le\rho$. The entire path lies in the budget ball. By \eqref{m:eq:objRem},
\[
\ell(0)-\ell(r(u_t))\ge\kappa t-g_0hk\rho^2-E_0\rho^2.
\]
Using $t\ge\rho-k\rho^2$ and $\kappa\le g_0h$ yields
\[
\Delta(\rho)\ge\kappa\rho-(E_0+2g_0hk)\rho^2.
\]
When $\kappa=0$, the zero action and the upper estimate give $0\le\Delta(\rho)\le(E_0+\tfrac12 g_0hk)\rho^2\le K_{\rm unif}\rho^2$. At $\rho=0$ the sole action is zero and $\Delta(0)=0$. The feasible set is nonempty and compact in each finite-dimensional closed ball, so the minimum in the statement is attained. Both inequalities are uniform over the entire family, proving \eqref{m:eq:uniform} and the comparison.
\end{proof}
\begin{remark}[Dependence on chart conditioning]\label{m:remark:7_3}
The protection conditioning and protection curvature enter through the product $k=aM_p$. It sets the certified radius $\rho_*=\min\{R/2,1/(8k)\}$ and contributes $2g_0hk$ to the quadratic error coefficient in \eqref{m:eq:Kexplicit}. At fixed $R$, once $k>1/(4R)$ the certified radius decreases as $1/(8k)$.
\end{remark}

\begin{example}[Loss of uniformity under an ill-conditioned protection map]\label{m:ex:eps}
For $\varepsilon>0$, take the smooth chart $p_\varepsilon(x,y)=\varepsilon x+y^2$, $r_\varepsilon(x,y)=y$, squared cost $x^2+y^2$, and linear loss $\ell(z)=-z$. For every fixed $\varepsilon$, $Dp_\varepsilon(0)=[\varepsilon\ 0]$ is surjective, $\ker Dp_\varepsilon(0)=\R(0,1)$, and the protected response Gram is \emph{identically} $G_\varepsilon=1$. Nevertheless exact protection requires $x=-y^2/\varepsilon$, and the best decrease under budget $\rho^2$ is
\[
\Delta_\varepsilon(\rho)=\left[\frac{\varepsilon^2}{2}\left(\sqrt{1+\frac{4\rho^2}{\varepsilon^2}}-1\right)\right]^{1/2}.
\]
For every \emph{fixed} $\rho>0$, $\Delta_\varepsilon(\rho)\to0$ as $\varepsilon\downarrow0$, although the linear coefficient $\kappa_\varepsilon=1$ at every fixed $\varepsilon$. Hence there is no family-uniform $K,\rho_*>0$ with $|\Delta_\varepsilon(\rho)-\rho|\le K\rho^2$ for all $\varepsilon,\rho\le\rho_*$. Here $\|Dp_\varepsilon(0)^\dagger\|=1/\varepsilon\to\infty$, violating uniform conditioning.
\end{example}

\begin{theorem}[Paid installation at a common write budget]\label{m:thm:commission}
From one parent compare a single write in chart $0$ with installation followed by a single write in chart $1$. The installation is admissible, preserves protection and the current objective baseline, and has record $(c,c)$ in an additive scalar reserve. In either chart a write $u$ has record $(\|u\|^2,\|u\|^2)$. Suppose both charts satisfy Theorem~\ref{m:thm:uniform} for the same objective on their common response space, with constants $K_{\rm unif},\rho_*$. For $0\le c\le B\le\rho_*^2$, put $r_0=\sqrt B$, $r_1=\sqrt{B-c}$. The optimal decreases on the two exact-protection endpoint sets satisfy
\begin{equation}\label{m:eq:installcriterion}
\boxed{\Delta_1(r_1)-\Delta_0(r_0)\ge
\kappa_1\sqrt{B-c}-\kappa_0\sqrt B-K_{\rm unif}(2B-c).}
\end{equation}
A positive right-hand side gives a strict advantage for installation. When $\kappa_0=0$ the criterion reduces to $\kappa_1\sqrt{B-c}>K_{\rm unif}(2B-c)$.
\end{theorem}
\begin{proof}
The reserve composition in Proposition~\ref{m:prop:record} leaves $B-c$ for the installed write. Both endpoint optima are attained on compact protected sets, and baseline preservation compares their decreases from the same value. Subtract the upper estimate $\Delta_0(r_0)\le\kappa_0 r_0+K_{\rm unif}r_0^2$ from the lower estimate $\Delta_1(r_1)\ge\kappa_1 r_1-K_{\rm unif}r_1^2$ of Theorem~\ref{m:thm:uniform}. The zero-radius value is $\Delta_i(0)=0$, including the installed write at $B=c$.
\end{proof}
If installation shifts the objective baseline from $\ell_{\rm old}$ to $\ell_{\rm new}$, the end-to-end decrease is the installed write's decrease plus $\ell_{\rm old}-\ell_{\rm new}$.

\Needspace{6\baselineskip}
\subsection{Nonasymptotic comparison of directed and random acquisition}
\begin{theorem}[Nonasymptotic directional comparison]\label{m:thm:random}
Let $R$ be an $m$-dimensional Euclidean space with $m\ge1$, $A\succ0$ and $a\ne0$. For unit $v$ define
\[
 \Gamma(v)=\frac{\langle a,v\rangle^2}{\langle v,Av\rangle},\quad
 \Gamma_{\rm opt}=\langle a,A^{-1}a\rangle,\quad
 \tau=\Tr A,\quad r_{\eff}(A)=\tau/\|A\|_{\op}.
\]
Then $\max_v\Gamma(v)=\Gamma_{\rm opt}$, achieved at $v\propto A^{-1}a$. If $v$ is uniform on $\mathbb S(R)$, then for all $\delta\in(0,1)$,
\begin{equation}\label{m:eq:smallball}
\boxed{\mathbb P\{\Gamma(v)\ge (\delta^2/4)\|a\|^2/\tau\}\ge1-\delta.}
\end{equation}
With probability at least $1-\delta-e^{-r_{\eff}(A)/16}$,
\begin{equation}\label{m:eq:ratio}
 \Gamma(v)\le4\log(2/\delta)\frac{\|a\|^2}{\tau},\qquad
 \frac{\Gamma_{\rm opt}}{\Gamma(v)}\ge\frac{r_{\eff}(A)}{4\log(2/\delta)}.
\end{equation}
Here $r_{\eff}(A)$ is the stable rank of $A^{1/2}$, and the ratio is taken as $+\infty$ on the null event $\Gamma(v)=0$.
\end{theorem}
\begin{proof}
Cauchy--Schwarz in the $A$ metric gives $\langle a,v\rangle^2\le\langle a,A^{-1}a\rangle\langle v,Av\rangle$, with equality at the stated direction. For the random bound write $v=Z/\|Z\|$, $Z\sim N(0,I_R)$, $u=a/\|a\|$, $z=\langle Z,u\rangle$, and $w=Z-zu$. The scalar $z\sim N(0,1)$ is independent of $w$. With $Q=\langle Z,AZ\rangle$ the quotient is $\|a\|^2z^2/Q$. If it is below $(\delta^2/4)\|a\|^2/\tau$, the triangle inequality and $\langle u,Au\rangle\le\tau$ imply
\[
 |z|<\frac{\delta}{2\sqrt\tau}\left(|z|\sqrt{\langle u,Au\rangle}+\sqrt{\langle w,Aw\rangle}\right)
 \quad\Longrightarrow\quad
 |z|<\frac{\delta}{2-\delta}\sqrt{\frac{\langle w,Aw\rangle}{\tau}}.
\]
Conditional on $w$, the Gaussian density bound gives $\mathbb P(|z|<t)\le\sqrt{2/\pi}\,t$. Jensen and $\mathbb E\langle w,Aw\rangle=\tau-\langle u,Au\rangle\le\tau$ bound the bad probability by $\sqrt{2/\pi}\delta/(2-\delta)<\delta$. For $m=1$ the bad event is empty; this also follows from the same inequality.

For the lower tail \cite[Lemma~1, (4.2)]{laurentmassart}, diagonalize $A$ with eigenvalues $\lambda_j>0$, and put $\nu=\sum_j\lambda_j^2$. The elementary inequality $\log(1+2x)\ge2x-2x^2$ for $x\ge0$ yields
\[
 \mathbb E e^{-\theta Q}=\prod_j(1+2\theta\lambda_j)^{-1/2}
 \le e^{-\theta\tau+\theta^2\nu}.
\]
Markov's inequality at $\theta=\tau/(4\nu)$ therefore gives
\[
 \mathbb P(Q\le\tau/2)\le e^{-\tau^2/(16\nu)}
 \le e^{-r_{\eff}(A)/16},
\]
since $\nu\le\|A\|_{\op}\tau$. The standard Gaussian bound gives $\mathbb P\{z^2>2\log(2/\delta)\}\le\delta$. The union bound on the two failure events, together with $\Gamma_{\rm opt}\ge\|a\|^2/\|A\|_{\op}$, gives \eqref{m:eq:ratio}.
\end{proof}
\begin{corollary}[Objective scaling versus Gram-capacity scaling]\label{m:cor:sqrt}
Suppose $\langle g,G_0g\rangle=0$ and a one-channel acquisition satisfies $\langle g,\Delta G_v g\rangle=c^2\Gamma(v)$ with the same $c>0$ for all directions. Write $v_*\propto A^{-1}a$ for a unit maximizing direction of Theorem~\ref{m:thm:random}. For $\Gamma(v)>0$, the ratio of the leading loss-decrease coefficients is $\sqrt{\Gamma_{\rm opt}/\Gamma(v)}$. If both charts satisfy Theorem~\ref{m:thm:uniform} with common constants $K_{\rm unif},\rho_*$, their finite-radius decreases obey
\[
\left|\Delta_{v_*}(\rho)-\Delta_v(\rho)
 -\rho c\bigl(\sqrt{\Gamma_{\rm opt}}-\sqrt{\Gamma(v)}\bigr)\right|
 \le 2K_{\rm unif}\rho^2,\qquad 0<\rho\le\rho_* .
\]
\end{corollary}
\begin{proof}
On the old plateau, $\kappa_v=\sqrt{\langle g,(G_0+\Delta G_v)g\rangle}=c\sqrt{\Gamma(v)}$. Division of the two positive coefficients gives the ratio. Subtraction of the two expansions in \eqref{m:eq:uniform} gives the finite-radius error bound.
\end{proof}
\begin{remark}[Complete compensating numerator]\label{m:remark:8_3}
For the finite-dimensional physical model $Cx+vy=0$ with $v\in R=\Ran C$, $\|v\|=1$, scalar $y$ and Euclidean cost $\|x\|^2+y^2$, let the possibly vector-valued response be $r(x,y)=Hx+(\Psi_0v)y$. For target gradient $g\in\calO$, define
\[
 A=I_R+(C^\dagger|_R)^*(C^\dagger|_R),\qquad
 a_g=P_R\{\Psi_0^*g-(C^\dagger)^*H^*g\}.
\]
The compensated column is $\ell_v=\Psi_0v-HC^\dagger v$, its metric is $1+\|C^\dagger v\|^2$, and inherited orthogonal compensation gives
\[
 \Delta G_v=\frac{\ell_v\ell_v^*}{1+\|C^\dagger v\|^2},\qquad
 \langle g,\Delta G_vg\rangle=\frac{|\langle a_g,v\rangle|^2}{\langle v,Av\rangle}.
\]
The probabilistic bounds of Theorem~\ref{m:thm:random} apply when $a_g\ne0$ and $v$ is uniform on the unit sphere of $R$. In the scalar case $\Psi_0v=\langle h_\Psi,v\rangle$ this restores $a=P_R(h_\Psi-(C^\dagger)^*H^*)$.
\end{remark}

\Needspace{7\baselineskip}
\section{Exact state descriptions and Hilbert-space limits}\label{m:sec:closure}
\subsection{Exact controlled-state quotients}
\begin{theorem}[Controlled congruence]\label{m:thm:congruence}
Let $\mathfrak d_{\emptyset}(q)$ collect the readouts of the empty protocol at $q\in\mathcal S$. Let $\mathcal U(q)$ be the available action labels. Each $a\in\mathcal U(q)$ has immediate data $\mathfrak d_a(q)$, including protected and obligation readouts, the resource record $(m_a,c_a)$ and guards, and a deterministic successor $F_a(q)$. A descriptor $\varphi:\mathcal S\to\mathcal M$ admits exact quotient readouts, action availability and reduced successors if and only if for every $q,q'$ with $\varphi(q)=\varphi(q')$,
\[
\mathfrak d_{\emptyset}(q)=\mathfrak d_{\emptyset}(q'),\qquad \mathcal U(q)=\mathcal U(q'),
\]
and, for every $a\in\mathcal U(q)$,
\[
\mathfrak d_a(q)=\mathfrak d_a(q'),\qquad \varphi(F_a(q))=\varphi(F_a(q')).
\]
Every finite admissible action word descends to the quotient with its ordered reserve record.
\end{theorem}
\begin{proof}
Each quotient function has the same value on representatives of a descriptor fiber, which gives the stated equalities. Conversely, use any representative to define the initial readout, available labels and each labelled transition on $\varphi(\mathcal S)$. The equalities make these definitions independent of the representative. The initial readout agrees for the empty word. Each further label preserves equality of the emitted data and reduced successors, and composes the same resource record on both paths. Induction gives the result for every composable word.
\end{proof}
On each label's domain the quotient commutes with the transition: $\bar F_a\circ\varphi=\varphi\circ F_a$. For the left shift on $\ell^2$ with current readout $h_1$, successive observations reveal every coordinate. An exact descriptor must therefore be injective. Restricting to a $(k+1)$-dimensional ball and applying invariance of domain rules out a continuous injection of any open Hilbert ball into $\mathbb R^k$.

For coordinatewise second-moment updates, the partition subspaces give the invariant linear classes containing $\mathbf1$ that are closed under squaring \cite[Proposition~2.1]{canogar}.
\begin{theorem}[Exact reduction on coordinate partitions]\label{m:thm:partition}
A linear subspace $S\subseteq\R^d$ containing $\mathbf1$ is closed under coordinatewise squaring if and only if it is closed under coordinatewise products, if and only if $S$ is the space of vectors constant on blocks of a partition $\{I_1,\ldots,I_k\}$ of $\{1,\ldots,d\}$. Restrict the state domain to the invariant class in which every vector register, including controller memory, belongs to $S$. If every input-dependent update and intervention preserves this class, and its availability, guards, costs and readouts respect the same fixed known partition, then block values together with all required scalar and discrete registers give an exact recursively closed input-conditioned description. For $v$ vector registers, $n_s$ scalar registers and fixed $k$, the real coordinate count is $vk+n_s$.
\end{theorem}
\begin{proof}
Polarization $2(x\odot y)=(x+y)^{\odot2}-x^{\odot2}-y^{\odot2}$ gives product closure from square closure. Define $i\sim j$ if $x_i=x_j$ for all $x\in S$. A generic linear combination of a basis of $S$ has distinct coordinates on different equivalence blocks, since forbidden coefficient choices lie in finitely many proper hyperplanes. Polynomial interpolation yields each block-indicator vector as a polynomial in this generic vector. Product closure and $\mathbf1\in S$ put every indicator in $S$; these indicators span exactly the block-constant space. The converse is immediate. The invariant-class and partition-compatibility hypotheses permit reconstruction of every register and successor from its block values, so Theorem~\ref{m:thm:congruence} applies. For a gradient-based instance, within-block permutation invariance of an objective makes its gradient block-constant at each block-constant state; compatible coordinatewise controller updates, including optimizer-memory updates, preserve the subspace. If $E$ is the block-indicator matrix, the physical squared norm is $\|E\bar x\|^2=\bar x^*W_I\bar x$ with $W_I=E^*E=\operatorname{diag}(|I_j|)$. The corresponding reduced gradient is $W_I^{-1}\nabla\bar\ell$; retaining block multiplicity preserves the physical cost and gradient dynamics.
\end{proof}
\begin{example}[Equal Grams with different regulated successors]\label{m:counterexample:9_4}
For states $(\theta,h,z)$, let the protected readout be $\theta$, the response readout be $z$, and a write $(x,u)$ increment them by $x$ and $\gamma(h)u$, where $\gamma(h)=1+h(h-1)$. Exact protection requires $x=0$ and the response Gram is $G(h)=\gamma(h)^2$ under unit quadratic cost. With the same $\theta,z$, the states $h=0$ and $h=1$ give identical current Grams and budget images. The same zero-cost regulation $h\mapsto h+1$ gives their successors Grams $1$ and $9$.
\end{example}
\begin{proposition}[Smooth refinement rank obstruction]\label{m:prop:rankgap}
Let $\varphi:U\to\R^r$ and the future response $\Lambda:U\to\R^p$ be smooth on a finite-dimensional state chart $U\subseteq\mathcal S$. At $q\in U$, define
\[
r_{\rm gap}(q)=\operatorname{rank}\begin{pmatrix}D\varphi_q\\ D\Lambda_q\end{pmatrix}-\operatorname{rank}D\varphi_q
=\operatorname{rank}(D\Lambda_q|_{\ker D\varphi_q}).
\]
If a refinement $\psi=(\varphi,\eta)$, $\eta:U\to\R^k$, permits a smooth local decoder $\Lambda=f\circ\psi$, then $k\ge r_{\rm gap}(q)$.
\end{proposition}
\begin{proof}
The rank identity follows from rank--nullity applied to $\ker D\varphi_q$. The chain rule gives $\ker D\psi_q\subseteq\ker D\Lambda_q$. Restrict to $\ker D\varphi_q$: the rank of $D\eta_q$ must cover the rank of $D\Lambda_q$ there, so $k\ge r_{\rm gap}(q)$.
\end{proof}
\subsection{Finite realization from an invariant observable span}
When the future word readouts span a finite-dimensional space, the transition pullbacks act on it by matrices \cite{psr,koopman}.
\begin{proposition}[Finite linear realization]\label{m:prop:observablespan}
For a fixed labelled system of \emph{globally defined} maps $F_a:\mathcal S\to\mathcal S$ on one common invariant state domain, and scalar readouts in $\beta$, let $\mathscr O$ be the real linear span of all coordinates of $\beta\circ F_w$ over finite words $w$, including the empty word. For partial maps, adjoin an absorbing failure state. To represent availability or action costs, include their indicators or values in the family of observables before taking its span. If $\dim\mathscr O=r<\infty$, there exist an exact descriptor $\varphi:\mathcal S\to\R^r$, matrices $A_a$ and readout $T$ such that
\[
\varphi(F_a(q))=A_a\varphi(q),\qquad \beta(q)=T\varphi(q).
\]
The minimal dimension among \emph{linear} realizations of all these probes is $r$. If $\mathscr O_H$ denotes the span of word readouts of length at most $H$, exact stabilization $\mathscr O_{H+1}=\mathscr O_H$ implies stabilization at every later horizon.
\end{proposition}
\begin{proof}
Choose a basis $f_1,\ldots,f_r$ of $\mathscr O$ and set $\varphi=(f_1,\ldots,f_r)$. Prefixing any observable by $F_a$ stays inside $\mathscr O$, hence is represented by a matrix $A_a$; the current readout lies in the span and yields $T$. Conversely a $k$-dimensional linear realization makes every probe a linear combination of at most $k$ coordinates, so $r\le k$. Since $\mathscr O_{H+1}=\mathscr O_0+\sum_a F_a^*\mathscr O_H$, a single exact equality $\mathscr O_{H+1}=\mathscr O_H$ makes $\mathscr O_H$ invariant under every pullback; induction gives all later equalities. Two states have the same basis descriptor precisely when all their finite-word readouts agree.
\end{proof}

\subsection{Hilbert limits and resource-topological closure}
\begin{proposition}[Two spectral limits]\label{m:prop:Hilbertlimits}
Let $T:K\to\mathcal O$ be a bounded protected linear response between Hilbert spaces and $G=TT^*$. Define the extended exact minimum squared cost $E(D)=\inf\{\|u\|^2:Tu=D\}$ with $E(D)=+\infty$ off $\operatorname{Ran}T$. Then
\[
\boxed{E(D)=\lim_{\alpha\downarrow0}\langle D,(G+\alpha I)^{-1}D\rangle.}
\]
If $J_h:\mathcal H\to\mathcal O$ is Hilbert--Schmidt, $\Pi=P_{\overline{\operatorname{Ran}T}}$, and $R_\alpha=\alpha(G+\alpha I)^{-1}$, then
\[
J_h^*R_\alpha J_h\longrightarrow J_h^*(I-\Pi)J_h\quad\text{in trace norm}.
\]
\end{proposition}
\begin{proof}
By polar decomposition, the minimum-norm preimage quadratic form is the inverse-square-root form of $G$, with finite domain $\operatorname{Ran}G^{1/2}=\operatorname{Ran}T$. Spectral monotone convergence of $(\lambda+\alpha)^{-1}$ to $\lambda^{-1}$ for $\lambda>0$ (and to $+\infty$ at $\lambda=0$) proves the first identity, including its extended value. Spectral calculus also gives $R_\alpha\to P_{\ker G}=I-\Pi$ strongly, with $0\preceq R_\alpha\preceq I$. For an orthonormal basis $(e_j)$ of $\mathcal H$, the traces of the positive differences satisfy dominated convergence because $\sum_j\|J_he_j\|^2<\infty$. Positivity identifies the trace with trace norm, giving the result.
\end{proof}
\begin{remark}[Exact range and closed bounded-budget images]\label{m:remark:10_2}
The finite-cost domain $\operatorname{Ran}G^{1/2}=\operatorname{Ran}T$ can strictly contain $\operatorname{Ran}G$ and can be nonclosed. At any fixed finite budget $b\ge0$, its affordable subset $T(\overline B_{\sqrt b})$ is weakly compact and norm closed.
\end{remark}

\begin{proposition}[Two-stage closedness with bounded witnesses]\label{m:prop:closed}
Let the complete-state spaces $X,Y,Z$ and the response space $\calO$ be reflexive Banach spaces, and let $V$ be a reflexive Banach lattice. Let
\[
 \mathcal R_1\subset X\times Y\times\calO\times V_+\times V,\quad
 \mathcal R_2\subset Y\times Z\times\calO\times V_+\times V
\]
be weakly sequentially closed relations, with each record satisfying $c_i\preceq m_i$. Fix $x\in X$ and $b\in V_+$. A terminal triple $(D,z,r)$ is attainable if there exist $y,D_1,D_2,m_1,c_1,m_2,c_2$ with records in these relations and
\[
 m_1\preceq b,\quad m_2\preceq b-c_1,\quad
 D=D_1+D_2,\quad r=b-c_1-c_2.
\]
Assume that each weakly convergent sequence of attainable terminal triples admits witnesses with uniformly bounded $y$, $D_1$, $c_1$ and $c_2$. Then the attainable terminal set is weakly sequentially closed and hence norm closed. If the two relations are convex, this set is also convex and weakly closed.
\end{proposition}
\begin{proof}
For a weakly convergent terminal sequence, $D$ and $z$ are norm bounded. Hence $D_2=D-D_1$ is bounded. The threshold inequalities give $0\preceq m_1\preceq b$ and $0\preceq m_2\preceq b-c_1$, so the lattice norm and boundedness of $c_1$ bound both thresholds. Thus every hidden witness variable lies in a bounded subset of a reflexive Banach space. Weak sequential compactness of bounded balls gives a common subsequence converging weakly in this finite product. The two full relations pass to the limit by hypothesis. The positive cone is norm closed and convex, hence weakly closed, so the linear order inequalities pass to the limit as well. The response and reserve equations are bounded linear relations and also pass to the limit. This yields a valid limit witness. Norm-convergent sequences are weakly convergent, giving norm closedness. With convex elementary relations, the affine feasibility equations and inequalities give convexity after projection; a norm-closed convex set is weakly closed.
\end{proof}
For the diagonal map $S_0$ in Remark~\ref{m:rem:historyboundaries}(b), the composition of $\{0\}\times\ell^2$ with the graph of $S_0$ is $\{0\}\times\Ran S_0$. The truncations of $(1/n)_{n\ge1}$ converge in $\ell^2$, while their unique intermediate witnesses have norm $\sqrt N$ at cutoff $N$.

\subsection{Cost-whitened compensation and history-supported irreducible response}\label{m:sec:historybridge}
A historical direction can produce a response outside the old output closure even at arbitrarily small positive budget. The finite-dimensional history decomposition of \citet{paper1} extends to this setting through orthogonal projection and range factorization \cite{douglas,anderson}.

\begin{theorem}[History-supported irreducible compensation bridge]\label{m:thm:historybridge}
In the synchronous chart of Theorem~\ref{m:thm:weighted}, the cost-whitened old and acquired channels are
\[
 S_0=H\iota_0 M_0^{-1/2}:K_0\longrightarrow\calO,
 \qquad \widetilde{\mathcal L}=\Lambda W_{\calE}^{-1/2}:\calE\longrightarrow\calO.
\]
Their joint Gram is $S_0S_0^*+\widetilde{\mathcal L}\widetilde{\mathcal L}^*$. Set $\Pi_0=P_{\overline{\Ran S_0}}$, $\calO_\perp=\ker\Pi_0$ and $N_{\rm irr}=(I-\Pi_0)\widetilde{\mathcal L}$. Let $\mathcal H$ be a real Hilbert space carrying a positive trace-class operator $R$. With $J_h:\mathcal H\to\calE$ a bounded map into the cost-whitened acquisition coordinates, put $Z_h=\widetilde{\mathcal L}J_hR^{1/2}$ and
\begin{equation}\label{m:eq:historicalbridge}
 \boxed{\mathcal B_{\rm hist}=Z_h^*(I-\Pi_0)Z_h
 =R^{1/2}J_h^*N_{\rm irr}^*N_{\rm irr}J_hR^{1/2}.}
\end{equation}
\begin{enumerate}[label=(\roman*)]
\item $\mathcal B_{\rm hist}$ is positive and trace class. Its trace equals the best history-weighted squared output error:
\begin{equation}\label{m:eq:historicalinf}
 \begin{split}
 \Tr\mathcal B_{\rm hist}
 &=\|(I-\Pi_0)Z_h\|_{\mathfrak S_2}^2
 =\min_{A\in\mathfrak S_2(\mathcal H,\overline{\Ran S_0})}\|Z_h-A\|_{\mathfrak S_2}^2\\
 &=\inf_{K\in\mathfrak S_2(\mathcal H,K_0)}\|Z_h-S_0K\|_{\mathfrak S_2}^2.
 \end{split}
\end{equation}
The minimum in the middle expression is attained at $A_*=\Pi_0Z_h$. The final infimum is attained if and only if $\Pi_0Z_h$ admits a Hilbert--Schmidt factorization through $S_0$; in particular this holds when $\Ran S_0$ is closed.
\item Let $\lambda_1\ge\lambda_2\ge\cdots\ge0$ be the eigenvalues of $\mathcal B_{\rm hist}$ with multiplicity, completed by zeros as needed. For every integer $m\ge0$,
\begin{equation}\label{m:eq:historicaltail}
 \inf_{\substack{K\in\mathfrak S_2(\mathcal H,K_0)\\
                M\in\mathfrak S_2(\mathcal H,\calO_\perp),\ \operatorname{rank}M\le m}}
 \|Z_h-S_0K-M\|_{\mathfrak S_2}^2=\sum_{j>m}\lambda_j.
\end{equation}
Singular-value truncation of $(I-\Pi_0)Z_h$ gives an optimal orthogonal component $M$. A minimizing pair in~\eqref{m:eq:historicaltail} exists if and only if $\Pi_0 Z_h=S_0K$ for some $K\in\mathfrak S_2(\mathcal H,K_0)$. Allowing $M$ of rank at most $m$ with values anywhere in $\calO$ gives the same infimum.
\item $\mathcal B_{\rm hist}=0$ exactly when $N_{\rm irr}J_hR^{1/2}=0$, equivalently $\Ran Z_h\subseteq\overline{\Ran S_0}$. For dense $\Ran(J_hR^{1/2})$ in $\calE$ this is equivalent to $N_{\rm irr}=0$.
\item Suppose the full linear compensated chart is admissible at each quadratic budget. If $\Tr\mathcal B_{\rm hist}>0$, there exists $y=J_hR^{1/2}h$ with $N_{\rm irr}y\ne0$. At each $b>0$ the compensated channel reaches the target
\[
 d_b=\frac{\sqrt b\,\widetilde{\mathcal L}y}{\|y\|},\qquad
 \operatorname{dist}(d_b,\overline{\Ran S_0})
 =\frac{\sqrt b\,\|N_{\rm irr}y\|}{\|y\|}>0,
\]
at squared write cost $b$.
\end{enumerate}
\end{theorem}
\begin{proof}
Because $R^{1/2}$ is Hilbert--Schmidt and $\widetilde{\mathcal L},J_h$ are bounded, $Z_h\in\mathfrak S_2(\mathcal H,\calO)$. Put $F=(I-\Pi_0)Z_h$. Then $\mathcal B_{\rm hist}=F^*F\succeq0$ is trace class with trace $\|F\|_{\mathfrak S_2}^2$. For any Hilbert--Schmidt $A$ taking values in $\overline{\Ran S_0}$, the Pythagorean identity in the Hilbert--Schmidt space gives
\[
 \|Z_h-A\|_{\mathfrak S_2}^2
 =\|F\|_{\mathfrak S_2}^2+\|\Pi_0Z_h-A\|_{\mathfrak S_2}^2.
\]
Hence $A_*=\Pi_0Z_h$ is the unique minimum in the closed old-output sector. The same lower bound applies to $A=S_0K$. To approximate it through the old operator, for $\alpha>0$ take
\[
 K_\alpha=S_0^*(S_0S_0^*+\alpha I)^{-1}Z_h\in\mathfrak S_2(\mathcal H,K_0).
\]
Then $S_0K_\alpha=S_0S_0^*(S_0S_0^*+\alpha I)^{-1}Z_h\to\Pi_0Z_h$ in Hilbert--Schmidt norm: the positive contractions converge strongly to $\Pi_0$, and convergence on a finite-rank approximation to $Z_h$ plus the uniform bound proves the claimed ideal-norm limit. An exact old factor achieves the infimum if and only if $S_0K=\Pi_0Z_h$. Closed $\Ran S_0$ gives a bounded Moore--Penrose inverse, hence the Hilbert--Schmidt factor $K=S_0^\dagger\Pi_0Z_h$.

For (ii), orthogonality splits the squared error into $\|\Pi_0 Z_h-S_0K\|_{\mathfrak S_2}^2+\|F-M\|_{\mathfrak S_2}^2$. The first term has infimum zero by (i). Truncating the singular values of $F$ gives the second term minimum $\sum_{j>m}\lambda_j$, since $F^*F=\mathcal B_{\rm hist}$. For an arbitrary output component $M$ of rank at most $m$, projection of the error onto $\calO_\perp$ gives the lower bound $\|F-(I-\Pi_0)M\|_{\mathfrak S_2}^2\ge\sum_{j>m}\lambda_j$. Orthogonal components already approach that bound jointly with $S_0K_\alpha$.

For (iii), $F^*F=0$ precisely when $F=0$; a bounded map vanishing on a dense subspace vanishes everywhere. For (iv), choose $h$ with $Fh\ne0$. Theorem~\ref{m:thm:weighted} gives the cost $\|y\|^2$ for the physical compensated write with acquisition coordinate $W_{\calE}^{-1/2}y$, where $y=J_hR^{1/2}h$. Scaling this input by $\sqrt b/\|y\|$ gives $d_b$ at cost $b$. Its projection onto $\calO_\perp$ is $\sqrt b\,Fh/\|y\|$, giving the stated distance.
\end{proof}

\begin{figure}[tbp]
\centering
\includegraphics[trim=4pt 4pt 4pt 4pt,clip]{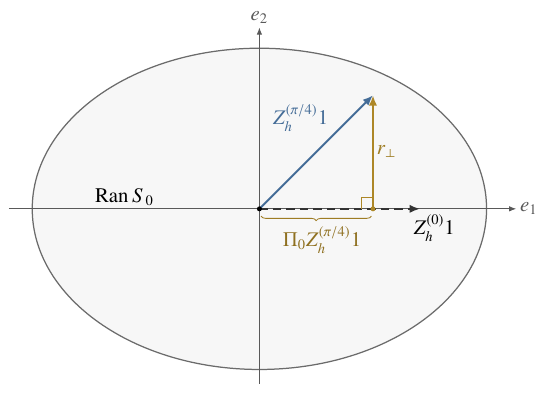}
\caption{Two historical directions in the same unit-budget response ellipse. The old channel is $S_0s=se_1$, the acquired channel $\widetilde{\mathcal L}_\theta$ rotates through $\theta$, and the unit historical input gives $Z_h^{(\theta)}1=\widetilde{\mathcal L}_\theta e_1$. Both angles give $G_{\rm joint}=\operatorname{diag}(2,1)$. For $\theta=0$ the dashed vector lies in $\Ran S_0$. For $\theta=\pi/4$ the solid vector splits into its orthogonal projection $\Pi_0 Z_h^{(\pi/4)}1$ and the residual $r_\perp$. The squared residual is $\sin^2\theta$, the history-weighted squared approximation error. On the scalable linear chart, every squared budget $b>0$ affords $\sqrt b\,\widetilde{\mathcal L}_{\pi/4}e_1$ outside the old range (Theorem~\ref{m:thm:historybridge}).}
\label{fig:history-bridge}
\end{figure}

Figure~\ref{fig:history-bridge} fixes the joint Gram and rotates the history-supported response relative to the old range.

\begin{remark}[Strict separation of historical residual, exact range, and costed Gram]\label{m:rem:historyboundaries}
The following three examples distinguish the information carried by the historical residual, the exact old response range, and the costed Gram.
\begin{enumerate}[label=(\alph*)]
\item \emph{Cost saving at zero historical residual.} With $S_0=\widetilde{\mathcal L}=1:\R\to\R$ and $J_h=R=1$, the old output spans all of $\R$, so $\mathcal B_{\rm hist}=0$, yet $G_0=1$, $G_{\rm joint}=2$. The minimum squared cost of output $d=1$ decreases from $1$ to $1/2$.
\item \emph{Exact range extension inside the old closure.} Let $S_0:\ell^2\to\ell^2$ be $(S_0x)_n=x_n/n$ and take $\widetilde{\mathcal L}:\R\to\ell^2$ with $\widetilde{\mathcal L}(1)=z=(1/n)_{n\ge1}$. Set $\mathcal H=\R$, $J_h=R=1$. Since $\overline{\Ran S_0}=\ell^2$, one has $\Pi_0=I$ and $\mathcal B_{\rm hist}=0$. But $z\notin\Ran S_0$ (its unique formal preimage is the non-square-summable sequence $(1,1,\ldots)$), whereas the new channel reaches $z$ at unit squared cost. For $m=1$, an unrestricted completion $M=Z_h$ and $K=0$ attain zero error. In~\eqref{m:eq:historicaltail}, $\calO_\perp=\{0\}$ and the same infimum is approached by truncations of $z$.
\item \emph{Equal present Grams with different historical residuals.} Set $\calO=\calE=\R^2$, $S_0t=te_1$, $J_h1=e_1$, $R=1$, and compare the acquisition channels $\widetilde{\mathcal L}_1=I_2$ and $\widetilde{\mathcal L}_2$ equal to a quarter-turn. Both give $\widetilde{\mathcal L}_i\widetilde{\mathcal L}_i^*=I_2$ and the same $G_{\rm joint}$, while $\mathcal B_{{\rm hist},1}=0$ and $\mathcal B_{{\rm hist},2}=1$.
\end{enumerate}

\end{remark}

\paragraph{From response geometry to execution.}
At a fixed parent $Z/\ker\varpi$ identifies cooperative response classes and the Gram determines the linear budget image. Across successive parents, controlled congruence identifies the state distinctions needed to reproduce labelled transitions. A physical decoder and a finite candidate menu specify the writes available to a controller.
\section{Physical decoders, tolerances and finite menus}\label{e:sec:writable}
At the execution state $(q,b)$ of Section~\ref{m:sec:protocols}, the decoder fixes a space of physical writes. Its protected part determines the local response Gram. Finite tolerance permits normal motion, while a proposal menu restricts the directions and amplitudes available within the budget.

\subsection{Physical decoders and protected capacity}\label{e:sec:decoder}
\begin{theorem}[Interface-faithful protected Gram order]\label{e:thm:decoder}
Let $\mathcal Z,\mathcal U,\mathcal O$ be real Hilbert spaces, $D_{\rm phys}:\mathcal Z\to\mathcal U$ a bounded linear decoder with closed range, and $H:\mathcal U\to\mathcal O$ a bounded linear response. Let $M$ be a bounded self-adjoint operator on $\mathcal U$ with $M\succeq mI$ for some $m>0$, and write $\|u\|_M^2=\langle u,Mu\rangle$. For a closed protection subspace $\mathcal K\subseteq\mathcal U$, the decoded protected subspace $\mathcal V=\operatorname{Ran}D_{\rm phys}\cap\mathcal K$ is closed in $\mathcal K$. At squared budget $b\ge0$, take the action domain in each $S\in\{\mathcal V,\mathcal K\}$ to be $S\cap\{u:\|u\|_M^2\le b\}$. Write $P_{\widetilde S}$ for the orthogonal projector onto $\widetilde S=M^{1/2}S$ and set
\[
G_S=H M^{-1/2}P_{\widetilde S}M^{-1/2}H^*.
\]
The response set and its support function in direction $g\in\mathcal O$ are
\[
\begin{aligned}
\mathcal A_S(b)&=\{Hu:u\in S,\ \langle u,Mu\rangle\le b\},\\
h_{\mathcal A_S(b)}(g)&=\sup_{y\in\mathcal A_S(b)}\langle g,y\rangle
=\sqrt b\sqrt{\langle g,G_Sg\rangle},
\end{aligned}
\]
and
\[
\boxed{0\preceq G_{{\mathcal V}}\preceq G_{{\mathcal K}}}.
\]
The two Grams coincide precisely when
\[
\big(P_{\widetilde{\mathcal K}}-P_{\widetilde{\mathcal V}}\big)M^{-1/2}H^*=0.
\]
\end{theorem}
\begin{proof}
Coercivity makes $M^{1/2}$ a bounded isomorphism, so $\widetilde{\mathcal V}\subseteq\widetilde{\mathcal K}$ are closed subspaces. Set $w=M^{1/2}u$. The restriction $u\in S$ is equivalent to $w\in\widetilde S$, and the cost bound is $\|w\|\le\sqrt b$. Maximizing $\langle g,HM^{-1/2}w\rangle$ over this ball gives $\sqrt b\|P_{\widetilde S}M^{-1/2}H^*g\|$. By nesting, $P_{\widetilde{\mathcal K}}-P_{\widetilde{\mathcal V}}$ is the orthogonal projector onto $\widetilde{\mathcal K}\cap\widetilde{\mathcal V}^{\perp}$. Consequently $\langle g,(G_{\mathcal K}-G_{\mathcal V})g\rangle=\|(P_{\widetilde{\mathcal K}}-P_{\widetilde{\mathcal V}})M^{-1/2}H^*g\|^2$. This proves both the order and the equality criterion.
\end{proof}

\paragraph{Physical metric calibration.}
If a latent coefficient $z\in\mathcal Z$ decodes as $u=D_{\rm phys}z$, then its induced \emph{physical} cost is
\[
\|z\|_{D_{\rm phys}^* M D_{\rm phys}}^2=\langle z,D_{\rm phys}^* M D_{\rm phys}z\rangle=\|D_{\rm phys}z\|_M^2.
\]
The pullback seminorm vanishes exactly on $\ker D_{\rm phys}$. The induced map $[z]\mapsto D_{\rm phys}z$ is an isometry from $\mathcal Z/\ker D_{\rm phys}$ in the pullback norm onto $\operatorname{Ran}D_{\rm phys}$ in the $M$-norm. This range is complete because it is closed and the $M$-norm is equivalent to the original norm on $\mathcal U$.

For $\mathcal U=\mathcal K=\mathbb R^2$, $M=I$ and $H(x_1,x_2)=x_2$, a decoder with range $\operatorname{span}\{e_1\}$ has $G_{\mathcal V}=0$ and $G_{\mathcal K}=1$. In the Euclidean specialization of Theorem~\ref{m:thm:weighted}, $G_+=NN^*+LW^{-1}L^*$ describes the full chart. Restriction to a physical source gives $G_{\mathcal V}$. Restriction to a finite set of proposals introduces the menu factor of Section~\ref{e:sec:menu}.

\subsection{Toleranced opportunity and paid installation}\label{e:sec:tolerance}

\begin{theorem}[Opportunity under a common protection tolerance]\label{e:thm:tolerance}
In physically whitened Euclidean coordinates, let $p_i,r_i$ be the protection and response maps of chart $i$, both zero at the origin, with $C_i=Dp_i(0)$, $H_i=Dr_i(0)$ and $P_i=P_{\ker C_i}$. Write $\Delta_i^\tau(\rho)$ for the supremum of $\ell(0)-\ell(r_i(u))$ over the chart's admissible endpoints $\|u\|\le\rho$, $\|p_i(u)\|\le\tau$, including zero. Let $g=\nabla\ell(0)$ and $\|g-\widehat g\|\le\varepsilon$. Define
\[
 \kappa_i(s)=\|P_iH_i^*s\|,\qquad \nu_i(s)=\|(C_i^\dagger)^*H_i^*s\|,\qquad
 h_i=\|H_iP_i\|,\qquad n_i=\|H_iC_i^\dagger\|.
\]

Suppose $p_0,r_0$ are $C^{1,1}$ on a neighborhood of the $R$-ball and throughout that ball
\[
 \|p_0(u)-C_0u\|\le\tfrac12 M_p\|u\|^2,\qquad
 |\ell(r_0(u))-\ell(0)-\langle g,H_0u\rangle|\le E_0\|u\|^2.
\]
For $0\le\rho\le R$, $\tau\ge0$, put $F(k,v)=k\rho+v\tau+(E_0+\tfrac12 M_pv)\rho^2$. Then
\begin{equation}\label{e:eq:tolerance-est}
 \Delta_0^\tau(\rho)\le F(\kappa_0(g),\nu_0(g))
 \le U_0(\rho,\tau):=F(\kappa_0(\widehat g)+\varepsilon h_0,\nu_0(\widehat g)+\varepsilon n_0).
\end{equation}

\end{theorem}
\begin{proof}
Write $u=P_0u+C_0^\dagger[p_0(u)-(p_0(u)-C_0u)]$. Pairing with $-H_0^*g$ bounds the tangent contribution by $\kappa_0(g)\rho$ and the normal contribution by $\nu_0(g)(\tau+M_p\rho^2/2)$. Add the objective remainder and take the supremum. The functions $\kappa_i$ and $\nu_i$ are Lipschitz with constants $h_i$ and $n_i$, yielding~\eqref{e:eq:tolerance-est}.
\end{proof}

Theorem~\ref{m:thm:commission} compares the two exact-repair charts at a common parent and objective, with protection and baseline preserved by installation. For common constants $K,\rho_*$, installation charge $c$ and squared-write budget $B$ with $0\le c\le B\le\rho_*^2$, substituting $|\kappa_i(g)-\kappa_i(\widehat g)|\le\varepsilon h_i$ gives, with $[z]_+=\max\{z,0\}$,
\begin{equation}\label{e:eq:robust}
\begin{aligned}
 \Delta_1^0(\sqrt{B-c})-\Delta_0^0(\sqrt B)\ge{}&
 [\kappa_1(\widehat g)-\varepsilon h_1]_+\sqrt{B-c}\\
 &-(\kappa_0(\widehat g)+\varepsilon h_0)\sqrt B-K(2B-c).
\end{aligned}
\end{equation}

\begin{corollary}[One installed candidate]\label{e:prop:tolerance-candidate}
Let the old chart satisfy Theorem~\ref{e:thm:tolerance}. At the common parent, let installation preserve protection and the objective baseline and consume $c$ of the squared-write budget $B$, where $0\le c\le B\le R^2$. Let $u_1$ be an admissible installed action with $\|u_1\|^2\le B-c$, $\|p_1(u_1)\|\le\tau$, and objective remainder bounded by $E_1\|u_1\|^2$ at that action. Let the modeled old domain at $(\sqrt B,\tau)$ contain the native old comparison set, and let $\delta_0,\delta_1$ bound the errors in modeled decrease on that set and at $u_1$. The candidate's native decrease strictly exceeds the old native supremum whenever
\begin{equation}\label{e:eq:tolerance-certificate}
 -\langle\widehat g,H_1u_1\rangle-\varepsilon\|H_1u_1\|-E_1\|u_1\|^2-\delta_1
 > U_0(\sqrt B,\tau)+\delta_0 .
\end{equation}
\end{corollary}
\begin{proof}
Cauchy--Schwarz bounds the gradient error by $\varepsilon\|H_1u_1\|$. Subtracting this term, the objective remainder and $\delta_1$ gives the lower bound on the candidate's native decrease in~\eqref{e:eq:tolerance-certificate}. The modeled old domain contains every compared old write, so Theorem~\ref{e:thm:tolerance} bounds the old native supremum by $U_0(\sqrt B,\tau)+\delta_0$.
\end{proof}
For exact protection, the old chart's regular-repair hypotheses and $\sqrt B\le\rho_*$ allow $(\kappa_0(\widehat g)+\varepsilon h_0)\sqrt B+KB$ in place of $U_0$.

If both measured guards use $\|\widehat p_i\|\le\tau_{\rm meas}$ with errors $\|p_i-\widehat p_i\|\le d_i$, set $\tau=\tau_{\rm meas}+\max_i d_i$ to cover both native residuals. A common hard ceiling $\tau_{\rm safe}$ holds on the joint error event when $\tau\le\tau_{\rm safe}$.

\subsection{Finite candidate menus and realized opportunity}\label{e:sec:menu}
A generator returns a finite subset of the physical write space. Its best linearized decrease is the support function of that subset in the loss-gradient direction.

\begin{figure}[tbp]
\centering
\includegraphics[trim=4pt 4pt 4pt 4pt,clip]{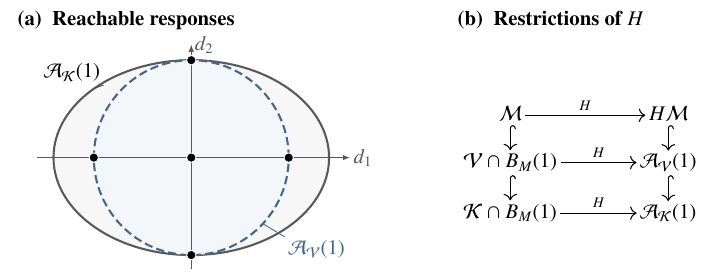}
\caption{A physical decoder and a finite menu restrict the response at one parent. The ambient action space is $\mathcal K=\mathbb R^3$ with $M=I$, the decoder writes $\mathcal V=\operatorname{span}\{e_1,e_2\}$, and $H(u_1,u_2,u_3)=(u_1+u_3,u_2)$. (a) At unit squared budget the ambient response image is the ellipse $d_1^2/2+d_2^2\le1$. The decoder gives the disk $d_1^2+d_2^2\le1$ with dashed boundary. The menu $\mathcal M=\{0,\pm e_1,\pm e_2\}$ gives the five marked points. (b) Both squares commute. The vertical arrows are inclusions and the horizontal arrows apply $H$ on the indicated domain.}
\label{fig:response-capacity}
\end{figure}

Figure~\ref{fig:response-capacity} separates the restriction to a physical source from the subsequent restriction to a finite menu.

\begin{theorem}[Exact budget-weighted finite-menu alignment]\label{e:thm:menu}
Use $M,H,\mathcal K,\mathcal V$ of Theorem~\ref{e:thm:decoder}. Fix $b>0$ and a finite executable menu $\mathcal M\subseteq{\mathcal V}$ with zero action included and $\|u\|_M^2\le b$ for every $u\in\mathcal M$. Let $g\in{\mathcal O}$ be the true gradient of one fixed current loss at the common parent. Put
\[
a_{\mathcal V}=P_{M^{1/2}{\mathcal V}}M^{-1/2}H^*g,\quad
\kappa_{\mathcal V}=\|a_{\mathcal V}\|=\sqrt{\langle g,G_{{\mathcal V}}g\rangle}.
\]
If $\kappa_{\mathcal V}>0$, set $v_*=-a_{\mathcal V}/\kappa_{\mathcal V}$ and define
\[
\gamma_\mathcal M(g;b):=
\max_{u\in\mathcal M}\left\langle\frac{M^{1/2}u}{\sqrt b},v_*\right\rangle\in[0,1].
\]
Then the menu's optimal linearized loss decrease has the \emph{exact} factorization
\begin{equation}\label{e:eq:menu-linear}
\boxed{\mathsf A_\mathcal M^{\rm lin}(g;b):=
\max_{u\in\mathcal M}\{-\langle g,Hu\rangle\}
=\gamma_\mathcal M(g;b)\sqrt{b\langle g,G_{{\mathcal V}}g\rangle}.}
\end{equation}
If $\kappa_{\mathcal V}=0$, the linearized menu opportunity is zero and we set $\gamma_\mathcal M=0$. When $\kappa_{\mathcal K}:=\sqrt{\langle g,G_{{\mathcal K}}g\rangle}>0$, the ideal-to-realized relative linear opportunity factors as
\begin{equation}\label{e:eq:menu-ratio}
\boxed{\frac{\mathsf A_\mathcal M^{\rm lin}(g;b)}{\sqrt b\,\kappa_{\mathcal K}}
=\frac{\kappa_{\mathcal V}}{\kappa_{\mathcal K}}\gamma_\mathcal M(g;b).}
\end{equation}
The first factor gives the fraction retained by the decoder. The second measures the signed directional alignment of each proposal weighted by the fraction of the radius it uses, maximized over the menu.
\end{theorem}
\begin{proof}
For $u\in{\mathcal V}$, $\langle g,Hu\rangle=\langle a_{\mathcal V},M^{1/2}u\rangle$. If $\kappa_{\mathcal V}>0$, negating and normalizing gives
$-\langle g,Hu\rangle=\sqrt b\,\kappa_{\mathcal V}\langle M^{1/2}u/\sqrt b,v_*\rangle$.
Maximizing over the finite menu proves \eqref{e:eq:menu-linear}. Zero belongs to the menu, so the maximum is nonnegative; Cauchy--Schwarz and $\|u\|_M\le\sqrt b$ bound it by one. Since ${\mathcal V}\subseteq{\mathcal K}$, Theorem~\ref{e:thm:decoder} gives $\kappa_{\mathcal V}\le\kappa_{\mathcal K}$; division yields \eqref{e:eq:menu-ratio}. If $\kappa_{\mathcal V}=0$, every menu action has zero linearized response against $g$.
\end{proof}

For each offered action let $D(u)$ and $D_{\rm native}(u)$ be the modeled and native decreases, with $|D(u)+\langle g,Hu\rangle|\le E_M\|u\|_M^2$, $|D_{\rm native}(u)-D(u)|\le\delta_{\rm exec}$ and $\|g-\widehat g\|\le\varepsilon$. Then $v(u)=-\langle\widehat g,Hu\rangle$ approximates the native decrease within $h(u)=\varepsilon\|Hu\|+E_M\|u\|_M^2+\delta_{\rm exec}$. Taking maxima gives
\begin{gather}
 \max_{u\in\mathcal M}\{v(u)-h(u)\}\le \max_{u\in\mathcal M}D_{\rm native}(u)
 \le\max_{u\in\mathcal M}\{v(u)+h(u)\},\label{e:eq:menu-observable-sharp}\\
 \big|\max_{u\in\mathcal M}D(u)-\mathsf A_{\mathcal M}^{\rm lin}(g;b)\big|\le E_M b.\label{e:eq:menu-finite}
\end{gather}
If selection uses the gradient estimate, the same inequalities hold on any event giving the bounds simultaneously over all eligible actions.

The factorization~\eqref{e:eq:menu-ratio} concerns actions in the protected subspace $\mathcal V$. For a menu with tolerance $\|p(u)\|\le\tau$, Theorem~\ref{e:thm:tolerance} also accounts for the normal response.

\section{Continuation and policy comparison}\label{e:sec:continuation}
A current protected action changes the starting point of subsequent dynamics. To compare its future effect, we specify the entire continuation, a common semantic readout and an error certificate valid for the selected action. Deterministic error tubes bound the terminal readout along matched paths. Stochastic comparisons bound the expected terminal advantage using the successor laws on each signature sector.

\subsection{Matched continuation and current-to-terminal transfer}\label{e:sec:contract}
\begin{definition}[Matched finite continuation]\label{e:def:continuation}
Use the execution state of Section~\ref{m:sec:protocols}, with a metric on each compatible signature sector. At a common parent $x=(q,b)$, let $a$ be an action constructed from pre-decision information $\mathcal I_t$, and let $a_0$ be the reference action. Their post-action execution states are $x_a=\Phi_a(x)$ and $x_0=\Phi_{a_0}(x)$ when both transitions are defined. The reference no-write action still includes all subsequently prescribed ordinary updates. Fix one continuation rule $\pi$, horizon $H$, matched exogenous trace $\omega$ and a \emph{common} real-valued terminal readout $J$; define
\begin{equation}\label{e:eq:continuation-value}
V_H^{\pi,\omega}(x,a;J)=J\bigl(F_H^{\pi,\omega}(x_a)\bigr)-J\bigl(F_H^{\pi,\omega}(x_0)\bigr).
\end{equation}
Here $F_H^{\pi,\omega}$ is the complete closed-loop update from a post-action cut. Under the shared trace $\omega$, the continuation rule selects each subsequent action from the state reached on its branch. The complete state includes the controller memory, randomization, retained information, policy and structural data required by that update. Across distinct signatures, the terminal comparison uses the common readout $J$. At a protected checkpoint $h$, the corresponding comparison uses the common readout $P_h$.
\end{definition}

A current loss decrease is a current contrast for $J=-\ell\circ r$. A different terminal utility requires its own current contrast or a bound transferring the decrease to that readout.

\subsection{Deterministic and matched-path continuation bounds}\label{e:sec:continuation-bounds}
\begin{theorem}[Closed-loop error transport in nominal tubes]\label{e:thm:transport}
For each branch $i\in\{0,a\}$ and checkpoint $k=0,\ldots,H$, fix a metric space $(\mathcal S_k^i,d_k^i)$ for the compatible true and nominal signatures at that checkpoint. At the post-action cut, let $x_0^i=x_i$ be the actual state and $\widehat x_0^i$ the nominal state. Certify a finite initial radius $E_0^i\ge d_0^i(x_i,\widehat x_0^i)$. For $0\le k<H$ set
\[
x_{k+1}^i=F_k^i(x_k^i;\omega),\qquad
\widehat x_{k+1}^i=\widehat F_k^i(\widehat x_k^i;\omega).
\]
Here $F_k^i$ includes feedback and structural changes along the signature sequence of branch $i$. Given nonnegative bounds $L_k^i,e_k^i$, form the radii
\begin{equation}\label{e:eq:errors}
E_{k+1}^i=L_k^iE_k^i+e_k^i.
\end{equation}
For every $i$ and $0\le k<H$, suppose the nominal-centered closed ball
\[
\overline B_{E_k^i}^{d_k^i}(\widehat x_k^i)\subseteq \operatorname{Dom}F_k^i
\]
lies within the executable domain of the closed-loop update. Assume $F_k^i$ maps this ball into $\mathcal S_{k+1}^i$, is $L_k^i$-Lipschitz there from $d_k^i$ to $d_{k+1}^i$, and satisfies
\[
d_{k+1}^i\big(F_k^i(\widehat x_k^i;\omega),\widehat F_k^i(\widehat x_k^i;\omega)\big)\le e_k^i.
\]
If the common terminal semantic readout $J$ is $L_J^i$-Lipschitz on the corresponding terminal ball in branch $i$, and $|J(\widehat x_H^i)-\widehat J(\widehat x_H^i)|\le d_J^i$, then
\begin{equation}\label{e:eq:terminal-error}
\boxed{\left|V_H^{\pi,\omega}-\left[\widehat J(\widehat x_H^a)-\widehat J(\widehat x_H^0)\right]\right|
\le\eta_a+\eta_0,\quad \eta_i=L_J^iE_H^i+d_J^i.}
\end{equation}
For a random matched trace $\omega$, the bound holds on the event where the stated hypotheses hold jointly for both branches at every checkpoint.
\end{theorem}
\begin{proof}
The initial certificate puts $x_0^i$ in the first ball. Assuming $d_k^i(x_k^i,\widehat x_k^i)\le E_k^i$, the true state lies in the certified domain of $F_k^i$. The Lipschitz bound and the one-step mismatch in the common successor metric give
\[
d_{k+1}^i(x_{k+1}^i,\widehat x_{k+1}^i)\le L_k^i E_k^i+e_k^i=E_{k+1}^i.
\]
Induction puts each true endpoint in its terminal ball. Its readout error is at most $L_J^iE_H^i+d_J^i$. Subtracting the scalar readouts and applying the triangle inequality gives \eqref{e:eq:terminal-error}.
\end{proof}

On a common metric chart let $\varpi_H$ be a modulus of continuity for the continuation increment $R(x)=J(F_H^{\pi,\omega}(x))-J(x)$. A current contrast $J(x_a)-J(x_0)\ge d$ then gives $V_H^{\pi,\omega}\ge d-\varpi_H(d_*(x_a,x_0))$. The nominal tubes also give a direct estimate of the terminal contrast in~\eqref{e:eq:terminal-error}.

\subsection{Relative and absolute future protection}\label{e:sec:future}
\begin{figure}[tbp]
\centering
\includegraphics[trim=4pt 4pt 4pt 4pt,clip]{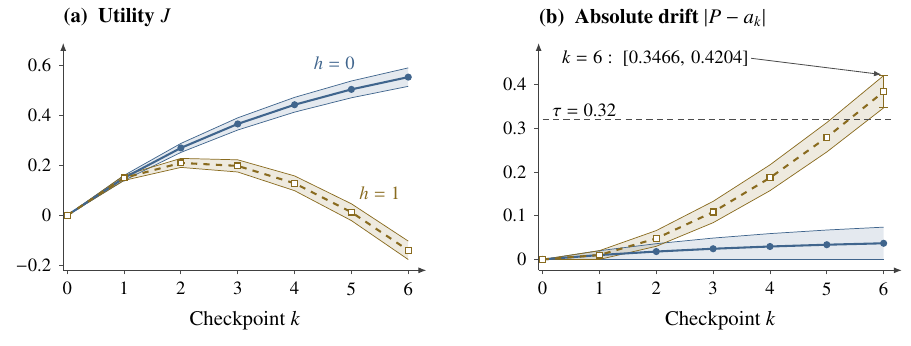}
\caption{A common continuation from two hidden states with identical current readouts. Both branches start at $w=(j,p)=(0,0)$ with fixed mode $h\in\{0,1\}$ and anchor $a_k=0$. A write $(u,v)$ changes $w$ by $(u,v)$ at squared cost $u^2+v^2$. For $J=j$ and $P=p$, current protection gives $v=0$ and both utility-response Grams equal one. The nominal transition is $\widehat F_k(w,h)=0.8w+(0.15-0.06kh,\,0.01+0.03kh)$. The true transition adds the same disturbance trace on both branches, with $\|\xi_k\|_\infty\le0.01$. The $0.8$-Lipschitz bound in the state sup norm gives $E_0=0$ and $E_{k+1}=0.8E_k+0.01$. At integer checkpoints the shaded intervals bound the true readouts (Theorems~\ref{e:thm:transport} and~\ref{e:thm:protect}). Circles mark $h=0$, squares $h=1$. The $h=0$ drift stays below $0.074$. For $h=1$ the upper bound at $k=5$ is $0.312992$. At $k=6$ the nominal drift is $0.3835008$ and even the lower bound, $0.346608$, exceeds $\tau=0.32$.}
\label{fig:continuation-tubes}
\end{figure}

\begin{theorem}[Relative and frozen-anchor checkpointed protection envelopes]\label{e:thm:protect}
Continue under Theorem~\ref{e:thm:transport}. At a checkpoint $0\le h\le H$, let $P_h$ be a common semantic protected readout into a normed space $\mathcal Y_h$ with a realization on each branch. On branch $i$, assume that $P_h$ is $L_{P,h}^i$-Lipschitz on $\overline B_{E_h^i}^{d_h^i}(\widehat x_h^i)$ and that the nominal readout $\widehat P_h^i$ satisfies
$\|P_h(\widehat x_h^i)-\widehat P_h^i(\widehat x_h^i)\|\le d_{P,h}^i$.
Then the \emph{relative} branch comparison satisfies
\begin{equation}\label{e:eq:guard}
\begin{split}
\|P_h(x_h^a)-P_h(x_h^0)\|\le{}&
\|\widehat P_h^a(\widehat x_h^a)-\widehat P_h^0(\widehat x_h^0)\|\\
&+L_{P,h}^aE_h^a+L_{P,h}^0E_h^0+d_{P,h}^a+d_{P,h}^0.
\end{split}
\end{equation}
Independently, let $a_h(x)\in\mathcal Y_h$ denote the \emph{frozen-parent historical obligation} transported to the checkpoint's output semantics by a fixed map. For a computed $\widehat a_h(x)$ with $\|a_h(x)-\widehat a_h(x)\|\le d_{\mathrm{anc},h}$, the action branch obeys the \emph{absolute} guard
\begin{equation}\label{e:eq:absolute-guard}
\boxed{\begin{aligned}
\|P_h(x_h^a)-a_h(x)\|\le{}&
\|\widehat P_h^a(\widehat x_h^a)-\widehat a_h(x)\|\\
&+L_{P,h}^aE_h^a+d_{P,h}^a+d_{\mathrm{anc},h}.
\end{aligned}}
\end{equation}
The absolute bound for the reference branch follows by setting the branch index to $0$ in \eqref{e:eq:absolute-guard}.

\end{theorem}
\begin{proof}
For each branch insert its nominal state and implemented nominal protected readout; this bounds its semantic readout error by $L_{P,h}^iE_h^i+d_{P,h}^i$. Taking the difference of the two branches gives \eqref{e:eq:guard}. For \eqref{e:eq:absolute-guard}, insert the nominal action readout and the independently estimated frozen anchor; three applications of the triangle inequality give the claim. Both estimates compare readouts in $\mathcal Y_h$.
\end{proof}
Figure~\ref{fig:continuation-tubes} follows two states with equal current readouts and equal protected response Grams through one closed-loop continuation.

A continuation can carry a hidden displacement into a later protected coordinate even when the current protected readouts agree. Equality of the two later readouts also permits a common displacement from the frozen anchor, which enters the absolute bound.

\subsection{Stochastic continuation and policy comparison}\label{e:sec:stochastic}
Assume the execution-state space is standard Borel with measurable signature sectors. Each sector used for transport carries a compatible Polish metric. Fix a finite horizon $T$. Write $Q_t^-$ and $Q_t^+$ for the state--reserve pairs immediately before certification and at the ensuing decision. Let $\mathscr F_t^-\subseteq\mathscr F_t^+$ be the corresponding histories, including the execution states. The selector observes $\mathcal G_t\subseteq\mathscr F_t^+$, including its random draws up to the decision.

Under a policy $\pi$, the prefix kernel $C_t^\pi$ describes certification and its charges, and the execution kernel $K_t^{D_t}$ describes the selected branch $D_t$. On the reached histories these satisfy
\[
 \mathcal L_\pi(Q_t^+\mid\mathscr F_t^-)=C_t^\pi(Q_t^-,\cdot),\qquad
 \mathcal L_\pi(Q_{t+1}^-\mid\mathscr F_t^+)=K_t^{D_t}(Q_t^+,\cdot).
\]
The state includes the history required for these identities. A stopped path enters an absorbing state carrying its payoff at zero subsequent cost. Under the reference response $K_t^\beta$, backward recursion~\cite{mdp} gives
\[
 V_T^{\beta,-}=J_T,\qquad V_t^{\beta,+}=K_t^\beta V_{t+1}^{\beta,-},\qquad
 V_t^{\beta,-}=C_t^\beta V_t^{\beta,+}.
\]

\begin{theorem}[Two-cut performance difference]\label{e:thm:two-cut}
Let $\alpha$ and $\beta$ share an initial law, an integrable terminal utility $J_T$ and the stopping convention. Suppose the reference kernels and continuation values above are defined on the states needed along both policies, and all displayed expectations are finite. Writing $\mathcal J(\pi)=\mathbb E_\pi J_T$, one has
\begin{equation}\label{e:eq:two-cut}
\begin{aligned}
\mathcal J(\alpha)-\mathcal J(\beta)
={}&\sum_{t=0}^{T-1}\mathbb E_\alpha\!\left[\big[(C_t^\alpha-C_t^\beta)V_t^{\beta,+}\big](Q_t^-)\right]\\
 &+\sum_{t=0}^{T-1}\mathbb E_\alpha\!\left[\big[(K_t^{D_t}-K_t^\beta)V_{t+1}^{\beta,-}\big](Q_t^+)\right].
\end{aligned}
\end{equation}
The two terms value the certification prefix at the pre-certification state and the selected action at the ensuing decision state.
\end{theorem}
\begin{proof}
Telescope the reference values along an $\alpha$ trajectory, splitting each increment at $Q_t^+$. Conditional expectation of $V_t^{\beta,+}(Q_t^+)-V_t^{\beta,-}(Q_t^-)$ given $\mathscr F_t^-$ gives the first integrand. Conditional expectation of $V_{t+1}^{\beta,-}(Q_{t+1}^-)-V_t^{\beta,+}(Q_t^+)$ given $\mathscr F_t^+$ gives the second. The common initial law gives $\mathbb E_\alpha V_0^{\beta,-}(Q_0^-)=\mathcal J(\beta)$, and the terminal reference value is $J_T$.
\end{proof}

\begin{theorem}[A branch-metric one-step advantage certificate]\label{e:thm:advantage}
At a reached post-certification execution state $x$, let $a$ be an admissible action and let the reference response be defined on both successor laws and their reference continuations. Define
\[
 A_t^\beta(x,a)=\big[(K_t^a-K_t^\beta)V_{t+1}^{\beta,-}\big](x),\qquad
 \widehat A_t(x,a)=\big[(\widehat K_t^a-\widehat K_t^\beta)\widehat V_{t+1}^{\beta,-}\big](x).
\]
For each $i\in\{a,\beta\}$, assume the true and approximate successor laws at $x$ have finite first moments on a common signature-local Polish space $(\mathcal S_i,d_i)$, where $\widehat V_{t+1}^{\beta,-}$ is $L_{t+1}^i$-Lipschitz. Suppose
\[
 |V_{t+1}^{\beta,-}-\widehat V_{t+1}^{\beta,-}|\le\eta_{t+1}^i
 \quad K_t^i(x,\cdot)\text{-almost surely},\qquad
 W_1^{d_i}(K_t^i(x,\cdot),\widehat K_t^i(x,\cdot))\le\delta_t^i.
\]
Then
\begin{equation}\label{e:eq:adv-bound}
 |A_t^\beta(x,a)-\widehat A_t(x,a)|\le
 \eta_{t+1}^a+\eta_{t+1}^\beta+L_{t+1}^a\delta_t^a+L_{t+1}^\beta\delta_t^\beta.
\end{equation}
\end{theorem}
\begin{proof}
On each branch the true kernel integrates $|V_{t+1}^{\beta,-}-\widehat V_{t+1}^{\beta,-}|$ to at most $\eta_{t+1}^i$. The Kantorovich--Rubinstein bound on the change of kernel applied to $\widehat V_{t+1}^{\beta,-}$ contributes at most $L_{t+1}^i\delta_t^i$~\cite{kr}. Add the two branchwise bounds.
\end{proof}

Both branches use the reference value $V^\beta$ on their respective state spaces. For a successor law spanning several signatures, one may bound the kernel term on the full mixture by total variation when the readout is bounded, or by coupling the two laws under a transport cost that bounds readout variation across the participating strata.

\begin{corollary}[Admissible-policy improvement with a certified uncertainty margin]\label{e:cor:policy}
Let $\alpha$ and $\beta$ satisfy the two-cut kernel property above, with the same $C_t$, initial law, terminal utility and stopping semantics. Write $x_t=Q_t^+$. At each $\alpha$-reached decision, $\alpha$ either invokes the reference response or selects an admissible $\mathcal G_t$-measurable action $a_t$. Assume the hypotheses of Theorem~\ref{e:thm:advantage} hold simultaneously almost surely for the selected branches and the following margin is $\mathcal G_t$-measurable and nonnegative:
\begin{equation}\label{e:eq:policy-margin}
\boxed{\begin{aligned}
\underline A_t:={}&\widehat A_t(x_t,a_t)
-\bigl[\eta_{t+1}^{a_t}(x_t,a_t)+\eta_{t+1}^\beta(x_t,a_t)\bigr]\\
&{}-L_{t+1}^{a_t}(x_t)\delta_t^{a_t}(x_t)
-L_{t+1}^\beta(x_t)\delta_t^\beta(x_t)\ge0.
\end{aligned}}
\end{equation}
Set $\underline A_t=0$ on reference branches. If the required expectations are finite, then $\mathbb E_\alpha J_T\ge\mathbb E_\beta J_T$. Strictly positive expected summed margins imply strict improvement. Equal branchwise value-error bounds give the margin with $2\eta_{t+1}$.
\end{corollary}
\begin{proof}
With a common prefix kernel the first sum in~\eqref{e:eq:two-cut} vanishes. In the second sum the fallback has zero advantage, and~\eqref{e:eq:adv-bound} bounds the selected action's advantage below by $\underline A_t$. Summing gives the two assertions.
\end{proof}

\section{Paid certification, selection and resource contracts}\label{e:sec:resources}
Certification acts on the state and reserve before the execution decision. Both execution and abstention carry its cost forward. An adaptive policy can bound its risk by conditioning on the information available after selection or by controlling all candidates on a common event.

\subsection{Paid certification and abstention}\label{e:sec:paid}
Certification is an executable protocol. Its action generator, measurements, forward models, numerical checks, memory and discarded proposals use the same recorded resource units as the compared adaptive policies.

\begin{definition}[Paid gateway and two distinct reference policies]\label{e:def:paidgate}
At a pre-decision execution state $x=(q,b)$, the certification rule selects probes using the available information, structural signature $\varsigma(q)$ and reserve $b$. It produces a post-certification state and record $r_{\rm cert}=(m_{\rm cert},c_{\rm cert})$. The rule is defined at every reached decision invoking it and along each reference continuation used for valuation. Its domain includes the prescribed fallback or stopping transition when a required probe is unavailable.

The primary policy $\alpha_{\rm gate}$ selects its action after certification. The paid reference $\beta_{\rm gate}$ uses the same prefix kernel at each of its own states and then takes the reference action. The direct reference $\beta_{\rm direct}$ takes that action with an identity prefix. The execute and fallback suffixes have records $r_{\rm act}$ and $r_\beta$.
\end{definition}

\begin{figure}[tbp]
\centering
\includegraphics[trim=4pt 4pt 4pt 4pt,clip]{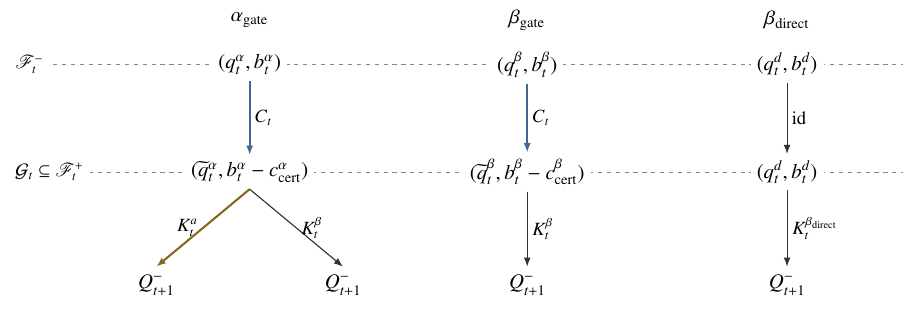}
\caption{Two decision cuts in a paid policy comparison. At a fixed decision time $t$, superscripts $\alpha,\beta,d$ identify the states and available reserves reached by the three policies. They share an initial law, terminal utility and stopping convention. The two gated policies apply the same certification kernel $C_t$ to their own states. The direct reference has an identity prefix. At the second cut the primary selector uses $\mathcal G_t$ to choose $K_t^a$ or the fallback $K_t^\beta$, both evaluated at the reached paid state. Each arrow to $Q^-_{t+1}$ denotes the next pre-certification state drawn from its edge kernel, including the remaining reserve. Subtracting the certification overhead from the paid-policy gain gives $\Delta_{\rm total}=\Delta_{\rm paid}-\Delta_{\rm overhead}$ (Corollary~\ref{e:thm:overhead}). The conditional risk bound for the paid comparison uses the selected action at the second cut (Theorem~\ref{e:thm:risk}).}
\label{fig:paid-gate}
\end{figure}

\begin{corollary}[Paid and direct reference comparisons]\label{e:thm:overhead}
Under Theorem~\ref{e:thm:two-cut}, suppose the three policies have a common initial law and integrable terminal utility. At a state where the relevant stages are executable, the complete decision records are
\begin{align}
 (m_{\rm exec},c_{\rm exec})&=\big(m_{\rm cert}\vee(c_{\rm cert}+m_{\rm act}),\ c_{\rm cert}+c_{\rm act}\big),\label{e:eq:exec-reserve}\\
 (m_{\rm abst},c_{\rm abst})&=\big(m_{\rm cert}\vee(c_{\rm cert}+m_\beta),\ c_{\rm cert}+c_\beta\big).\label{e:eq:abst-reserve}
\end{align}
Set $\Delta_{\rm paid}=\mathcal J(\alpha_{\rm gate})-\mathcal J(\beta_{\rm gate})$ and $\Delta_{\rm overhead}=\mathcal J(\beta_{\rm direct})-\mathcal J(\beta_{\rm gate})$. Then
\begin{equation}\label{e:eq:overhead-policy}
 \mathcal J(\alpha_{\rm gate})-\mathcal J(\beta_{\rm direct})
 =\Delta_{\rm paid}-\Delta_{\rm overhead}.
\end{equation}
Consequently, jointly valid bounds $L_{\rm paid}\le\Delta_{\rm paid}$ and $U_{\rm overhead}\ge\Delta_{\rm overhead}$ imply a gain over the direct reference of at least $L_{\rm paid}-U_{\rm overhead}$.
\end{corollary}
\begin{proof}
Proposition~\ref{m:prop:record} composes the certification record with each suffix at the intermediate reserve $b-c_{\rm cert}$. The shared prefix cancels in the paid comparison in~\eqref{e:eq:two-cut}. On a fallback the action advantage also vanishes, while the record retains the certification cost. Adding the contrast between paid and direct references gives~\eqref{e:eq:overhead-policy}.
\end{proof}

\begin{proposition}[Stopped adaptive search under a fixed contract]\label{e:prop:adaptive-search}
Starting at $(q,b)$ under $\mathfrak C_g$, let $N$ be a stopping time bounded by $N_{\max}$ for the trial filtration $(\mathcal I_n)_{n=0}^{N_{\max}}$, with $\mathcal I_0=\mathcal I_t$ and $\mathcal G_t=\mathcal I_N$. Each trial $n$ uses an $\mathcal I_{n-1}$-measurable candidate, follows an admissible path and incurs a reserve record. At each possible stop $n$, an $\mathcal I_n$-measurable rule selects a commit or fallback from the transcript. A commit requires the stored joint candidate to satisfy the full predicate $\mathcal Q_g$ of Definition~\ref{m:def:contract-promotion} and the physical write budget. Suppose the ordered records fit $b$ along every allowed search path, including its terminal branch and a fallback available at exhaustion.
Then the stopped decision is $\mathcal G_t$-measurable, its record is the ordered product of all incurred records, and each committed successor belongs to $\mathcal Q_g$.
\end{proposition}
\begin{proof}
On $\{N=n\}$ the terminal choice equals the $n$th measurable rule, giving $\mathcal I_N$-measurability. Proposition~\ref{m:prop:record} composes successive records at the remaining reserve. The fixed contract supplies the commit predicate.
\end{proof}
For fixed-parent trials, isolation or path-admissible restoration returns the physical and controller registers used by the candidate to their common values (Section~\ref{m:sec:protocols}). The transcript retains observations and each trial starts with the reserve left by earlier probes and restorations. A change in the registers defining the chart gives the next local model a new parent.

\subsection{Selection-safe probabilistic risk}\label{e:sec:risk}
For an advantage estimated from finite data, selection enters the risk calculation. The first result conditions the failure law on the information used by the selector. The second controls all selected actions on one common event.

\begin{theorem}[Post-selection conditional-risk admission and paid-reference improvement]\label{e:thm:risk}
Compare $\alpha_{\rm gate}$ and $\beta=\beta_{\rm gate}$ under the two-cut controlled kernels of Section~\ref{e:sec:stochastic}, with the same certification kernel $C_t$, stopping semantics, finite horizon, initial law and integrable terminal utility. Write $x_t=Q_t^+$, and let $\mathcal G_t$ be the information available after certification and selection and before execution. The indicator $\xi_t\in\{0,1\}$, admissible action $a_t$, signed lower margin $\underline\mu_t\in\mathbb R$, loss envelope $M_t\ge0$ and $0\le\zeta_t<1$ are $\mathcal G_t$-measurable. Set $a_t=\beta$ on abstention or after stopping, where $\xi_t=0$. Write $A_t^\beta:=A_t^\beta(x_t,a_t)=[(K_t^{a_t}-K_t^\beta)V_{t+1}^{\beta,-}](x_t)$ for the realized advantage. Assume all terms entering expectations are integrable under the reached primary law and reference values are defined on the required continuation domains. On $\{\xi_t=1\}$ suppose $A_t^\beta\ge-M_t$ always, $A_t^\beta\ge\underline\mu_t$ on $\mathcal E_t$, and
\begin{equation}\label{e:eq:postsel-risk}
{\mathbb P}(\mathcal E_t^c\mid\mathcal G_t)\le\zeta_t\qquad \text{a.s. on }\{\xi_t=1\}.
\end{equation}
On abstention the gateway invokes the reference post-certification kernel, and an absorbed path keeps its recorded payoff. Then, under $\alpha_{\rm gate}$,
\begin{equation}\label{e:eq:risk-lower}
\boxed{{\mathbb E}\left[\xi_t A_t^\beta\mid\mathcal G_t\right]
\ge \xi_t\big((1-\zeta_t)\underline\mu_t-\zeta_tM_t\big).}
\end{equation}
In particular, executing only when
\begin{equation}\label{e:eq:risk-threshold}
\boxed{\underline\mu_t\ge\frac{\zeta_t}{1-\zeta_t}M_t}
\end{equation}
makes the right-hand side nonnegative and gives $\mathcal J(\alpha_{\rm gate})\ge\mathcal J(\beta_{\rm gate})$. A strictly positive expected sum of the displayed lower bounds gives strict improvement over the paid reference. Corollary~\ref{e:thm:overhead} supplies the comparison with the direct reference.
\end{theorem}
\begin{proof}
On $\{\xi_t=1\}$, put $\mu_t^\star=\max\{\underline\mu_t,-M_t\}$ and $p_t=\mathbb P(\mathcal E_t^c\mid\mathcal G_t)$. The advantage is at least $\mu_t^\star$ on $\mathcal E_t$ and at least $-M_t$ on its complement. Hence its conditional mean is at least $\mu_t^\star-p_t(\mu_t^\star+M_t)$. Since $\mu_t^\star+M_t\ge0$ and $p_t\le\zeta_t$, this is at least $(1-\zeta_t)\mu_t^\star-\zeta_t M_t\ge(1-\zeta_t)\underline\mu_t-\zeta_t M_t$. Multiplying by the measurable indicator gives \eqref{e:eq:risk-lower}, with zero advantage on the reference response. The tower property and the two-cut identity~\eqref{e:eq:two-cut} give the policy contrast, including the certification costs already carried by both trajectories.
\end{proof}

\begin{corollary}[Simultaneous advantage and cumulative checkpoint risk]\label{e:cor:simultaneous-risk}
Under the common-prefix policy comparison of Theorem~\ref{e:thm:risk}, suppose a single event with failure probability at most $\zeta$ guarantees $A_t^\beta\ge\underline\mu_t$ at all executed decisions over the finite horizon. If $A_t^\beta\ge-M_t$ there and $[\underline\mu_t+M_t]_+\le L_t$ for deterministic finite $L_t\ge0$, then
\begin{equation}\label{e:eq:risk-uniform}
 \mathcal J(\alpha_{\rm gate})-\mathcal J(\beta_{\rm gate})
 =\sum_t\mathbb E[\xi_t A_t^\beta]
 \ge\sum_t\mathbb E[\xi_t\underline\mu_t]-\zeta\sum_t L_t .
\end{equation}
Independently, on a finite or countable horizon with countably many possible checkpoints $h$ at each $t$, let $S_{t,h}\in\mathcal G_t$ indicate that the action executes and checkpoint $h$ is selected before stopping. Suppose on $S_{t,h}$ a violation $\mathcal V_{t,h}$ entails failure of a certificate whose conditional failure probability given $\mathcal G_t$ is at most the measurable allocation $\zeta_{t,h}\ge0$. Under the reached policy law,
\begin{equation}\label{e:eq:failure-budget}
 \mathbb P\!\left(\bigcup_{t,h}(S_{t,h}\cap\mathcal V_{t,h})\right)
 \le \mathbb E\!\sum_{t,h}\mathbf1_{S_{t,h}}\zeta_{t,h}.
\end{equation}
\end{corollary}
\begin{proof}
The positive shortfall $[\xi_t(\underline\mu_t-A_t^\beta)]_+$ is zero on the joint event and at most $L_t$ on its complement. Take expectations and use~\eqref{e:eq:two-cut} for~\eqref{e:eq:risk-uniform}. For~\eqref{e:eq:failure-budget}, conditioning on $\mathcal G_t$ bounds each joint violation probability by $\mathbb E[\mathbf1_{S_{t,h}}\zeta_{t,h}]$. The union bound and Tonelli's theorem sum these terms.
\end{proof}

\subsection{Conditional adaptation frontiers}\label{e:sec:frontier}
The attainable frontier compares policy families at a common initial state and reserve.

\begin{definition}[Conditional adaptation frontier]\label{e:def:frontier}
Fix a probability space $(\Omega,\mathcal F,\mathbb P_{\rm env})$ carrying the exogenous trace and any policy random draws. Give each policy $\pi$ its decision filtration $(\mathcal G_t^\pi)$ generated by its available observations and draws. Fix a common initial state $q_0$, an independent evaluation law $\mu_{\rm eval}$ whose held-out realizations and outcomes are unavailable to policy selection, a total-resource budget $\mathbf b$, and a family $\mathfrak P$ of admissible policies, each adapted to its decision filtration. Let $\mathsf{Out}:\mathcal S\to\mathcal Y$ be a terminal output map and $\ell:\mathcal Y\times\mathcal Z\to\mathbb R$ an evaluation loss, with the displayed expectations well-defined. Define the attainable adaptation loss
\[
\mathsf L^*_{\mathfrak P}(q_0,\mathbf b;\mu_{\rm eval})
:=\inf_{\substack{\pi\in\mathfrak P\\m(\pi,\omega)\preceq\mathbf b\quad\mathbb P_{\rm env}\text{-a.s.}}}
\mathbb E_{\omega\sim\mathbb P_{\rm env}}\mathbb E_{z\sim\mu_{\rm eval}}
\big[\ell(\mathsf{Out}(q_T^\pi(\omega)),z)\big].
\]
Set $\inf\varnothing=+\infty$. The pathwise threshold includes the resources required for certification and execution.

\end{definition}

\begin{proposition}[Conditional frontier monotonicity and non-nested renewal]\label{e:prop:frontier}
For two families of policies whose final evaluation law, full parent state, information and resource semantics agree, if $\mathfrak P_0\subseteq\mathfrak P_1$ as affordable policy sets, then $\mathsf L^*_{\mathfrak P_1}\le\mathsf L^*_{\mathfrak P_0}$. Fixed-capacity replacement can make protected action subspaces nonnested and their Grams incomparable.
\end{proposition}
\begin{proof}
Taking an infimum over the enlarged affordable family gives the first inequality. For the second, in $\mathbb R^2$ with identity response, let one writable subspace be $\operatorname{span}\{e_1\}$ and its fixed-capacity replacement be $\operatorname{span}\{e_2\}$. Their response Grams are $\operatorname{diag}(1,0)$ and $\operatorname{diag}(0,1)$, neither of which is larger in Loewner order. Moreover a paid installation consumes reserve and can break equal-budget inclusion.
\end{proof}
\section{Causal Transformer state and physical decoders}\label{t:sec:source}
The controller writes to two internal projections of a causal Transformer. Its physical decoder acts on the weights that produce the current logits, while its optimizer registers help determine the next write. Table~\ref{t:tab:six} pairs the general objects with their realization in this system. Appendix~\ref{app:specializations} treats a final-head repair bank.

\begin{table}[ht]\centering\small
\caption{Objects of the controlled dynamics and their causal Transformer realization.}\label{t:tab:six}
\begin{tabularx}{\linewidth}{@{}L{0.16\linewidth}L{0.37\linewidth}Y@{}}\toprule
Component & Mathematical object & Causal Transformer realization\\\midrule
State and information & Execution state $(q,b)$, signature $\varsigma(q)$, information at the two cuts & Weights, optimizer memory, arrived tokens, reference logits and remaining reserve\\
Physical write & Decoder $D$, metric $M$, pullback $D^*MD$, admissible menu & Historical core and acquisition factor on each selected projection\\
Protection & Contract $\mathfrak C_g$, readouts and tolerances & Reference top token, per-position KL ceiling and frozen coefficients\\
Local response & Protection and response derivatives, cost metric, curvature bounds & Causal-prefix Jacobians, budgeted candidates and finite stored writes\\
Continuation & Branchwise transition maps or kernels and common terminal readout & AdamW moments and clocks, matched arrivals, subsequent losses and protected logits\\
Resources & Ordered reserve records $(m,c)$ over the executed stages & Sensing, native scans, write cost, resident and temporary storage\\\bottomrule
\end{tabularx}
\end{table}

\subsection{Causal model and complete implementation state}
Let $\mathfrak T$ denote a family of decoder-only causal Transformers. An instance specifies a finite sequence length $T$, block count $L_{\rm blk}$, hidden width $d$, vocabulary size $V$, projection modules $\mathcal M_{\rm mod}$ and its choices of attention, positional encoding, normalization, feed-forward activation and embedding tying. At a prefix $x_{\le j}$ its output logits are
\begin{equation}\label{t:eq:logits}
 h_{\vartheta}(x_{\le j};q)\in\R^d,\qquad
 z_j(q)=W_{\rm eff}(q)h_{\vartheta}(x_{\le j};q)+b_{\rm out}\in\R^V.
\end{equation}
The two co-located sources act inside $h_{\vartheta}$, with $W_{\rm eff}=W_{\rm out}$ at the head. The execution state is $x_t=(q_t,b_t)$, where $b_t$ is the available reserve and the physical and controller component has the form
\begin{equation}\label{t:eq:state}
 \begin{aligned}
 q_t=\big(&W_0,W_{\rm out},\{(A_\mu,B_\mu,s_\mu)\}_{\mu\in\mathcal M_{\rm LoRA}},
 \{(U_\mu,K_\mu,V_\mu)\}_{\mu\in\mathcal M_{\rm rep}},\\
 &\sigma_t,\mathfrak C_{g_t},\Xi^{\rm opt}_t,
 \mathsf{cursor}_t,\mathsf{store}_t,\mathsf{rng}_t,\mathsf{policy}_t,
 \mathsf{clock}_t\big).
 \end{aligned}
\end{equation}
Here $\sigma_t=\varsigma(q_t)$ identifies the installed cores, frames and writable slots. The optimizer register $\Xi^{\rm opt}_t$ contains the moments and per-array counters. The contract $\mathfrak C_{g_t}$ identifies the protected anchors and fixes their reference targets, tolerances and frozen registers. The data cursor and store determine which inputs have arrived and which the controller can reuse. The policy register retains observations used in later decisions. Each branch updates its reserve $b_t$ by the net consumption of its executed stages.
With the causal mask and randomness fixed, the logit at position $j$ depends only on embeddings in its prefix. Its Jacobian is block lower triangular:
\begin{equation}\label{t:eq:causality}
 D_{e_\tau}z_j=0\qquad(\tau>j).
\end{equation}
Across examples, the arrival filtration $\cI_t$ determines the information available for a probe. The ensuing action is measurable at the post-certification cut $\mathcal G_t$.

\subsection{Finite LoRA writes and their tangent decoder}
For any writable projection $W\in\R^{d_o\times d_i}$, use the conventional LoRA parameterization \citep{lora}
\begin{equation}\label{t:eq:lora}
 W_{\rm eff}=W_0+sBA,\qquad A\in\R^{r\times d_i},\quad B\in\R^{d_o\times r},\quad s=\alpha/r.
\end{equation}
Expanding both factors at $(A,B)$ separates the tangent decoder from the finite write:
\begin{equation}\label{t:eq:lora-finite}
 D_{\rm L}(\Delta A,\Delta B)=s(B\Delta A+\Delta B A),\qquad
 \Delta W_{\rm eff}=D_{\rm L}(\Delta A,\Delta B)+s\Delta B\Delta A.
\end{equation}
The bilinear term has Frobenius norm at most $|s|(\|\Delta A\|_\F^2+\|\Delta B\|_\F^2)/2$. In column coordinates $\theta=(\vecop A,\vecop B)$ the Jacobian and tangent metric are
\begin{equation}\label{t:eq:jac-lora}
 J_{\rm L}=s[\,I_{d_i}\otimes B\quad A^\top\otimes I_{d_o}\,],
 \qquad M_{\rm pull}=J_{\rm L}^*M_WJ_{\rm L}.
\end{equation}
For $M_W\succ0$ this metric is positive definite on the quotient by $\ker J_{\rm L}$. A finite write costs $\|\Delta W_{\rm eff}\|_{M_W}^2$. Factor values and optimizer registers remain part of the continuation state. At $B=0$, for example, changing $A$ leaves the current weight fixed and changes the coefficient gradient in Proposition~\ref{t:prop:neutral}.

\subsection{Co-located historical and acquisition sources}\label{t:subsec:two-interfaces}
\begin{figure}[tbp]
\centering
\includegraphics[trim=4pt 4pt 4pt 4pt,clip]{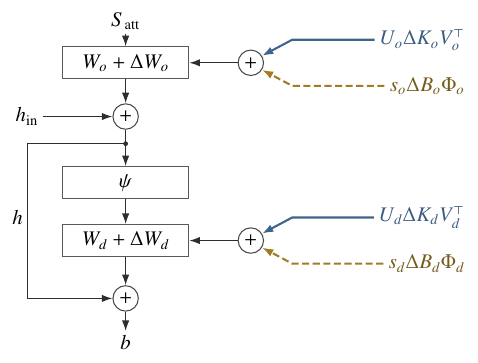}
\caption{Co-located writes in the attention-output and MLP-output projections of one block. In each projection the historical core (solid blue) and acquisition factor (dashed ochre) add to the physical write in~\eqref{t:eq:colocated-decoder}. For fixed attention features $S_{\rm att}$, the first sum gives $h=h_{\rm in}+(W_o+\Delta W_o)S_{\rm att}$. The second gives $b=h+(W_d+\Delta W_d)\psi(h)$. Thus the downstream module acts on features already changed by the upstream write, producing the mixed response in~\eqref{t:eq:block-mixed-finite}.}\label{fig:two-interfaces}
\end{figure}

For any selected module $\mu\in\mathcal M_{\rm write}$, let $U_\mu\in\R^{d_o\times r_o}$ and $V_\mu\in\R^{d_i\times r_o}$ have orthonormal columns. The writable coordinates are the old repair core $K_\mu\in\R^{r_o\times r_o}$ and the new LoRA output factor $B_\mu\in\R^{d_o\times r_n}$. At the parent, fix its input frame $\Phi_\mu\in\R^{r_n\times d_i}$. The exact \emph{linear effective-weight} source map on the fixed-basis chart is
\begin{equation}\label{t:eq:colocated-decoder}
 \Delta W_\mu=U_\mu\Delta K_\mu V_\mu^\top+s_\mu\Delta B_\mu\Phi_\mu,
 \qquad D_X=V_\mu\otimes U_\mu,\quad
 D_Y=s_\mu(\Phi_\mu^\top\otimes I_{d_o}).
\end{equation}
The effective weight is linear in $(\Delta K_\mu,\Delta B_\mu)$. The full-network guard evaluates the downstream logit response at each finite candidate. Basis formation, installation, renewal and changes to $\Phi_\mu$ carry their own action and resource records.
Figure~\ref{fig:two-interfaces} shows the source at the two writable projections of one block.

\subsection{Exact cross-cost whitening and source-restricted joint compensation}
\label{t:subsec:joint-physical}
Let a linear source slice have old/new coordinates $(x,y)$ and decoder
$u=D_Xx+D_Yy$ in a common real physical Hilbert space. Assume a frozen coercive physical metric $M$ on $u$, and for the following finite-dimensional chart suppose the pullback
\begin{equation}\label{t:eq:cross-pullback}
 \Gamma=D^*MD=\begin{bmatrix}M_{XX}&M_{XY}\\M_{XY}^{\top}&M_{YY}\end{bmatrix}\succ0,\qquad D=[D_X\ D_Y],
\end{equation}
so redundant codes have already been resolved in a chart preserving the old and new source blocks. For a co-located module with a separable metric $M_W=M_{\mathrm{in}}\otimes M_{\mathrm{out}}$ and decoder (\ref{t:eq:colocated-decoder}), its exact blocks are
\begin{align}\label{t:eq:cross-blocks}
 M_{XX}&=(V_{\mu}^{\top}M_{\mathrm{in}}V_{\mu})\otimes(U_{\mu}^{\top}M_{\mathrm{out}}U_{\mu}),\\
 M_{YY}&=s_{\mu}^2(\Phi_{\mu}M_{\mathrm{in}}\Phi_{\mu}^{\top})\otimes M_{\mathrm{out}},\nonumber\\
 M_{XY}&=s_{\mu}(V_{\mu}^{\top}M_{\mathrm{in}}\Phi_{\mu}^{\top})\otimes(U_{\mu}^{\top}M_{\mathrm{out}}).\nonumber
\end{align}
If $K_{\mu}$ is constrained diagonal, replace $D_X$ by $(V_{\mu}\otimes U_{\mu})S_{\rm diag}$ and accordingly restrict the $XX$ and $XY$ blocks. The Hadamard formula $M_{XX}= (V_{\mu}^{\top}M_{\mathrm{in}}V_{\mu})\odot(U_{\mu}^{\top}M_{\mathrm{out}}U_{\mu})$ follows.

\begin{proposition}[Injectivity of the full co-located core/output-factor chart]\label{t:prop:source-injectivity}
Suppose $s_\mu\ne0$, $U_\mu,V_\mu$ have full column rank, $K_\mu$ ranges over all $r_o\times r_o$ matrices, and $\Phi_\mu$ has full row rank while $\Delta B_\mu$ ranges over all $d_o\times r_n$ matrices. Then the co-located decoder $(\Delta K_\mu,\Delta B_\mu)\mapsto\Delta W_\mu$ in (\ref{t:eq:colocated-decoder}) is injective \emph{if and only if}
\begin{equation}\label{t:eq:right-overlap}
 \Ran(V_\mu)\cap\Ran(\Phi_\mu^\top)=\{0\}.
\end{equation}
For every positive definite physical metric $M$, this is equivalent to positive definiteness of $\Gamma=D^*MD$ for this full source chart. For $M=M_{\rm in}\otimes M_{\rm out}$, the condition $V_\mu^\top M_{\rm in}\Phi_\mu^\top=0$ gives $M_{XY}=0$. Restricting $K_\mu$ to diagonal matrices preserves the sufficiency of~\eqref{t:eq:right-overlap} for injectivity.
\end{proposition}
\begin{proof}
Each component map is separately injective under the rank assumptions. If their images meet in a nonzero matrix, its row space lies in both $\Ran(V_\mu)$ and $\Ran(\Phi_\mu^\top)$, hence the intersection is nontrivial. Conversely, for nonzero $w=V_\mu a=\Phi_\mu^\top b$ in that intersection and nonzero $u\in\Ran(U_\mu)$, the rank-one matrix $uw^\top$ lies in the image of each component, yielding a nonzero cancellation pair. A positive-definite metric has $\langle Dz,MDz\rangle=0$ exactly when $Dz=0$. Restricting the core to its diagonal only reduces its image, so disjoint right spans still ensure injectivity. Orthogonality of right spans under $M_{\rm in}$ makes the cross-block in~\eqref{t:eq:cross-blocks} vanish.
\end{proof}

\begin{proposition}[Cross-metric whitening and source-restricted compensated response]\label{t:prop:cross-whiten}
Assume (\ref{t:eq:cross-pullback}) is positive definite and set
\begin{equation}\label{t:eq:schur-whiten}
 S_c=M_{YY}-M_{XY}^{\top}M_{XX}^{-1}M_{XY}\succ0,\quad
 \bar x=M_{XX}^{1/2}(x+M_{XX}^{-1}M_{XY}y),\quad \bar y=S_c^{1/2}y.
\end{equation}
Then the exact physical squared cost is $\langle u,Mu\rangle=\|\bar x\|^2+\|\bar y\|^2$. In these coordinates
\begin{equation}\label{t:eq:two-white-decoders}
 u=\bar D_X\bar x+\bar D_Y\bar y,\quad
 \bar D_X=D_XM_{XX}^{-1/2},\quad
 \bar D_Y=(D_Y-D_XM_{XX}^{-1}M_{XY})S_c^{-1/2}.
\end{equation}
Let $J_P$ be the differentiable parent-to-protected-logit Jacobian and $J_R$ the parent-to-target Jacobian, both on a fixed causal frame. Define
$C=J_P\bar D_X$, $B=J_P\bar D_Y$, $H=J_R\bar D_X$ and $\Psi=J_R\bar D_Y$. For the full linear source chart, the \emph{linearized} protected-response problem
\begin{equation}\label{t:eq:source-joint-block}
 C\bar x+B\bar y=0,\qquad H\bar x+\Psi\bar y=d,
\end{equation}
has the following compensated Gram in physically whitened source coordinates (Theorem~\ref{m:thm:weighted})
\begin{equation}\label{t:eq:source-joint-gram}
 \begin{gathered}
 P_C=I-C^\dagger C,\quad N=HP_C,\quad \Ran Q=\Ker((I-CC^\dagger)B),\quad Q^\top Q=I,\\
 L=(\Psi-HC^\dagger B)Q,\quad W=I+Q^\top B^\top(C^\dagger)^\top C^\dagger BQ,\\
 G_{\mathrm{src},+}=NN^\top+LW^{-1}L^\top.
 \end{gathered}
\end{equation}
A target $d$ is algebraically reachable in this full linear source chart exactly when $d\in\Ran G_{\mathrm{src},+}$, and its minimum squared physical \emph{write} cost equals $d^\top G_{\mathrm{src},+}^\dagger d$. If the writable set contains the entire whitened ball of squared radius $b\ge0$, reachability within that budget is equivalent to both $d\in\Ran G_{\mathrm{src},+}$ and $d^\top G_{\mathrm{src},+}^\dagger d\le b$.
\end{proposition}
\begin{proof}
Complete the square in $(x,y)^\top\Gamma(x,y)$ to obtain (\ref{t:eq:schur-whiten}) and invert the block triangular change of variables, obtaining (\ref{t:eq:two-white-decoders}). Every solution to the first constraint of (\ref{t:eq:source-joint-block}) has $\bar y=Qv$ and $\bar x=-C^\dagger BQv+z$, $z\in\Ker C$; hence the second constraint is $Hz+Lv=d$ and the cost is $\|z\|^2+v^\top Wv$. Whitening $v$ gives the constrained minimum-norm response operator $[HP_C\ \ LW^{-1/2}]$, whose Gram is (\ref{t:eq:source-joint-gram}). The Moore--Penrose minimum-norm formula gives the range and cost statements. Comparison of this minimum with $b$ gives budget feasibility.
\end{proof}
\begin{remark}[Singular charts and executable menus]\label{t:rem:authorization}
For singular $\Gamma$, the joint quotient by $\Ker\Gamma=\Ker D$ represents physical writes. The readouts $J_P D$ and $J_R D$ descend because they vanish on $\Ker D$. To use the two-block whitening~\eqref{t:eq:schur-whiten}, first quotient the old and new source spaces by $\Ker D_X$ and $\Ker D_Y$, respectively. The resulting product has a positive-definite pullback exactly when $\Ran D_X\cap\Ran D_Y=\{0\}$. If the ranges overlap, retain the old block and restrict the reduced new block to the preimage of a complement of this intersection in $\Ran D_Y$. This product chart preserves the full joint physical image and assigns shared directions to the old block. Source-specific action domains and finite menus refer to their original writes (Theorem~\ref{e:thm:menu} and inequality~\eqref{e:eq:menu-observable-sharp}).

\end{remark}
\begin{example}[Cross-cost versus source collision]\label{t:ex:cross-cost}
Take $M=I_2$, $D_X=e_1$, $D_Y=e_1+e_2$. Then $M_{XX}=1$, $M_{XY}=1$, $M_{YY}=2$ and $S_c=1$. At fixed $y$, the cost is $(x+y)^2+y^2$ and its minimizer is $x=-y$. If $D_Y=e_1$, then $S_c=0$ and the joint quotient identifies $(x,y)$ with $x+y$. Keeping the old block leaves no independent new direction, although either original source can write $e_1$.
\end{example}

\subsection{Two modules in one block: separable cost and joint finite response}\label{t:subsec:two-modules}
Suppose the selected architecture exposes an attention output projection \texttt{o\_proj} and an MLP output projection \texttt{down\_proj} in the same block. Denote them by $o$ and $d$. They are distinct physical weight arrays; each may contain the co-located core and LoRA sources of (\ref{t:eq:colocated-decoder}). Their interaction through the causal computation graph is captured by the mixed derivatives below.

\begin{proposition}[Cross-module cost and response]\label{t:prop:two-modules}
Let $\Delta w_\mu=\vecop\Delta W_\mu=D_\mu v_\mu$, $\mu\in\{o,d\}$, be fixed linear physical source charts on disjoint, unaliased weight arrays, and put $\Delta w=(\Delta w_o,\Delta w_d)$. If the common physical metric is separable, $M_\oplus=\diag(M_o,M_d)$, then
\begin{equation}\label{t:eq:module-pullback}
 \Gamma_\oplus=
 \begin{bmatrix}D_o^*M_oD_o&0\\0&D_d^*M_dD_d\end{bmatrix},\qquad
 \|\Delta w\|_{M_\oplus}^2
 =\|D_ov_o\|_{M_o}^2+\|D_dv_d\|_{M_d}^2.
\end{equation}
This removes the \emph{cross-module} physical cost term; the old-core/new-factor cross term inside either $D_\mu$ may remain nonzero.

Let $P(v_o,v_d)$ be the full-network protected-logit map on a $C^2$ branch containing the whole segment $tv$, $0\le t\le1$, where $v=(v_o,v_d)$. Write $J_\mu=D_{v_\mu}P(0)$ and $H_{\mu\nu}(tv)=D^2_{v_\mu v_\nu}P(tv)$, with fixed action norms. Taylor's integral formula gives
\begin{equation}\label{t:eq:module-cross-remainder}
 \begin{aligned}
 P(v)-P(0)-J_ov_o-J_dv_d
 =\int_0^1(1-t)\big(&H_{oo}(tv)[v_o,v_o]
 +2H_{od}(tv)[v_o,v_d]\\
 &+H_{dd}(tv)[v_d,v_d]\big)\,dt.
 \end{aligned}
\end{equation}
If $\|H_{\mu\nu}(tv)\|_{\op}\le M_{\mu\nu}$ throughout this segment, then
\begin{equation}\label{t:eq:module-cross-bound}
 \|P(v)-P(0)-J_ov_o-J_dv_d\|
 \le\tfrac12 M_{oo}\|v_o\|^2
 +M_{od}\|v_o\|\|v_d\|
 +\tfrac12 M_{dd}\|v_d\|^2.
\end{equation}
A separable physical cost is compatible with a nonzero mixed Hessian.
\end{proposition}
\begin{proof}
The block decoder $\diag(D_o,D_d)$ and separable metric give (\ref{t:eq:module-pullback}) by multiplication. The second derivative of $t\mapsto P(tv_o,tv_d)$ contains the two equal mixed terms displayed in (\ref{t:eq:module-cross-remainder}); integrate its second derivative and apply the three uniform bounds. The downstream graph may depend on the changed output of both modules, independently of how the physical cost is assigned.
\end{proof}
For a sequential attention/MLP block, this interaction is already visible before the downstream network. Hold the attention features $S_{\rm att}$ and all other weights fixed, let $h=h_0+\Delta W_oS_{\rm att}$, and write
$b(\Delta W_o,\Delta W_d)=h+(W_d+\Delta W_d)\psi(h)$,
where $\psi$ includes the fixed post-attention normalization and MLP operations preceding \texttt{down\_proj}. The exact mixed finite difference is
\begin{equation}\label{t:eq:block-mixed-finite}
 \begin{aligned}
 &b(\Delta W_o,\Delta W_d)-b(\Delta W_o,0)-b(0,\Delta W_d)+b(0,0)\\
 &\hspace{1em}=\Delta W_d\big[\psi(h_0+\Delta W_oS_{\rm att})-\psi(h_0)\big].
 \end{aligned}
\end{equation}
On a chart where $\psi$ is $L_\psi$-Lipschitz in Frobenius norm, the norm of this mixed difference is at most
$L_\psi\|\Delta W_d\|_\F\|\Delta W_oS_{\rm att}\|_\F$.
Further downstream nonlinearities contribute the additional curvature already included in the full-network map $P$.
A coercive source slice or a verified norm conversion relates the coordinate bounds to physical cost (Proposition~\ref{t:prop:remainder}). The joint native guard tests the mixed response of the stored candidate against the full-vocabulary and top-token obligations.

\section{Native protection, response, and precision}\label{t:sec:guards}
\subsection{Full-vocabulary KL and its global logit bound}
For $T_\mathrm{KL}>0$, let $p=\softmax(z/T_\mathrm{KL})$ and $p_\delta=\softmax((z+\delta)/T_\mathrm{KL})$. Then, for any fixed reference probability vector $q$,
\begin{equation}\label{t:eq:KL-exact}
 \KL(q\|p_\delta)-\KL(q\|p)
 =\log\sum_{v=1}^{V}p_v e^{\delta_v/T_\mathrm{KL}}-\frac{q^\top\delta}{T_\mathrm{KL}}.
\end{equation}
When the parent is the reference ($q=p$), the first derivative vanishes at $\delta=0$. Since the log-sum-exp Hessian is $T_\mathrm{KL}^{-2}(\diag(p)-pp^\top)$ and has operator norm at most $1/(2T_\mathrm{KL}^2)$, the descent inequality yields
\begin{equation}\label{t:eq:kl-bound}
 0\le\KL\left(p\middle\|p_\delta\right)\le\frac{\|\delta\|_2^2}{4T_\mathrm{KL}^2}.
\end{equation}
Constant shifts preserve softmax probabilities, and rare-token directions can carry arbitrarily little KL even after centering. For a persistent reference $q\ne p$, the linear term of~\eqref{t:eq:KL-exact} is $(p-q)^\top\delta/T_\mathrm{KL}$. The identity and bound place the reference distribution in the first argument at the fixed temperature $T_{\rm KL}$.

\subsection{Fixed low-dimensional anchor protection versus full-vocabulary KL}
A fixed semantic logit projection reduces the number of protection-Jacobian rows. At each parent anchor token $n=1,\dots,N_{\rm anc}$, freeze the parent logit $z_{0,n}$ and a matrix $U_n\in\R^{V\times k_n}$ with orthonormal columns in $\mathbf1^\perp$. Put $\Pi_\circ=I-\mathbf1\mathbf1^\top/V$ and $\delta_n=z_n(u)-z_{0,n}$. Define
\begin{equation}\label{t:eq:proj-anchor}
 p_{\rm proj}(u)=N_{\rm anc}^{-1/2}\big(U_n^\top\Pi_\circ\delta_n\big)_{n=1}^{N_{\rm anc}},\quad
 p_{\rm full}(u)=N_{\rm anc}^{-1/2}\big(\Pi_\circ\delta_n\big)_{n=1}^{N_{\rm anc}}.
\end{equation}
Both vanish at the parent and are differentiable whenever the chart is. A projected guard controls exactly its $\sum_n k_n$ coordinates.
\begin{proposition}[Full-KL certificate from projected logits and an explicit remainder]\label{t:prop:projected-kl}
At a fixed temperature $T_{\rm KL}>0$, let $p_{0,n}=\softmax(z_{0,n}/T_{\rm KL})$. Define the omitted centered-logit energy
\begin{equation}\label{t:eq:residual-energy}
 R_{\rm omit}(u)^2=\frac1{N_{\rm anc}}\sum_{n=1}^{N_{\rm anc}}\|(I-U_nU_n^\top)\Pi_\circ\delta_n\|_2^2.
\end{equation}
Then for all finite logits with the stated reference temperature,
\begin{equation}\label{t:eq:proj-to-kl}
 \frac1{N_{\rm anc}}\sum_n\KL\big(p_{0,n}\,\|\,\softmax((z_{0,n}+\delta_n)/T_{\rm KL})\big)
 \le\frac{\|p_{\rm proj}(u)\|_2^2+R_{\rm omit}(u)^2}{4T_{\rm KL}^2}.
\end{equation}
In particular, $\|p_{\rm full}(u)\|_2\le\tau$ implies mean parent-forward KL at most $\tau^2/(4T_{\rm KL}^2)$ for arbitrary finite logits.
\end{proposition}
\begin{proof}
Apply~\eqref{t:eq:kl-bound} to the centered displacement at each position. Orthogonality splits its squared norm into the projected and omitted energies. Averaging gives~\eqref{t:eq:proj-to-kl}.
\end{proof}
\begin{corollary}[Maximum-position KL and top-token protection]\label{t:cor:max-kl-margin}
In the setting of Proposition~\ref{t:prop:projected-kl}, put
$K_j(u)=\KL(p_{0,j}\|\softmax((z_{0,j}+\delta_j)/T_{\rm KL}))$.
At every position $j$,
\begin{equation}\label{t:eq:max-position-kl}
 K_j(u)\le\frac{\|\Pi_\circ\delta_j\|_2^2}{4T_{\rm KL}^2}.
\end{equation}
Thus $\max_j\|\Pi_\circ\delta_j\|_2\le 2T_{\rm KL}\sqrt{\tau_{\rm KL}}$ suffices for a ceiling $\tau_{\rm KL}\ge0$ on every position.

Suppose $V\ge2$ and the parent has a unique top token $v_j^\star$ at position $j$, with margin
\begin{equation}\label{t:eq:parent-top-margin}
 m_j=\min_{v\ne v_j^\star}(z_{0,j,v_j^\star}-z_{0,j,v})>0.
\end{equation}
It remains the unique top token after the write exactly when
\begin{equation}\label{t:eq:top-pairwise}
 \delta_{j,v}-\delta_{j,v_j^\star}
 <z_{0,j,v_j^\star}-z_{0,j,v}\qquad(v\ne v_j^\star).
\end{equation}
In particular, $\|\delta_j\|_\infty<m_j/2$ suffices; it also suffices to replace $\delta_j$ by its centered version $\Pi_\circ\delta_j$. These are per-position conditions on the full vocabulary.
\end{corollary}
\begin{proof}
The KL bound holds before averaging in~\eqref{t:eq:proj-to-kl}. The new margin against $v$ changes by $\delta_{j,v_j^\star}-\delta_{j,v}$. Centering preserves this difference, whose magnitude is at most $2\|\delta_j\|_\infty$.
\end{proof}

\begin{figure}[tbp]
\centering
\includegraphics[trim=4pt 4pt 4pt 4pt,clip]{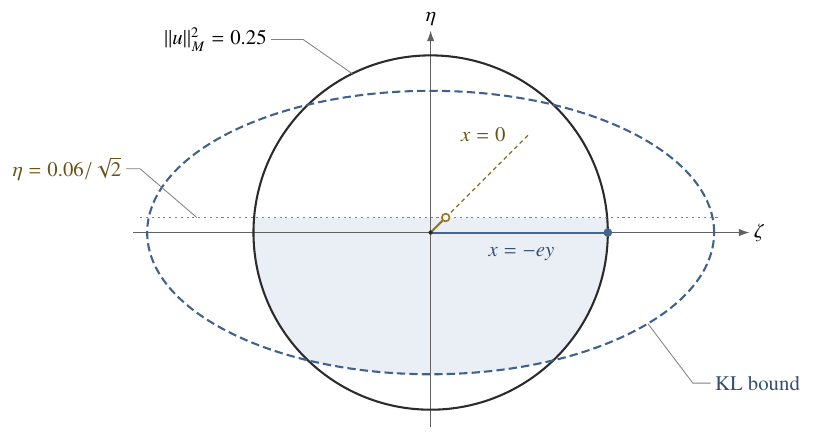}
\caption{Response under a narrow decision margin on the weighted minimum-cost slice. Take $M_X=\operatorname{diag}(2,8)$, $C=(1,2)$, $\xi=Cx$, $\zeta=(y-\xi)/\sqrt2$ and $\eta=(\xi+y)/\sqrt2$. For $z_0=(0.06,0,-4)$ at temperature one, the logit increment is $\delta=(-\eta/\sqrt2+\zeta/(2\sqrt6),\,\eta/\sqrt2+\zeta/(2\sqrt6),\,-\zeta/\sqrt6)$. Writing $e=M_X^{-1}C^\top(CM_X^{-1}C^\top)^{-1}=(1/2,1/4)^\top$ and $x=e\xi+v$ with $v\in\ker C$ gives squared cost $\zeta^2+\eta^2+v^\top M_Xv$. Thus $v=0$ gives the unique least-cost lift for each plotted point. Shading marks the intersection of the budget disk, the sufficient KL ellipse $\eta^2+\zeta^2/4\le0.16$ (parent-forward KL at most $0.04$ by Corollary~\ref{t:cor:max-kl-margin}), and the strict top-token condition $\eta<0.06/\sqrt2$. The third token stays below the first throughout the disk. With response $r=y$ and $y\ge0$, the uncompensated ray $x=0$ reaches a tie at $y=0.06$ (open circle). The compensated ray $x=-ey$ satisfies $Cx+y=0$, keeps the first-pair margin at $0.06$, and reaches $y=0.5/\sqrt2$ at the budget boundary (filled circle).}
\label{fig:margin-compensation}
\end{figure}

Figure~\ref{fig:margin-compensation} locates the decision boundary, a sufficient KL ellipse and the weighted minimum-cost repair on one control slice.

\begin{example}[Positive full KL at zero projected displacement]\label{t:ex:proj-not-full}
For a three-token vocabulary, let $z_0=(0,0,0)$, $U=(1,-1,0)^\top/\sqrt2$ and $\delta=(1,1,-2)$. Then $U^\top\Pi_\circ\delta=0$, while $\delta$ is not a constant vector, so $\KL(\softmax(z_0)\|\softmax(z_0+\delta))>0$.
\end{example}
\paragraph{First-order protection.}
Differentiating~\eqref{t:eq:proj-anchor} gives $C_i=N_{\rm anc}^{-1/2}(U_n^\top\Pi_\circ D_u z_{i,n}(0))_n$. Omitting $U_n^\top$ gives the full centered-logit Jacobian. At the reference parent the derivative of scalar forward KL vanishes by~\eqref{t:eq:KL-exact}, so its kernel is the entire source. The projected or centered logit rows constrain the search to first order. Corollary~\ref{t:cor:max-kl-margin} certifies a finite write from its full centered displacement at each position. The implementation evaluates the stored candidate on the full vocabulary.

To apply Theorem~\ref{e:thm:tolerance} to a KL-admitted old menu, its logit tolerance must cover every compared old endpoint. At a reference parent, suppose the softmax covariance along every compared logit segment is bounded below by $\gamma I$ on its centered displacement direction, for a common $\gamma>0$. Taylor's integral formula gives $K_j(u)\ge\gamma\|\Pi_\circ\delta_j\|^2/(2T_{\rm KL}^2)$. A per-position KL ceiling then places each centered displacement in the ball of radius $T_{\rm KL}\sqrt{2\tau_{\rm KL}/\gamma}$.

\subsection{Adaptive local charts with fixed native obligations}\label{t:subsec:adaptive-guards}
Within the paid search of Proposition~\ref{e:prop:adaptive-search}, a rejected trial can reveal new violating positions. The next trial may add their Jacobian rows or rebuild its Fisher or top-token search basis. Each candidate uses a fixed parent, source, metric and remainder domain. Admission requires a stored write within its physical budget and, at each protected position, a KL divergence within the persistent-reference ceiling and the same reference top token. Restoring the model and optimizer between trials preserves the common parent while the probe transcript and charges accumulate.

\paragraph{Persistent references at a new parent.}
At one protected position, let $z_t$ be the new parent's logits and $\delta_t=z_t(u)-z_t$ the candidate increment. Write $p_t=\softmax(z_t/T_{\rm KL})$ and $p_t(u)=\softmax(z_t(u)/T_{\rm KL})$. For any fixed reference distribution $p^{\rm ref}$, subtracting the two instances of~\eqref{t:eq:KL-exact} gives the finite identity
\[
 \KL(p^{\rm ref}\|p_t(u))
 =\KL(p^{\rm ref}\|p_t)
 +\frac{(p_t-p^{\rm ref})^\top\delta_t}{T_{\rm KL}}
 +\KL(p_t\|p_t(u)).
\]
Equation~\eqref{t:eq:kl-bound} bounds the final term. When $p^{\rm ref}=\softmax(z^{\rm ref}/T_{\rm KL})$ for reference logits $z^{\rm ref}$, Corollary~\ref{t:cor:max-kl-margin} applies directly at $z^{\rm ref}$ with the \emph{total} offset $(z_t-z^{\rm ref})+\delta_t$. Its pairwise inequalities then use the reference top token and its reference margins. With a distribution-only KL obligation, the finite identity above directly evaluates the absolute divergence.

For a fixed linear logit readout $\mathsf R_{\rm obs}$, write $e_t=\mathsf R_{\rm obs}(z_t-z^{\rm ref})$ and $J_{\rm obs,t}=D_u[\mathsf R_{\rm obs}z_t(u)]_{u=0}$. The homogeneous search condition $J_{\rm obs,t}u=0$ preserves the current residual to first order. A search targeting restoration of the reference readout solves $J_{\rm obs,t}u=-e_t$ on the writable source when feasible. For either search, admission uses the finite stored candidate and the contract tolerance.

\subsection{Shared weights couple the complete causal prefix}\label{t:subsec:cross-token}
For a writable projection $p$ in block $\ell$, the factor-to-logit Jacobian is
\begin{equation}\label{t:eq:network-J}
 J_{z,\ell p}=\mathcal J_{\mathrm{down},\ell p}\,\mathcal J_{\mathrm{insert},\ell p}\,J_{\mathrm L,\ell p}.
\end{equation}
Here $J_{\mathrm L,\ell p}$ decodes a coefficient increment into a physical weight increment. The insertion map applies that weight increment to the fixed input at every use of the projection. The downstream differential transports the resulting activations through the masked causal graph, including any grouped-query sharing.

Let $a_{\mu,\tau}$ be the output of the projection $W_\mu$ with input $h^{\mathrm{in}}_{\mu,\tau}$ at position $\tau$. Writing $G^{\mathrm{out}}_{\mu,n\leftarrow\tau}=D_{a_{\mu,\tau}}z_n$ at the parent gives
\begin{equation}\label{t:eq:causal-weight-jac}
 D_{W_\mu}z_n[\Delta W_\mu]
 =\sum_{\tau\le n}G^{\mathrm{out}}_{\mu,n\leftarrow \tau}(\Delta W_\mu h^{\mathrm{in}}_{\mu,\tau})
 =\sum_{\tau\le n}\big((h^{\mathrm{in}}_{\mu,\tau})^\top\otimes G^{\mathrm{out}}_{\mu,n\leftarrow \tau}\big)\vecop\Delta W_\mu.
\end{equation}
Terms with $\tau<n$ carry earlier perturbations into the target position. If the weight is also tied across depths, the sum includes each use at each depth. A full-prefix vector--Jacobian product gives the corresponding projected protection row:
\begin{equation}\label{t:eq:vjp-anchor}
 C_{0,j:}=N_{\rm anc}^{-1/2}D_x\Big[\big\langle u_j,\Pi_\circ(z_{n(j)}(D_Xx)-z_{0,n(j)})\big\rangle\Big]_{x=0},
\end{equation}
where $u_j$ is a column of $U_{n(j)}$ and the normalization is that of~\eqref{t:eq:proj-anchor}. Row blocks can be computed without forming the full vocabulary-by-parameter Jacobian. For a sketch-based approximation $\widehat C$ on the same source and readout spaces as $C_0$, Proposition~\ref{t:prop:numeric-null} bounds true leakage using a certified operator error $\|C_0-\widehat C\|_{\op}\le\delta_C$ and bounds on the computed projector norm and residual. By the singular-value perturbation bound, the thresholded ranks agree if every singular value of $\widehat C$, including zeros, lies more than $\delta_C$ from the cutoff. A truncated spectrum requires a certified upper bound on the omitted tail below the cutoff by the same margin.

\subsection{Loss curvature, LoRA remainder, and numerical projection}\label{t:sec:precision}
For one supervised token with label $y$, write $p=\softmax(z)$ and $\ell(z)=-\log p_y$. The logit gradient and Hessian are
\begin{equation}\label{t:eq:cekl-hess}
 \nabla_z\ell=p-e_y,\qquad \nabla_z^2\ell=\diag(p)-pp^\top.
\end{equation}
The Hessian has operator norm at most $1/2$. Consequently an affine logit response $z(u)=z(0)+Hu$ satisfies the global bound
\begin{equation}\label{t:eq:cekl-E0}
 0\le \ell(z(u))-\ell(z(0))-\langle\nabla_z\ell(z(0)),Hu\rangle
 \le\tfrac14\|Hu\|_2^2\le\tfrac14\|H\|_{\op}^2\|u\|^2.
\end{equation}
On a chart where the adapter response $z(u)$ is $C^2$, the composition Hessian is
$D^2(\ell\circ z)[u,u]=\langle J_z u,(\nabla_z^2\ell)J_zu\rangle+\langle\nabla_z\ell,D^2z[u,u]\rangle$. The second term introduces the curvature of the logit map.

\begin{proposition}[Conditional LoRA-to-protection finite remainder]\label{t:prop:remainder}
For $u=(\Delta A,\Delta B)$ let $\|u\|_{\rm fac}^2=\|\Delta A\|_\F^2+\|\Delta B\|_\F^2$. Fix the Frobenius norm on effective-weight perturbations and a norm on protected logits. At a frozen parent $(A,B)$ of scale $s\in\mathbb R$, put
\[
 \mathcal W(u)=s\big[(B+\Delta B)(A+\Delta A)-BA\big],\qquad p(u)=f(\mathcal W(u)).
\]
Assume $f$ is $C^2$ on a branch containing $\mathcal W(tu)$ for $0\le t\le1$, with $\|Df\|_{\op}\le L_f$ and $\|D^2f\|_{\op}\le M_f$ along this image in the stated norms. If $\|D\mathcal W(0)\|_{\op}=\|J_{\rm L}\|_{\op}\le J_*$, then
\[
 \begin{aligned}
 J_{\rm seg}&:=\sup_{0\le t\le1}\|D\mathcal W(tu)\|_{\op}
       \le J_*+|s|\,\|u\|_{\rm fac},\\
 \|p(u)-p(0)-Dp(0)u\|&\le
       \tfrac12\big(M_fJ_{\rm seg}^2+|s|L_f\big)\|u\|_{\rm fac}^2 .
 \end{aligned}
\]
If the same bounds hold for every such segment in a factor ball of radius $\rho$, a uniform remainder constant there is $M_f(J_*+|s|\rho)^2+|s|L_f$.
\end{proposition}
\begin{proof}
For $h=(h_A,h_B)$, the derivative difference applied to $h$ is $st(\Delta B h_A+h_B\Delta A)$. Frobenius submultiplicativity and Cauchy--Schwarz bound its norm by $t|s|\,\|u\|_{\rm fac}\|h\|_{\rm fac}$. Also $D^2\mathcal W[u,u]=2s\Delta B\Delta A$ has norm at most $|s|\,\|u\|_{\rm fac}^2$. The second derivative of $t\mapsto p(tu)$ is $D^2f[D\mathcal W(tu)u,D\mathcal W(tu)u]+Df[D^2\mathcal W[u,u]]$, with $f$ differentiated at $\mathcal W(tu)$. Integrating against $1-t$ gives the remainder.
\end{proof}
On a writable slice where $\|u\|_{\rm fac}\le\kappa\|u\|_{\rm cost}$ with $\kappa>0$, a factor-norm remainder constant $M$ on the ball of radius $\rho$ gives the cost-norm constant $\kappa^2 M$ on the ball of radius $\rho/\kappa$.

\begin{proposition}[Guard residual under a computed numerical null projector]\label{t:prop:numeric-null}
Let $C$ be the real-arithmetic protected Jacobian on a fixed source chart and $\widehat C$ its computed approximation, with $\|C-\widehat C\|_{\op}\le\delta_C$. If a computed projector $\widehat P$ satisfies $\|\widehat P\|_{\op}\le1+\delta_P$ and $\|\widehat C\widehat P\|_{\op}\le\epsilon_P$, then
\begin{equation}\label{t:eq:projector-leak}
 \|C\widehat P\|_{\op}\le \epsilon_P+\delta_C(1+\delta_P).
\end{equation}
For the exact orthogonal projector $\widehat P$ onto $\Ker\widehat C$ in the source inner product, $\|C\widehat P\|_{\op}\le\delta_C\|\widehat P\|_{\op}\le\delta_C$. If $\Ker\widehat C=\{0\}$, then $\widehat P=0$ and the residual vanishes.
\end{proposition}
\begin{proof}
Write $C\widehat P=\widehat C\widehat P+(C-\widehat C)\widehat P$ and apply the triangle inequality and operator submultiplicativity.
\end{proof}
\Needspace{8\baselineskip}
\section{Optimizer-native continuation and execution cost}\label{t:sec:execution}
\label{t:sec:consumer}
\subsection{Rank-coordinate symmetry at a neutral birth}
\begin{proposition}[Universal rank symmetries of the first write]\label{t:prop:neutral}
Attach an additive slot $Y=Y_0+sUA^\top B^\top$ with $B=0\in\R^{d_o\times r}$, $s\ne0$, and a full-row-rank carrier $A\in\R^{r\times d_i}$, where $d_o\ge1$ and $1\le r\le d_i$. Its output agrees with the parent on every input. For input rows $U\in\R^{n_{\rm in}\times d_i}$ ($n_{\rm in}\ge1$), output differential $R_Y=\nabla_Y\ell(Y_0)$ and a step of size $\eta>0$,
\begin{equation}\label{t:eq:birth-grad}
 \nabla_A\ell=0,\qquad G_B:=\nabla_B\ell=sR_Y^\top UA^\top .
\end{equation}
For $T\in GL(r,\mathbb R)$, compare $A$ and $TA$ at the same $U,Y_0,R_Y$, clocks and hyperparameters. A universal symmetry preserves the first physical write for every $U,R_Y$ under the unclipped update. Let $A^+=c_A A$, with $c_A=1$ for a frozen carrier or $c_A=1-\eta_A\lambda_A$ for a carrier stepped under decoupled decay with zero gradient and zero moments. Assume $c_A\ne0$.
\begin{enumerate}[label=(\roman*),leftmargin=*]
\item An SGD step on $B$ has universal group $O(r)$. With $B^+=-\eta G_B$, its weight and output increments are
\begin{equation}\label{t:eq:birth-sgd}
 \Delta W=-\eta s c_A G_B A,\qquad \Delta Y=U\Delta W^\top=-\eta s^2 c_A UA^\top AU^\top R_Y .
\end{equation}
\item For elementwise AdamW take zero $B$ moments, completed array clock $n$, $0\le\beta_1,\beta_2<1$ and $\epsilon>0$. Equation~\eqref{t:eq:adam} gives
\begin{equation}\label{t:eq:birth-adam}
 B^+=-\eta f_n(G_B),\qquad f_n(x)=\frac{a_n x}{b_n|x|+\epsilon},\qquad
 a_n=\frac{1-\beta_1}{1-\beta_1^{n+1}},\quad
 b_n=\sqrt{\frac{1-\beta_2}{1-\beta_2^{n+1}}}.
\end{equation}
Here $f_n$ acts entrywise. The universal group is the signed permutation group $B_r=\{\pm1\}^r\rtimes S_r$. For either update, the degenerate choice $c_A=0$ gives $\Delta W=0$ for all $T$.
\end{enumerate}
\end{proposition}
\begin{proof}
Differentiation at $B=0$ gives~\eqref{t:eq:birth-grad}. Under $A\mapsto TA$ the gradient becomes $G_B T^\top$. Every test row $g$ is realizable with one input row satisfying $uA^\top=g$ and one nonzero entry $1/s$ in $R_Y$. For SGD, equality of weight increments for all $g$ gives $T^\top T=I$ by right cancellation of $A$. For AdamW the decay of $B=0$ vanishes and bias correction gives~\eqref{t:eq:birth-adam}. Since $\Delta W=s c_A B^+ A$, universal equality for $c_A\ne0$ requires $f_n(gT^\top)T=f_n(g)$ for all $g$. Linearizing at $g=0$ gives $T^\top T=I$ because $f_n'(0)=a_n/\epsilon>0$. Right multiplication by $T^\top$ gives $f_n(gT^\top)=f_n(g)T^\top$. Taking $g=e_k^\top$ yields $f_n(T_{jk})=f_n(1)T_{jk}$ in every entry of column $k$. Thus each nonzero entry has magnitude one. Orthogonality gives one per row and column, and every such matrix commutes with the odd entrywise update.
\end{proof}
At fixed $U,R_Y$, two carriers $A_1,A_2$ give the same SGD output exactly when $U(A_1^\top A_1-A_2^\top A_2)U^\top R_Y=0$. For AdamW put $E_T=f_n(G_B T^\top)T-f_n(G_B)$. When $c_A\ne0$, equality of physical writes is equivalent to $E_T=0$, and equality of batch outputs to $UA^\top E_T^\top=0$. All test rows are available at fixed $U$ if $UA^\top$ has rank $r$. A particular $G_B$ may have additional symmetries, including rotations among its zero columns. For example, with $n=0$, $c_A=s=R_Y=1$ and $U=(1,0)$, the identity carrier $I_2$ and its $\pi/4$ rotation give the scalar output increments $-\eta/(1+\epsilon)$ and $-\eta/(1/\sqrt2+\epsilon)$.

\subsection{Exact finite AdamW state recursion and structural renewal}\label{t:subsec:adam-native}
For the trainable array $\theta$, let $n_t$ count its completed AdamW steps just before decision $t$, and let $g_t$ be its clipped gradient when the array is stepped. That update uses bias-correction index $n_t+1$ and sets $n_{t+1}=n_t+1$. A skipped array retains its counter under this convention. The elementwise recursion \citep{adamw} in real arithmetic is
\begin{align}
 m_{t+1}&=\beta_1m_t+(1-\beta_1)g_t, &v_{t+1}&=\beta_2v_t+(1-\beta_2)g_t^{\odot2},\label{t:eq:adam}\\
 \widehat m_{t+1}&=m_{t+1}/(1-\beta_1^{n_t+1}),&\widehat v_{t+1}&=v_{t+1}/(1-\beta_2^{n_t+1}),\nonumber\\
 \theta_{t+1}&=(1-\eta_t\lambda)\theta_t-\eta_t\,\widehat m_{t+1}\oslash(\sqrt{\widehat v_{t+1}}+\epsilon).\nonumber
\end{align}
The native transition also includes storage rounding (Section~\ref{t:sec:precision}). For branches with the same constants and clock, writing $g_t^a=g_t^0+\delta g_t$ in~\eqref{t:eq:adam} gives
\[
 \delta m_{t+1}=\beta_1\delta m_t+(1-\beta_1)\delta g_t,\qquad
 \delta v_{t+1}=\beta_2\delta v_t+(1-\beta_2)(2g_t^0+\delta g_t)\odot\delta g_t .
\]
The coordinatewise square in the second-moment register also appears in the invariant-subspace criterion of Theorem~\ref{m:thm:partition}. At birth $a_0=b_0=1$ for a fresh array. For a renewed zero-moment slice in an existing array, the same entrywise formula uses its retained counter and the gradient delivered after clipping.

\paragraph{Fixed-capacity slot replacement.}
A rank-$r$ projection is $s\sum_{j=1}^r b_ja_j^\top$. Replacing physical slot $j$ with $(\widetilde a_j,0)$ changes the effective weight by $-sb_ja_j^\top$ at staging, even though the \emph{newborn} component contributes zero. A slice reset changes the moments owned by the renewed slot and preserves the other slots and the array clock. The protocol checks the resulting finite write before committing the replacement \citep[Section~6.2]{paper4}.

\subsection{Atomic promotion to a second generation}\label{t:subsec:generation-promotion}
Promotion is a paid transition $(q^-,b^-)\leadsto(q^+,b^--c_{\rm prom})$ (Proposition~\ref{m:prop:contract-promotion}). It retains the existing references, tolerances and frozen values, and fixes new reference logits at positions from previously consumed data. At handoff the controller freezes the first-generation coefficients and optimizer state, including their ownership and write masks. A new carrier with $B_{\rm gen2}=0$ adds zero effective weight. Its probe checks the native logits at every protected position and verifies a nonzero coefficient gradient via~\eqref{t:eq:birth-grad}. The controller commits the new generation and contract together after restoring the probe state. A rejected transaction restores the previous checkpoint and retains its charge. The promoted carrier spans new input directions, changing the row space available to its first physical write. Within that space, Proposition~\ref{t:prop:neutral} identifies the signed permutations as the rank-coordinate changes preserving the first unclipped AdamW write for every input and loss differential at zero moments.

\begin{table}[htbp]\centering\small
\caption{Protection across two generations. The 32 Phase A contrasts define the fixed-parent diagnostic. In Phase B both online policies include all 40 fixed contrasts alongside adaptively selected rows. The native guard checks every listed position.}\label{t:tab:promotion-contract}
\begin{tabularx}{\linewidth}{@{}L{0.14\linewidth}L{0.16\linewidth}L{0.23\linewidth}Y@{}}\toprule
Phase & Fixed token contrasts & Native protected positions & Write-state contract\\\midrule
A & $32$ & Original $2{,}048$ & First-generation source and optimizer\\
B, after promotion & $40$ & Original $2{,}048$ plus $512$ consumed-data positions: $2{,}560$ total & First generation frozen, second generation installed with $B_{\rm gen2}=0$ and a verified writable learning signal\\\bottomrule
\end{tabularx}
\end{table}
The online solver rebuilds position-specific rows at the current parent. The guard of Section~\ref{t:subsec:adaptive-guards} applies to the full expanded anchor set. At fixed parent, source and metric, adding rows restricts the linear kernel. Installing the new frame also changes the source on which that kernel is computed.

\begin{figure}[tbp]
\centering
\includegraphics[trim=4pt 4pt 4pt 4pt,clip]{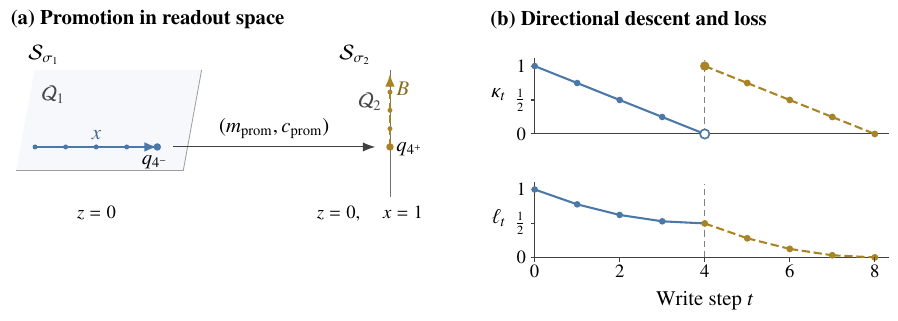}
\caption{A neutral promotion changes the writable direction at fixed output.
(a) In the analytic model $f(h)=xh_1+Bh_2+zh_3$, the common coordinates are the readouts $(f(e_1),f(e_2),f(e_3))=(x,B,z)$. The first contract protects $z=0$. Four writes of $1/4$ take $x$ from zero to one. The promotion retains $z=0$, fixes the newly protected value $x=1$, and installs $B=0$ with fixed input row $e_2^\top$. Thus $\mathcal Q_2\subseteq\mathcal Q_1$ in the common readout space. The two labeled states lie in different strata and represent the same function. The arrow carries the promotion's required reserve and net charge (Table~\ref{t:tab:promotion-contract}).
(b) With Euclidean write cost, the response Gram on $(x,B)$ changes from $e_1e_1^\top$ to $e_2e_2^\top$. For $\ell=\tfrac12[(1-x)^2+(1-B)^2]$, the directional descent coefficient $\kappa_t=\sqrt{\langle\nabla\ell,G_t\nabla\ell\rangle}$ jumps from zero just before commit (open circle) to one at commit (filled circle). The loss is $1/2$ on both sides and $\partial_B\ell=-1$ at birth. Four subsequent writes of $1/4$ to $B$ reduce the loss to zero. Lines join successive writes within each generation.}
\label{fig:versioned-protection}
\end{figure}

Figure~\ref{fig:versioned-protection} pairs this promotion contract with an exact example of protected descent across a function-neutral installation.

\subsection{Physical resource coordinates}\label{t:subsec:stages}
The net squared write cost uses the physical metric in~\eqref{t:eq:cross-pullback}. Operation counts accumulate across decisions, while a transient workspace contributes its peak requirement and releases its storage at the end of the stage. With a frozen acquisition frame of rank $r_n$, module $\mu$ has $d_o r_n$ writable acquisition coefficients and $r_o^2$ historical coefficients. Persistent reference logits require $n_C V$ scalar entries when stored in full. Table~\ref{t:tab:stage-ledger} assigns each resource to the stage that incurs it. The stage records compose by Proposition~\ref{m:prop:record}.

\begin{table}[H]\centering\footnotesize\renewcommand{\arraystretch}{1.12}
\caption{Resource accounting for a causal Transformer decision. Each stage lists its incurred work and the reserve or measurement needed for execution.}\label{t:tab:stage-ledger}
\begin{tabularx}{\linewidth}{@{}L{0.18\linewidth}Y L{0.29\linewidth}@{}}\toprule
Stage & Incurred work & Resource coordinates \\\midrule
Parent checkpoint & Save writable coefficients, optimizer moments and clocks, RNG and data cursor & $F_{\rm checkpoint}$, snapshot storage, copy bandwidth\\
Paid sensor and menu & Probes, gradient/VJP rows, rank tests, candidates and conditional simulation & $F_{\rm sense}$, scratch VRAM, retained caches\\
Native candidate checks & Stage each stored write, evaluate the protected logits and restore the parent between candidates & $F_{\rm check}$ for all trials, policy scans and separate audit scans\\
Commit or abstain & Slot installation and moment-slice transition, or reference action & $F_{\rm branch}$, installation bytes and rollback state\\
Ordinary continuation & Subsequent optimizer updates along the branch & $F_{\rm continue}$ and peak activation memory\\
Held-out evaluation & Terminal and path readouts independent of online selection & $F_{\rm eval}$ and peak device memory\\\bottomrule
\end{tabularx}
\end{table}
\section{Closed-loop evolution in a causal Transformer}\label{exp:protected-adaptation}
In Phase A, four policies processed the same 200 arrivals on Qwen2.5-7B Base in FP32 \citep{qwenbase}. The model offered a historical channel and an acquisition channel under a fixed protection reference. We first replay individual writes to identify the effect of each channel, then follow the four policies from a common initial model. Phase B takes the terminal Joint checkpoint as the starting point for two continuations under an expanded contract.

\paragraph{Data and evaluation.}
For the first stream, we extracted changed and newly added material from version 2.67.0 of the public Pydantic AI repository relative to version 1.107.0 \citep{pydanticai}. It includes source, tests, examples and documentation. We grouped files and near-duplicates before assigning data roles. Following initialization, the learner saw one 129-token window from each of 200 distinct source files in the registered order and optimized next-token cross entropy on the 128 targets per arrival. The development panel contained one window from each of 64 other files and was held out of training and online candidate selection. We averaged cross entropy over the 512 targets in each window and then over the 64 windows. Historical protection used 16 separate contexts.

\subsection{Executing a protected write}\label{exp:native-write}
Two writable modules in Transformer block 14 act on the attention output and MLP down projections. Across them, $K$ has 32 historical-channel coordinates and $B$ has 28,672 acquisition coordinates. Write $d_i,d_o$ for the input and output widths of module $\mu$. With fixed feature frames and unit acquisition scale, its weight increment $\Delta W_\mu\in\mathbb R^{d_o\times d_i}$ and the total squared cost are
\[
 \Delta W_\mu=U_\mu\Delta K_\mu V_\mu^{\mathsf T}+\Delta B_\mu\Phi_\mu,
 \qquad c=\sum_\mu\lVert\Delta W_\mu\rVert_{M_\mu}^{2}.
\]
We calibrated the metric once using all 33 token positions in each of the 16 historical contexts. Let $r_{\mu j}$ be the mean square of module input coordinate $j$ across these 528 positions and $\bar r_\mu$ its average over the $d_i$ input coordinates. Both phases use the diagonal input metric $M_{{\rm in},\mu}=\operatorname{diag}_j(\tfrac14+r_{\mu j}/(\bar r_\mu+10^{-12}))$. Stacking the columns of $\Delta W_\mu$ gives
\[
 M_\mu=M_{{\rm in},\mu}\otimes I_{d_o},\qquad
 \lVert\Delta W_\mu\rVert_{M_\mu}^{2}
 =\operatorname{tr}\!\left(\Delta W_\mu M_{{\rm in},\mu}\Delta W_\mu^{\mathsf T}\right).
\]
Applying this norm to the combined write includes the signed cross term between its two channels. The controller tests the stored FP32 update. At every position in the 16 historical sequences of 128 tokens, the native argmax must return the reference top token and the candidate must satisfy $D_{\rm KL}(p_i^{\rm ref}\Vert p_i^{\rm trial})\leq 2\times10^{-5}$ over the full vocabulary. The reference stays fixed across arrivals. The hard ceiling on $c$ is 0.0625 per write.

In a separate fixed-parent test with 32 linearized token contrasts, the Jacobians $C_K$ and $C_B$ on the vectorized historical and acquisition coefficients had dimensions $32\times32$ and $32\times28{,}672$, respectively. The measured $C_K$ was invertible, so $-C_K^{-1}C_B\Delta B$ cancels the first-order contrast response of any acquisition increment $\Delta B$. Imposing the physical budget on this joint write reduced its squared cost from 0.760318 to 0.061249994. The stored acquisition increment had Frobenius norm 0.016925774 across the two modules, or 99.9958\% of the nominal AdamW acquisition norm 0.016926482. Scalar backtracking retained 25.0\%. The stored FP32 write had joint contrast residual $\lVert C_K\Delta K+C_B\Delta B\rVert_2=1.87\times10^{-10}$. The native forward pass preserved all 2,048 reference top tokens with maximum KL $1.15549\times10^{-5}$.

The online solver tracks the current acquisition proposal and minimizes its deviation subject to the physical budget, local persistent-KL bounds and predicted token-margin bounds. For a modeled pair with parent margin $m$, the post-write lower bound is $\min(m+20\times10^{-6},30\times10^{-6})$ logit. A full native check decides acceptance. If it finds new KL violations, the solver adds them together and also uses near-threshold positions from the same scan. B-only sets $\Delta K=0$. Joint includes all 32 pure $K$ directions and weights deviations of $B$ and $K$ from their targets by $\gamma_B=1$ and $\gamma_K=0.01$, with zero nominal $K$ increment.

At current logits $z_t$, the change in reference KL under a logit increment $\delta z$ contains a drift term $(p_t-p^{\rm ref})^{\mathsf T}\delta z$ in addition to $D_{\rm KL}(p_t\Vert p_{t+1})$. A matched 32-arrival ablation tested its practical value. Under identical physical tracking and a common native guard, the Joint controller whose predictor retained the reference-drift term reduced mean development cross entropy by 0.000213316. Using an incremental KL predictor gave 0.000164192. Persistent-reference prediction yielded a larger reduction on all 64 paired documents and a 29.9\% larger mean reduction.

\subsection{Identifying compensation at a fixed parent}\label{exp:fixed-parent}
We hold the parent and the stored acquisition increment $\Delta B$ fixed and toggle $\Delta K$ to measure its contribution to the native readout. Before arrival 11, the reference token on document 10 at position 9 had a positive margin of $+2^{-20}$ logit over its competitor. We replayed an accepted Joint write from the hard-cap condition at this parent. Its acquisition increment alone gave a margin of $-2.00272\times10^{-5}$. With the same $\Delta B$, restoring $\Delta K$ raised the margin to $+4.76837\times10^{-6}$. The historical increment thus restored the reference ordering by adding $2.47955\times10^{-5}$ logit to the margin of this pair. The combined write passed the full native guard.

Starting from one Joint parent saved just before arrival 32, three controllers used the largest permitted radius, 28 modeled KL positions and the same margin rule (Table~\ref{exp:sameparent}). All passed the full guard and reduced CE on all 64 development documents.

\begin{table}[H]\centering\small
\caption{Finite writes from the same arrival-32 parent. Gain is the CE reduction on the current example or the mean reduction on 64 development documents.}\label{exp:sameparent}
\begin{tabular}{@{}lrrrr@{}}\toprule
Writable channels & Write cost $c$ & Max KL $\times10^5$ & Arrival gain $\times10^4$ & Dev. gain $\times10^5$\\\midrule
B-only &0.002871&1.95432&6.94275&6.54440\\
Joint, equal tracking weights &0.002871&1.95349&6.94275&6.54105\\
Joint, asymmetric tracking &0.003148&1.94789&6.94752&6.50007\\\bottomrule
\end{tabular}
\end{table}
Lowering the $K$ tracking weight increased the squared cost of the pure $K$ component from $2.65\times10^{-7}$ to $2.81\times10^{-4}$, with a signed cross term in the net cost. The asymmetric candidate had the lowest maximum KL and the largest gain on the arrival. B-only achieved the largest mean development gain at this parent. At a smaller radius ($k=8$), two candidates with identical stored $B$ increments passed the full guard. Restoring $\Delta K$ lowered the maximum KL from $1.99829\times10^{-5}$ to $1.99622\times10^{-5}$, adding $2.07\times10^{-8}$ of headroom to an admissible write.

\subsection{Closed-loop trajectories over 200 arrivals}\label{exp:continuous}
Four policies shared the initial model, ordered stream (seed 31162026) and 64 development documents. Both reference policies applied ordinary AdamW updates to $K$ alone. The paid reference ($\beta_{\rm gate}$) evaluated Joint candidates and discarded them before applying its reference update. The direct reference ($\beta_{\rm direct}$) applied the same update directly. At initialization and after each of the 200 updates, the two references had matching snapshot hashes covering the writable model blocks, optimizer, training cursor and random-number state. The paid reference reused 199 retained-state observations in subsequent decisions. All 200 updated reference states passed the native audit.

\begin{table}[H]\centering\small
\caption{Four policies at arrival 200. Each scan covers 2,048 positions. Additional audits include terminal checks requiring a separate pass. Parentheses count initial-radius commits for the projected policies.}\label{exp:paired200}
\begin{tabular}{@{}lrrrr@{}}\toprule
Outcome & B-only & Joint & $\beta_{\rm gate}$ & $\beta_{\rm direct}$\\\midrule
Development CE reduction & 0.003592379 & 0.003579004 & 0.000002945 & 0.000002945\\
Cumulative squared write cost & 0.244428 & 0.249510 & $8.535\times10^{-6}$ & $8.535\times10^{-6}$\\
Commits (at initial radius) & 197 (176) & 198 (179) & 200 & 200\\
Policy scans & 669 & 663 & 468 & 0\\
Additional audit scans & 1 & 1 & 1 & 200\\
Endpoint maximum KL ($\times10^5$) & 1.994775 & 1.999596 & 0.004682 & 0.004682\\\bottomrule
\end{tabular}
\end{table}

All 395 commits by the projected policies preserved the 2,048 reference tokens and met the per-position KL ceiling. Each projected policy improved held-out CE more than historical-channel-only AdamW on all 64 development documents. Averaged over these documents, the B-only and Joint reductions were about $1{,}220$ and $1{,}215$ times the reference reduction of 0.000002945. The reference used 32 writable coordinates, B-only used 28,672, and Joint used both channels.

On each projected arm the first full-radius candidate failed the native check on 198 of 200 arrivals. Among these 198, B-only committed a revised write on 174 at the original radius and on 21 after shrinking. Joint committed on 177 at the original radius and on 19 after shrinking. Pooling the arm-specific arrivals gives 391 admissions after 396 initial rejections, including 351 admissions at the original radius. The paid reference's 468 policy scans comprised 269 candidate scans and 199 observations that also served as audits.

For acquisition-only writes at a fixed parent, Theorem~\ref{e:thm:decoder} gives the largest first-order decrease in the protected unit ball as $\|P_{\ker\overline C_B}\overline H_B^*g\|$, where $\overline C_B,\overline H_B$ are the cost-whitened protection and response derivatives and $g$ is the output loss gradient. For a finite acquisition write, the replays in Section~\ref{exp:fixed-parent} measure what adding $\Delta K$ changes at fixed $\Delta B$: it restored the reference token ordering at arrival 11 and increased KL headroom at arrival 32.

Over 200 arrivals with margin and KL rows rebuilt at each parent, the acquisition channel alone reduced mean held-out CE by 0.003592379 and Joint by 0.003579004, a 0.37\% difference in favor of B-only. Terminal CE favored B-only on 37 documents and Joint on 27. Joint used the historical channel in all 198 commits, with cumulative pure-$K$ squared cost $0.005359$ (2.15\% of net squared write cost). The policies committed together on 195 arrivals. From their respective parents, Joint alone committed on arrivals 92, 99 and 100, all at the initial radius in 10 scans against B-only's 31. B-only alone committed on arrivals 150 and 179.
\subsection{Continuation across a change of contract}\label{exp:promotion-continuation}
We branched from the Joint checkpoint at arrival 200 to test the version change of Section~\ref{t:subsec:generation-promotion}. Both arms retained the 2,048 original reference distributions. For each of four previously consumed training contexts containing a verified declaration or import, we added 128 positions and fixed their full-vocabulary reference logits at the Joint checkpoint. Every accepted update had to preserve the reference top token and satisfy the same per-position KL bound at all 2,560 positions.

\begin{figure}[H]
\centering
\includegraphics[trim=4pt 4pt 4pt 4pt,clip]{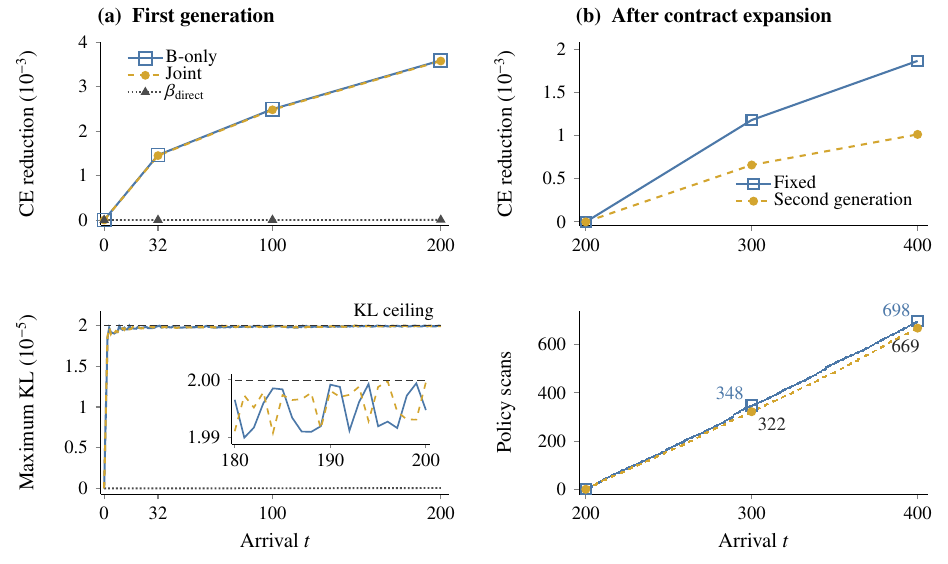}
\caption{Controlled trajectories. (a) CE reductions connect measurements at arrivals 0, 32, 100 and 200. The lower panel tracks maximum retained-state KL over 2,048 positions, carrying the parent forward at abstention. The inset enlarges arrivals 180--200. The two AdamW references followed the same parameter trajectory. (b) From the common Joint checkpoint, CE reductions use the new panel at arrivals 200, 300 and 400. The staircases count full policy scans over 2,560 positions.}
\label{fig:empirical-trajectories}
\end{figure}

The next stream used Python call sites from successive revisions of \texttt{pydantic-ai}. Static resolution identified 706 direct calls in revision $R_B$ (v2.70.0) to declarations present in $R_A$ (v2.67.0) and absent under the same qualified name in the baseline (v1.107.0). The registered stream contains 200 training windows from 40 source files and 64 development windows from 36 other files, with each call-template cluster confined to one partition. Both arms processed the same arrivals in the same order, with development evaluations at arrivals 200, 300 and 400 (Figure~\ref{fig:empirical-trajectories}(b)). Table~\ref{exp:phaseb-paired} gives cumulative outcomes at arrivals 300 and 400 against the common checkpoint. Each development window supplied 512 next-token targets. We report the mean of the 64 window losses, and separately the corresponding mean restricted to tokens spanning the callee expression.

The fixed arm continued with its first-generation coefficients and AdamW moments. The second-generation arm performed one atomic state transition. It froze the first-generation coefficients and optimizer state. In module $\mu$, the historical channel retained $U_\mu,V_\mu$. The four rows of the new acquisition frame $\Phi_{\mu,2}$ were constructed in the $M_{{\rm in},\mu}$-orthogonal complement of the columns of $[V_\mu\ \Phi_{\mu,1}^{\mathsf T}]$. This construction makes the historical and new acquisition writes orthogonal in the physical metric in exact arithmetic (Proposition~\ref{t:prop:source-injectivity}). The new coefficients and their AdamW moments started at zero. Subsequent writes took the form
\[
  \Delta W_\mu
  = U_\mu\,\Delta K_{\mu,2}\,V_\mu^{\mathsf T}
    + \Delta B_{\mu,2}\,\Phi_{\mu,2}.
\] The handoff verified identical FP32 logits at all protected positions at neutral installation, a nonzero gradient on a writable $B_2$ coefficient, and restoration of the starting state after the probe. Both arms had 28,704 writable coefficients. Renewal added 28,704 FP32 coefficients to the model while retaining the frozen first generation.

\begin{table}[htbp]\centering\small
\caption{Continuation from the same arrival-200 Joint checkpoint (development CE 2.023129042). The two column groups cover the first 100 arrivals and all 200 arrivals after that checkpoint. CE reductions and window improvements use arrival 200 as the baseline. Write costs, commits and full policy scans accumulate from arrival 201. Each scan checks all 2,560 positions.}\label{exp:phaseb-paired}
\setlength{\tabcolsep}{4pt}
\begin{tabular}{@{}lrr@{\hspace{1.6em}}rr@{}}\toprule
& \multicolumn{2}{c}{Arrivals 201--300} & \multicolumn{2}{c}{Arrivals 201--400}\\
\cmidrule(lr){2-3}\cmidrule(l){4-5}
Outcome & Fixed & Gen.\ 2 & Fixed & Gen.\ 2\\\midrule
Endpoint development CE & 2.021949435 & 2.022469912 & 2.021264235 & 2.022115670\\
Development CE reduction & 0.001179608 & 0.000659131 & 0.001864808 & 0.001013372\\
Callee-token CE reduction & 0.006064463 & 0.003949835 & 0.010230045 & 0.005144832\\
Windows improved (of 64) & 63 & 62 & 63 & 62\\
Commits (at initial radius) & 100 (86) & 100 (93) & 200 (178) & 200 (190)\\
Full policy scans & 348 & 322 & 698 & 669\\
Cumulative squared write cost & 0.105172156 & 0.108356486 & 0.201127476 & 0.190533684\\
Development CE reduction per unit cost & 0.011215968 & 0.006082982 & 0.009271770 & 0.005318598\\\bottomrule
\end{tabular}
\end{table}

By arrival 400, fixed continuation had yielded 1.74 times the development CE reduction per unit squared write cost of neutral renewal. Its terminal CE was lower on 51 of the 64 windows and on 30 of the 36 source files after averaging windows within file.

Let $Z(q)$ collect the full-vocabulary logits at the 2,560 protected positions. At arrival 200, $Z(q_F)=Z(q_N)$ and $\varsigma(q_F)\ne\varsigma(q_N)$ (Definition~\ref{m:def:frtfat}). The branches thus share an observation fiber at the fork and occupy different structural strata. At the common effective weights $W$, branch decoder $D_j$ pulls the objective differential back to $\mathrm d\ell_W\circ D_j$ and each differentiable guard row to $\mathrm d(g_i)_W\circ D_j$ (Proposition~\ref{t:prop:cross-whiten}). The resulting coefficient gradient and the inherited or fresh moments determine the nominal AdamW step by~\eqref{t:eq:adam}. On admission the controller stores the projected coefficients together with the moments and clock of that nominal step. On abstention it restores the pre-proposal model and optimizer registers and charges the branch for the scans performed.

At each arrival $t=201,\ldots,400$, both branches received the same training window and formed proposals from their own states. Each initial squared budget was the smaller of the nominal write cost and the hard cap of 0.0625. Renewal's budget was larger at arrivals 201--212 and smaller at arrivals 213--400. The median paired ratio of renewal to fixed budget over all 200 arrivals was 0.645. Let $e_{j,t}$ be the maximum native KL of the first locally projected candidate on branch $j$, minus the ceiling $2\times10^{-5}$. Despite the change in budget ordering, $0<e_{N,t}<e_{F,t}$ at all 200 arrivals. The median excess was $1.56\times10^{-6}$ for renewal and $3.22\times10^{-6}$ for fixed continuation.

Re-solving admitted all 200 writes on each branch. Renewal admitted 190 at its initial radius with 669 full policy scans, and fixed continuation admitted 178 at its initial radius with 698 scans. Every commit preserved all 2,560 reference top tokens and met the per-position KL ceiling. Among the accepted writes, renewal had the smaller maximum KL at 80 of the 200 arrivals. At arrival 400 the maximum KL was $1.990836\times10^{-5}$ for fixed continuation and $1.998259\times10^{-5}$ for renewal. Figure~\ref{fig:empirical-trajectories} follows held-out CE and verification work from the common checkpoint.
\FloatBarrier
\section{Synthesis}\label{sec:discussion}
An intervention changes both the system and the reserve available for its next action. Protocols with equal net consumption can require different initial reserves because each stage must be affordable when it begins. At one state the weighted Gram gives the affordable linear responses, and the cooperative quotient identifies directions gained by joint control. Two regular charts with the same affine data can differ in the cubic term of minimum cost. When blocks move sequentially, even a protected endpoint can require a detour whose peak cost exceeds its endpoint cost.

In the full compensated linear chart, a positive trace of the history-weighted operator $\mathcal B_{\rm hist}$ gives a direction and a constant $c>0$ such that, at every squared budget $b>0$, a reachable target lies at distance $c\sqrt b$ from the old response closure (Theorem~\ref{m:thm:historybridge}). Across successive states, a descriptor gives an exact quotient of deterministic labelled dynamics precisely when its fibers preserve readouts, enabled actions, immediate records and successor classes. For stochastic policies the two-cut identity evaluates verification at the state where it begins and the selected action at the state it leaves. Over all 200 arrivals the gated and direct AdamW references visited identical parameter checkpoints and incurred 468 and zero policy scans, respectively.

At promotion, the two Transformer controllers agreed on every protected logit and differed in the acquisition frame and optimizer moments available for the next update. Both admitted the entire subsequent stream. The inherited state delivered more held-out descent per unit squared write cost. The renewed state produced a smaller KL excess on every first native check and required fewer scans. Equality of protected readouts places the two states on one observation fiber. From each such state, the decoder and controller specify the available actions, and the reserve determines which of them the system can afford under its contract.
\appendix
\needspace{8\baselineskip}
\section{Specialization to fixed-parent affine and head-bank charts}\label{app:specializations}

\paragraph{Affine repair and compensation.}
In the finite Euclidean chart with full synchronous writes, $M=I$ and $m(u)=c(u)=\|u\|^2$. Theorems~\ref{m:thm:weighted} and~\ref{m:thm:support} specialize to fixed-parent regulation (FRT) and acquisition (FAT) \citep[Theorems~3.2--3.3 and 5.3--5.4]{paper4}.
For $Cx=a$ with $a\in\Ran C$, set $x_0=C^\dagger a$, $N=H(I-C^\dagger C)$ and $Z_0=\ker C\cap\ker H$. The additional target $Hx=d$ is attainable exactly when $d-Hx_0\in\Ran N$. Its minimum-norm repair and the residual response ellipsoid of a linear readout $R$ at budget $b\ge\|x_\star\|^2$ are
\[
 x_\star=x_0+N^\dagger(d-Hx_0),\qquad
 R x_\star+R\{w\in Z_0:\|w\|^2\le b-\|x_\star\|^2\}.
\]
Indeed, all repairs lie in $x_\star+Z_0$ with $x_\star\perp Z_0$. For $a=0$ this gives $x_\star=N^\dagger d$.

For the synchronous pair $Cx+By=0$, $Hx+\Psi y=d$, take an isometry $Q$ onto $\ker[(I-CC^\dagger)B]$ and put $T=C^\dagger BQ$, $L=\Psi Q-HT$, $W=I+T^*T$. With $G_0=NN^*$ and $G_+=G_0+LW^{-1}L^*$, the attainable targets form $\Ran G_+$ and each has minimum squared cost $\langle d,G_+^\dagger d\rangle$. If $G_0\succ0$, the Woodbury identity gives
\[
 G_0^{-1}-G_+^{-1}
 =G_0^{-1}L(W+L^*G_0^{-1}L)^{-1}L^*G_0^{-1}.
\]
Hence the extra channel lowers the minimum cost exactly when $L^*G_0^{-1}d\ne0$. For the additive readout $Z=Z_{\rm base}+\Phi A^\top$ at $A=0$, the SGD formula in Proposition~\ref{t:prop:neutral} gives $\nabla_A\ell=D^\top\Phi$ and $\Delta Z=-\eta\Phi\Phi^\top D$, where $D=\nabla_Z\ell(Z_{\rm base})$ \citep[Theorem~4.2]{paper4}.

\paragraph{Coupling of internal features and a final-head bank.}
In the final-head interface of Wang~\cite[Sections~6.1--6.5]{paper4}, let $W_{\rm head}\in\R^{V\times d}$ and fix bank rows $S\in\R^{m\times V}$. A writable head increment has transpose $X=ZS\in\R^{d\times V}$ with $Z\in\R^{d\times m}$. On fixed protected and target feature rows $C,H$, the exact logit increments are $CZS,HZS$. With columnwise vectorization their coefficient maps are $S^\top\otimes C,S^\top\otimes H$. A positive physical metric $M_V\otimes M_d$ pulls back to $(SM_VS^\top)\otimes M_d$ on coefficients and defines a norm on the quotient by the decoder kernel.
If an internal write also changes the protected features by $E$, the finite response is
\[
 (C+E)(W_{\rm head}^\top+X)-CW_{\rm head}^\top
 =CX+EW_{\rm head}^\top+EX .
\]
Exact joint repair is possible precisely when a writable $X$ solves $(C+E)X=-EW_{\rm head}^\top$.
At a single anchor whose logits a head-only write preserves, let $r_z=\nabla_z\ell(z)$ and $J_{h,\vartheta}=D_\vartheta h_\vartheta$. Holding $X$ fixed, the same cotangent $r_z$ gives the exact gradient change $\Delta\nabla_\vartheta\ell=J_{h,\vartheta}^{\top}Xr_z$. 
\needspace{5\baselineskip}

\end{document}